\documentclass[onefignum,onetabnum]{siamart190516}
\usepackage{amscd}
\usepackage{enumerate,amssymb,amsmath,graphics,epsfig}
\usepackage{tcolorbox}
\usepackage{mathrsfs}
\usepackage{epstopdf}
\usepackage{booktabs}
\usepackage{subfigure}
\theoremstyle{plain}
\newtheorem{thm}{Theorem}

\newtheorem{ass}{Assumption}
\theoremstyle{remark}

\newtheorem{remark}{Remark}

\newcommand{\m}[1]{\mathsf{#1}}
\newcommand{\RE}{\mathbb{R}}
\newcommand\Au{\m{A}_1}
\newcommand\Ad{\m{A}_2}
\newcommand\Bu{\m{B}_1}
\newcommand\Bd{\m{B}_2}
\newcommand\Aul{\m{A}^\ell_1}
\newcommand\Adl{\m{A}^\ell_2}
\newcommand\Bul{\m{B}^\ell_1}
\newcommand\Bdl{\m{B}^\ell_2}
\newcommand\KA{K_{\m{A}_1}}
\newcommand\KB{K_{\m{B}_1}}
\newcommand\NA{n_{\m{A}_1}}
\newcommand\NB{n_{\m{B}_1}}
\newcommand\Huo{H^1_0(\Omega)}

\newcommand\T{\mathcal{T}}
\newcommand\E{\mathcal{E}}
\newcommand\Pinabla{\Pi_k^\nabla}
\newcommand\Pio{\Pi^0_k}
\renewcommand\P{\mathbb{P}}
\newcommand\VemP{V_h^k(P)}
\newcommand\V{V_h^k}
\newcommand\tb{\tilde{b}_h}
\newcommand\tl{\tilde{\lambda}_h}
\newcommand\tu{\tilde{u}_h}

\DeclareMathOperator{\diag}{diag}

\headers{Parameter dependent eigenvalue problems}
{D. Boffi, F. Gardini, and L. Gastaldi}

\title{Approximation of PDE eigenvalue problems involving parameter dependent
matrices}

\author{Daniele Boffi\thanks{King Abdullah University of Science and
Technology (KAUST), Saudi Arabia, Dipartimento di Matematica ``F. Casorati'',
Universit\`a di Pavia, Italy, and Department of Mathematics and System
Analysis, Aalto University, Finland 
(\email{daniele.boffi@kaust.edu.sa}, 
\url{https://cemse.kaust.edu.sa/people/person/daniele-boffi}).}
\and
Francesca Gardini\thanks{
Dipartimento di Matematica ``F. Casorati'', Universit\`a di Pavia, Italy
(\email{francesca.gardini@unipv.it}, 
\url{http://www-dimat.unipv.it/gardini/}).}
\and
Lucia Gastaldi\thanks{DICATAM, Universit\`a di Brescia, Italy
(\email{lucia.gastaldi@unibs.it}
\url{http://lucia-gastaldi.unibs.it}).}
}%\footnotemark[3]}

\begin{document}

\maketitle

\begin{abstract}
We discuss the solution of eigenvalue problems associated with partial
differential equations that can be written in the generalized form
$\m{A}x=\lambda\m{B}x$, where the matrices $\m{A}$ and/or $\m{B}$ may depend
on a scalar parameter.
Parameter dependent matrices occur frequently when stabilized formulations are
used for the numerical approximation of partial differential equations. With
the help of classical numerical examples we show that the presence of one (or
both) parameters can produce unexpected results.
\end{abstract}

\section{Introduction}
\label{se:intro}
Several schemes for the approximation of eigenvalue problems arising from
partial differential equations lead to the algebraic form: find
$\lambda\in\RE$ and $x\in\RE^n$ with $x\ne 0$ such that
\begin{equation}
\label{eq:eig}
\m{A}x=\lambda\m{B}x,
\end{equation}
where $\m{A}$ and $\m{B}$ are matrices in $\RE^{n\times n}$.

We consider the case when the matrices $\m{A}$ and $\m{B}$ are symmetric and
positive semidefinite and may depend on a parameter. This is a typical
situation found in applications where elliptic partial differential equations
are approximated by schemes that require suitable parameters to be tuned (for
consistency and/or stability reasons). In this paper we discuss in particular
applications arising from the use of the Virtual Element Method (VEM),
see~\cite{MRR,BMRR,GV,MRV,MV,GMV,CGMMV}, where suitable parameters have to
be chosen for the correct approximation.
Similar situations are present, for instance, when a parameter-dependent
stabilization is used for the approximation of discontinuous Galerkin
formulations and when a penalty penalty term is added to the discretization of
the eigenvalue problem associated with Maxwell's
equations~\cite{CoDaMax,CoDaDurham,CoDareg,bfg,2006,WarburtonEmbree,2010,BoGue,BaCo}

In general, it may be not immediate to describe how the matrices $\m{A}$ and
$\m{B}$ depend on the given parameters. For simplicity, we
consider the case when the dependence is linear: under suitable assumptions
it is easy to discuss how the computed spectrum varies with respect to the
parameters.

The description of the spectrum in the linear case is not surprising and is
well known to a broad scientific
community~\cite{ElsnerSun,StewartSun,MR1066108,LiStewart,greenbaum2019firstorder}.
Nevertheless, the main focus of perturbation theory for eigenvalues and
eigenvectors is usually centered on the asymptotic behavior when the
parameters tend to zero. In our case, the asymptotic parameter is usually the
mesh size $h$ and we are interested in the convergence when $h$ goes to zero,
that is when the size of the involved matrices tends to infinity.
The presence of additional parameters makes the convergence more difficult to
describe and can produce unexpected results in the pre-asymptotic regime. For
this reason, we start by recalling how the spectrum of problem~\eqref{eq:eig}
is influhenced by the parameter, without considering $h$, and we translate
those results to an example of interest in Section~\ref{se:subVEM} where the
discretization parameter $h$ is considered as well.

We assume that the matrices $\m{A}$ and $\m{B}$
satisfy the following condition for $\m{C}=\m{A},\m{B}$.

\begin{ass}
\label{ass:C}
The matrix $\m{C}$ can be split into the sum
\begin{equation}
\label{eq:C}
\m{C}=\m{C}_1+\gamma\m{C}_2,
\end{equation}
where $\gamma$ is a non negative real number and $\m{C}_1$ and $\m{C}_2$ are
symmetric.
The matrices $\m{C}_1$ and $\m{C}_2$ satisfy the following properties:
\begin{enumerate}[a)]
\item $\m{C}_1$ is positive semidefinite with kernel $K_{\m{C}_1}$;
\item $\m{C}_2$ is positive semidefinite and positive definite on $K_{\m{C}_1}$;
\item $\m{C}_2$ vanishes on $K_{\m{C}_1}^\perp$, the orthogonal complement of
$K_{\m{C}_1}$ in $\RE^n$.
\end{enumerate}
\end{ass}

In~\cref{se:param} we describe the spectrum of~\eqref{eq:eig} as a
function of the parameters, in various situations that mimic the behavior of
matrices $\m{A}$ and $\m{B}$ originating from several discretization schemes.

\Cref{se:VEM}, which is the core of this paper, discusses the influence
of the parameters on the VEM approximation of eigenvalue problems. Several
numerical examples complete the papers, showing that the parameters have to be
carefully tuned and that wrong choices can produce useless results.

\section{Parametric algebraic eigenvalue problem}
\label{se:param}
Given two symmetric and positive semidefinite matrices $\m{A}$ and $\m{B}$
that can be written as
\begin{equation}
\label{eq:A}
\m{A}=\Au+\alpha\Ad
\end{equation}
and
\begin{equation}
\label{eq:B}
\m{B}=\Bu+\beta\Bd,
\end{equation}
with nonnegative parameters $\alpha$ and $\beta$, we consider the
eigensolutions to the generalized problem~\eqref{eq:eig}.

We assume that the splitting of the matrices $\m{A}$ and $\m{B}$ is obtained
with symmetric matrices and satisfies~\cref{ass:C} for
$\m{C}_1=\Au,\Bu$ and $\m{C}_2=\Ad,\Bd$. Moreover we denote by $\NA$ and $\NB$
the dimension of $\KA$ and $\KB$, respectively.

\begin{remark}
\label{re:geneig}
Problem~\eqref{eq:eig} has $n$ eigenvalues if and only if
$\mathrm{rank}\m{B}=n$, see~\cite{GVL}. If $\m{B}$ is singular the spectrum can
be finite, empty, or infinite (if $\m{A}$ is singular too). If $\m{A}$ is non
singular, usually one can
circumvent this difficulty by computing the eigenvalues of $\m{B} x=\mu\m{A}x$
and setting $\lambda=1/\mu$. The kernel of $\m{B}$ is the eigenspace associated
with the vanishing eigenvalue with multiplicity $m$, and 
the original problem has exactly $m$ eigenvalues conventionally set to
$\infty$.
\end{remark}

We want to study the behavior of the eigenvalues as the parameters $\alpha$ and
$\beta$ vary. We consider three cases.
\subsection{Case 1}
\label{se:caso1}
We fix $\beta>0$ so that $\m{B}$ is positive definite. This implies
that the eigenvalues of~\eqref{eq:eig} are all non negative.
Let us consider first $\alpha=0$ so that~\eqref{eq:eig} reduces to
\begin{equation}
\label{eq:eigA1}
\Au x=\lambda \m{B}x.
\end{equation}
Since $\Au$ is positive semidefinite, $\lambda=0$ is an eigenvalue
of~\eqref{eq:eigA1} with multiplicity equal to $\NA=\dim(\KA)$ and $\KA$ is the
associated eigenspace. 
In addition, we
have $m_A=n-\NA$ positive eigenvalues $\{\mu_1\le\dots\le\mu_{m_A}\}$ counted
with their multiplicity (since we are dealing with a symmetric problem, we do
not distinguish between geometric and algebraic multiplicity). We denote by
$v_j\in\KA^\perp$ the eigenvector associated with $\mu_j$, that is
\[
\Au v_j=\mu_j \m{B}v_j.
\]
Thanks to property c) of~\cref{ass:C} when $\m{C}=\m{A}$, we observe
that
\[
\m{A}v_j=\Au v_j+\alpha\Ad v_j=\Au v_j=\mu_j\m{B} v_j.
\]
Therefore $(\mu_j,v_j)$, for $j=1,\dots,m_A$, are eigensolutions of the
original system~\eqref{eq:eig}. 

On the other hand, the eigensolutions of
\[
\Ad w=\nu\m{B}w
\]
are characterized by the fact that $\NA$ eigenvalues $\nu_i$ ($i=1,\dots,\NA$)
are strictly positive with corresponding eigenvectors $w_i$ belonging to
$\KA$, while the remaining $m_A$ eigenvalues vanish and have $\KA^\perp$ as
eigenspace. 
Thus, property a) of~\cref{ass:C}, for $\m{C}=\m{A}$, yields
\[
\m{A}w_i=\Au w_i+\alpha\Ad w_i=\alpha\Ad w_i=\alpha\nu_i \m{B}w_i,
\]
which means that $(\alpha\nu_i,w_i)$, for $i=1,\dots,\NA$, are eigensolutions
of~\eqref{eq:eig}.

\begin{figure}
\begin{center}
\includegraphics[width=.48\textwidth]{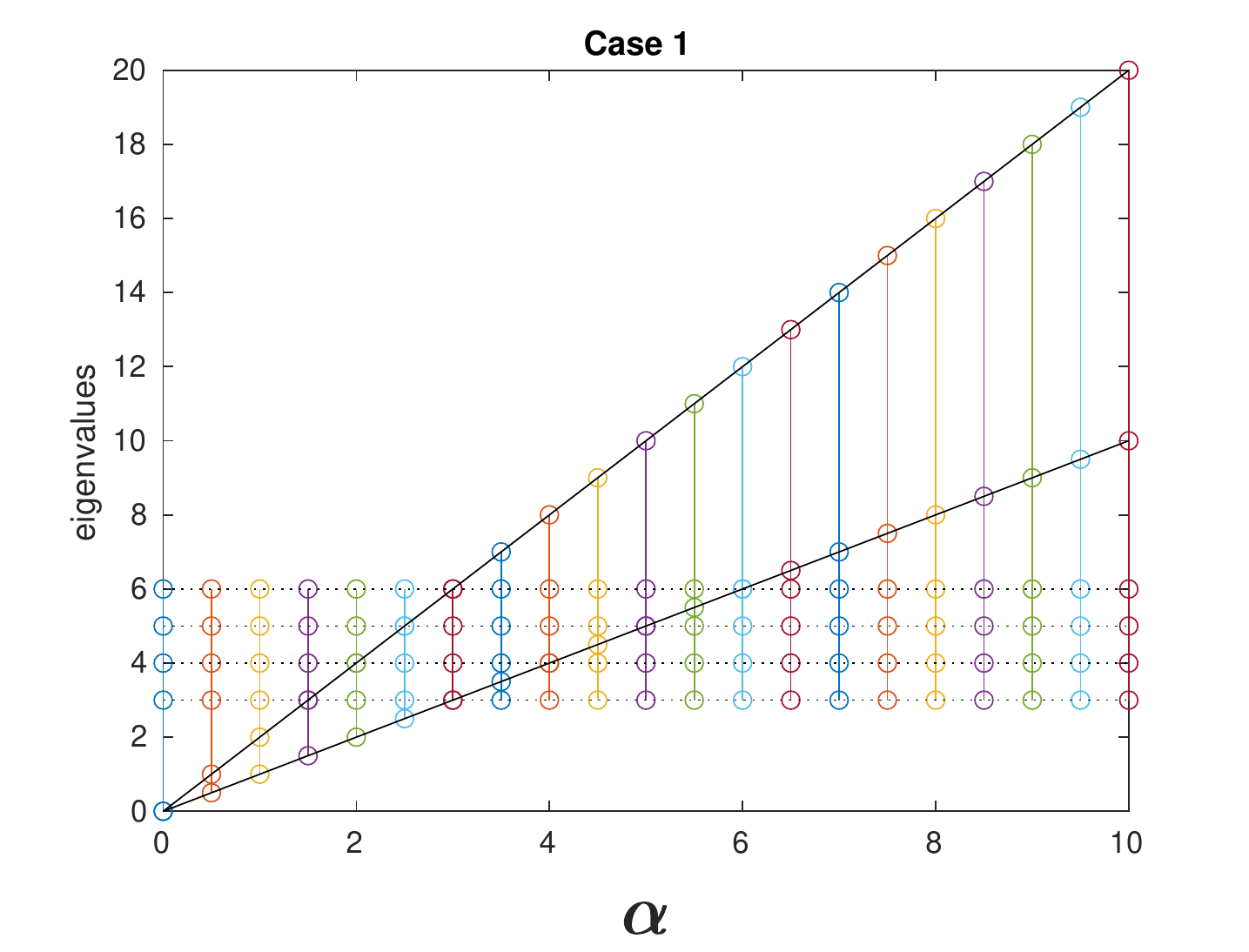}
\includegraphics[width=.48\textwidth]{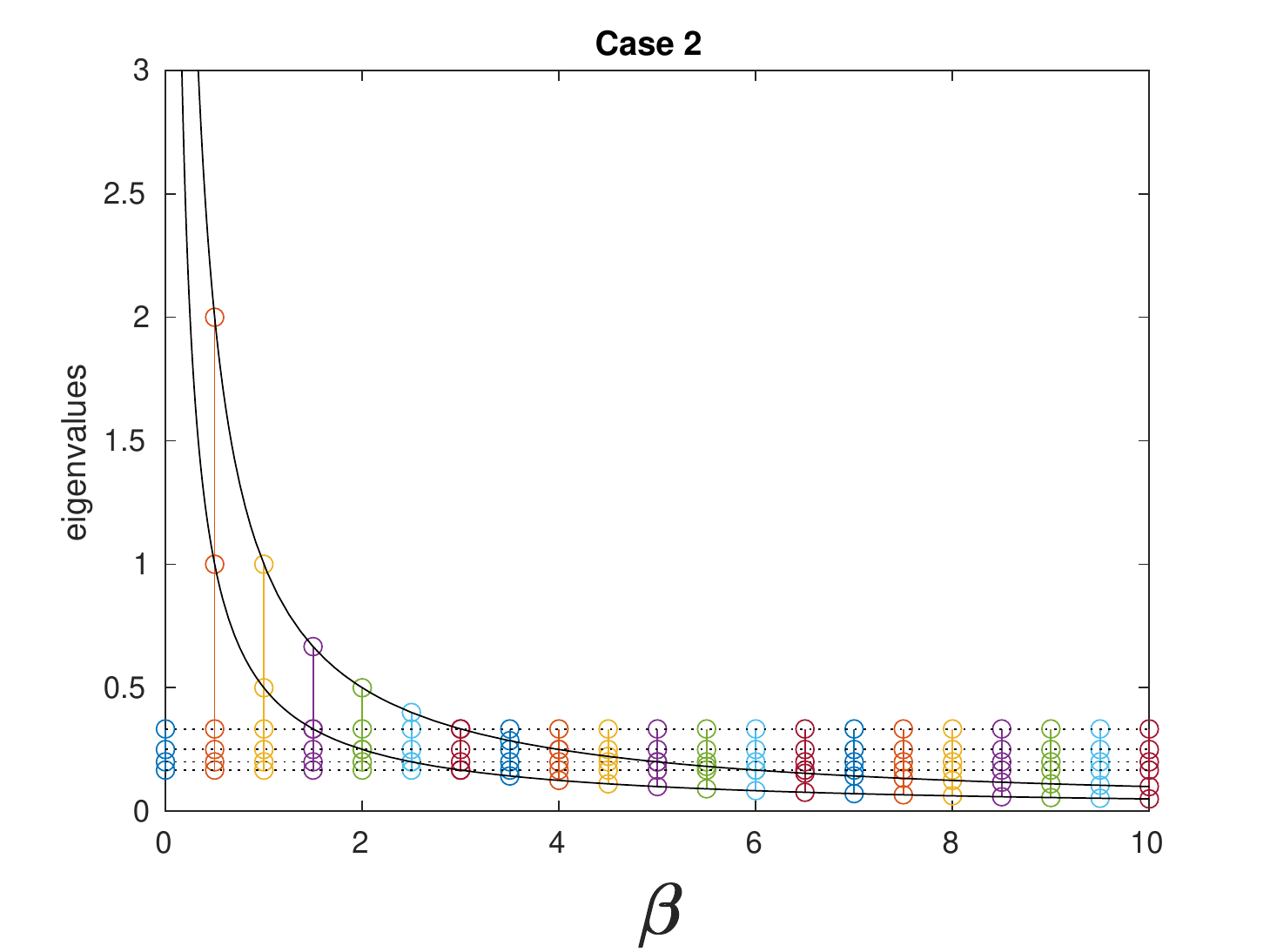}
\caption{Dependence of the eigenvalues on the parameters $\alpha$ (Case~1) and
$\beta$ (Case~2), respectively}
\label{fig:case1}
\end{center}
\end{figure}

Summarizing the eigenvalues of~\eqref{eq:eig} are:
\begin{equation}
\label{eq:eigA}
\lambda_k=
\left\{
\begin{array}{ll}
\alpha\nu_k & \text{if }1\le k\le\NA\\
\mu_{k-\NA} & \text{if } \NA+1\le k\le n.
\end{array}
\right.
\end{equation}
The left panel in~\cref{fig:case1} shows the eigenvalues
of a simple example where $\m{A}\in\RE^{6\times6}$ is
obtained by the combination of diagonal matrices with entries
\begin{equation}
\label{eq:matrici}
\diag(\Au)=[3,4,5,6,0,0],\quad\diag(\Ad)=[0,0,0,0,1,2].
\end{equation}
and $\m{B}=\mathbb{I}_6$ is the identity matrix.
 
Along the vertical lines we see the eigenvalues corresponding to a fixed value
of $\alpha$. The eigenvalues $3,4,5,6$ are associated with eigenvectors in
$\KA^\perp$ and do not depend on $\alpha$. The solid lines starting at the origin
display the eigenvalues $1,2$ multiplied by $\alpha$.  

\begin{remark}
\label{re:intersection}
We observe that if $\Ad$ is not positive definite on $\KA$, its kernel has a
nonempty intersection with $\KA$. Let $n_{12}$ be the dimension of this
intersection, then problem~\eqref{eq:eig} admits $n_{12}$ vanishing eigenvalues
which appear in the first case of~\eqref{eq:eigA}. 
\end{remark}

\subsection{Case 2}
\label{se:caso2}
Let us now fix $\alpha>0$, so that $\m{A}$ is positive definite. We have
that all the eigenvalues are positive. We observe that when $\beta=0$, the
matrix $\m{B}=\Bu$ may be singular, therefore it is convenient to consider the
following problem:
\begin{equation}
\label{eq:eig2}
\m{B}x=\chi\m{A}x,
\end{equation}
where $\chi=\frac 1\lambda$. If $\chi=0$, we conventionally set $\lambda=\infty$.
Problem~\eqref{eq:eig2} reproduces the same situation we had in Case~1, with the
matrices $\m{A}$ and $\m{B}$ switched.
Repeating the same arguments as before, we obtain that problem~\eqref{eq:eig2}
has two families of eigenvalues
\[
\chi_k=
\left\{
\begin{array}{ll}
\beta\xi_k & \text{if }1\le k\le\NB\\
\zeta_{k-\NB} & \text{if } \NB+1\le k\le n,
\end{array}
\right.
\]
where
\[
\aligned
&\Bu r_j=\zeta_j \m{A} r_j,\ j=1,\dots,n-\NB &&\text{ with }r_j\in\KB^\perp\\
&\Bd s_i=\xi_i\m{A}s_i,\ i=1,\dots,\NB &&\text{ with }s_i\in\KB.
\endaligned
\]
Going back to the original problem~\eqref{eq:eig}, we can conclude
that the eigensolutions of~\eqref{eq:eig} are the following ones:
\begin{equation}
\label{eq:eigB} 
\aligned
&\left(\frac1{\beta\xi_k},s_k\right) &&\text{ for } k=1,\dots,\NB\\
&\left(\frac1{\zeta_{k-\NB}},r_{k-\NB}\right)&&\text{ for }k=\NB+1,\dots,n.
\endaligned
\end{equation}
In the right panel of~\cref{fig:case1}, we report the eigenvalues
of a simple example where $\m{A}=\mathbb{I}_6$ and $\m{B}$ is obtained by
combining $\Bu=\Au$ and $\Bd=\Ad$ defined in~\eqref{eq:matrici}. We see that
the eigenvalues $\frac13,\frac14,\frac15,\frac16$ are independent of $\beta$
and that the remaining two eigenvalues lie along the hyperbolas $\frac1\beta$
and $\frac1{2\beta}$, plotted with solid line.
\subsection{Case 3}
\label{se:caso3}
We consider now the case when $\alpha$ and $\beta$ can vary independently from
each other. We have different situations corresponding to the relation between
$\KA$ and $\KB$. 
To ease the reading, let us introduce the following notation:
\begin{subequations}
\begin{alignat}{1}
&\Au v=\mu\Bu v
\label{eq:notation11}\\
&\Au w=\nu \Bd w
\label{eq:notation12}\\
&\Ad y=\chi\Bu y 
\label{eq:notation21}\\
&\Ad z=\eta\Bd z.
\label{eq:notation22}
\end{alignat}
\end{subequations} 
In this case the space
$\RE^n$ can be decomposed into four mutually orthogonal subspaces
\[
\RE^n=(\KA\cap\KB)\oplus(\KA\cap\KB^\perp)\oplus(\KA^\perp\cap\KB)\oplus
(\KA^\perp\cap\KB^\perp).
\]
Let us denote by $n_{\Au\cap\Bu}$ the dimension of $\KA\cap\KB$.
If $\KA\cap\KB\ne\emptyset$,
for $x\in\KA\cap\KB$ the eigenproblem to be solved is
$\alpha\Ad x=\lambda\beta\Bd x$, hence the eigenvalues are given by
$\frac{\alpha}{\beta}\eta_i$ $i=1,\dots,n_{\Au\cap\Bu}$,
see~\eqref{eq:notation22}. Next, if $x\in\KA\cap\KB^\perp$ we have to solve
$\alpha\Ad x=\lambda\Bu x$, which admits $(\alpha\chi_i,y_i)$
$i=1,\dots,\NA-n_{\Au\cap\Bu}$ as eigensolutions where $(\chi_i,y_i)$ are
defined in~\eqref{eq:notation21}. Similarly, if $x\in\KA^\perp\cap\KB$, we find
that the eigensolutions are $\left(\frac1{\beta}\nu_i,w,_i\right)$ 
$i=1,\dots,\NB-n_{\Au\cap\Bu}$ with $(\nu_i,w_i)$ given
by~\eqref{eq:notation12}. In the last case, $x\in\KA^\perp\cap\KB^\perp$, the
matrices $\m{A}$ and $\m{B}$ are non singular and thanks to property c) in
~\cref{ass:C}, for $\m{C}=\m{A}$ and $\m{C}=\m{B}$, we obtain that
the eigenvalues are positive and independent of $\alpha$ and $\beta$ and
correspond to those of~\eqref{eq:notation11}. In conclusion, we have
\[
\lambda_k=\left\{
\begin{array}{ll}
\displaystyle \frac{\alpha}{\beta}\eta_k
&\quad\text{ if }1\le k\le n_{\Au\cap\Bu}\\
\alpha\chi_{k-n_{\Au\cap\Bu}} &\quad\text{ if }n_{\Au\cap\Bu}+1\le k\le\NA\\
\displaystyle\frac1{\beta}\nu_{k-\NA}
&\quad\text{ if }\NA+1\le k\le\NA+\NB-n_{\Au\cap\Bu}\\
\mu_{k-\NA+\NB-n_{\Au\cap\Bu}}
&\quad\text{ if }\NA+\NB-n_{\Au\cap\Bu}+1\le k\le n.
\end{array}
\right.
\]

\begin{figure}[ht]
\begin{center}
\includegraphics[width=.98\textwidth]{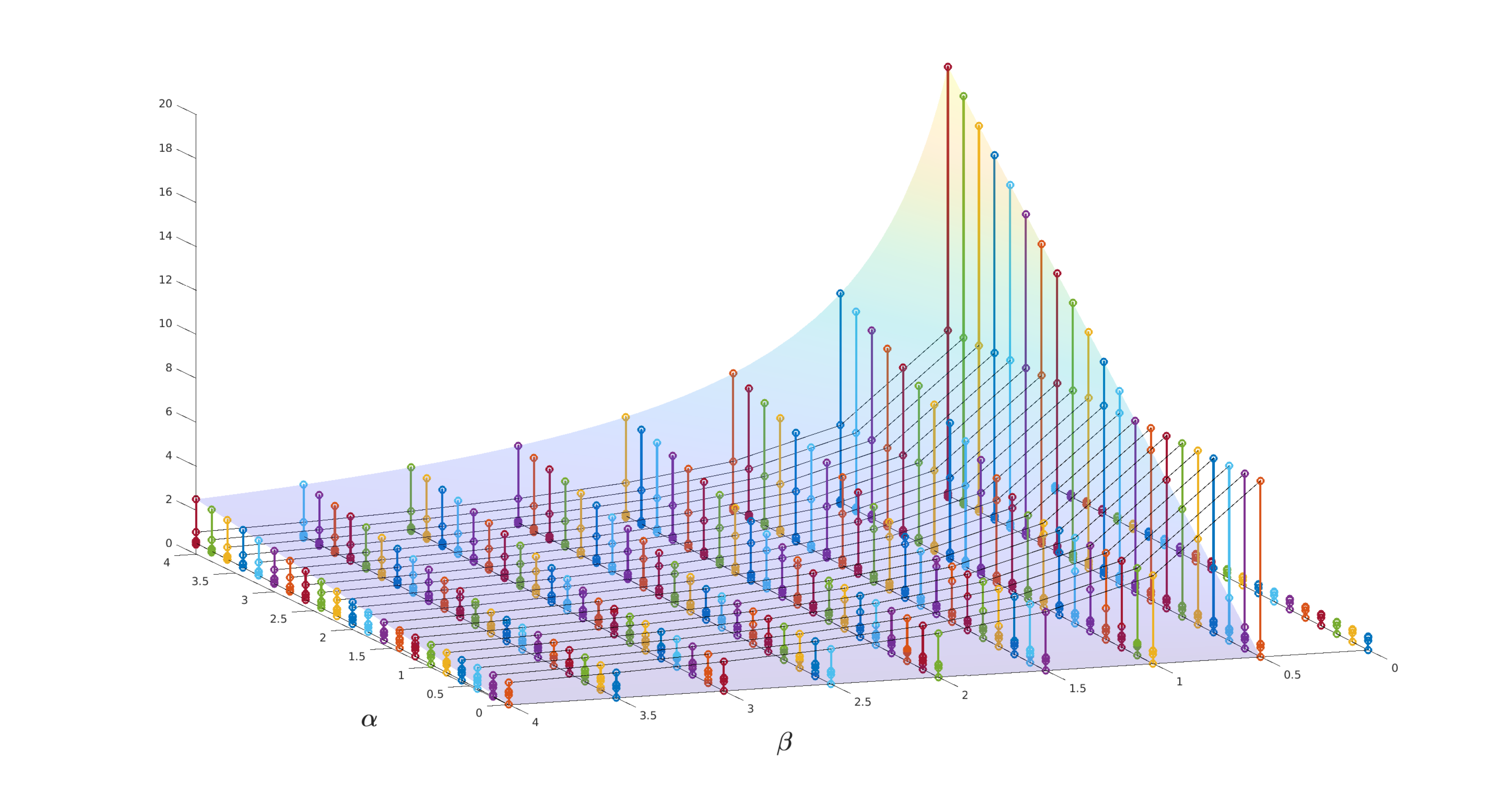}
\caption{Eigenvalues when $\KA\cap\KB\ne\emptyset$ as a function of $\alpha$
and $\beta$}
\label{fig:ab7}
\end{center}
\end{figure}
We report in~\cref{fig:ab7} the eigenvalues illustrating this last case when
we have diagonal matrices given by
\[
\aligned
&\diag(\Au)=[3,0,0,4,5,6] &&\quad\diag(\Ad)=[0,1,2,0,0,0]\\
&\diag(\Bu)=[7,8,0,0,9,10] &&\quad\diag(\Bd)=[0,0,0.8,1,0,0].
\endaligned
\]
The surface contains the
eigenvalues depending on both $\alpha$ and $\beta$, the hyperbolas those
depending only on $\beta$ and the straight lines those depending only on
$\alpha$. If we cut the three dimensional picture with a plane at $\beta>0$
fixed we recognize the behavior analyzed in~\cref{se:caso1} and shown
in~\cref{fig:case1} left. Analogously, taking a plane with $\alpha>0$
fixed, we recover Case 2 (see~\cref{se:caso2}).

If $\KA\cap\KB=\emptyset$, we set $n_{\Au\cap\Bu}=0$, hence the eigenvalues
are
\[
\lambda_k=\left\{
\begin{array}{ll}
\alpha\chi_k & \text{ if }1\le k\le \NA\\
\displaystyle\frac1\beta{\nu_{k-\NA}} & \text{ if } \NA+1\le k\le\NA+\NB\\
\mu_{k-\NA-\NB} & \text{ if } \NA+\NB+1\le k\le n
\end{array}.
\right.
\]

\begin{figure}
\begin{center}
\includegraphics[width=.98\textwidth]{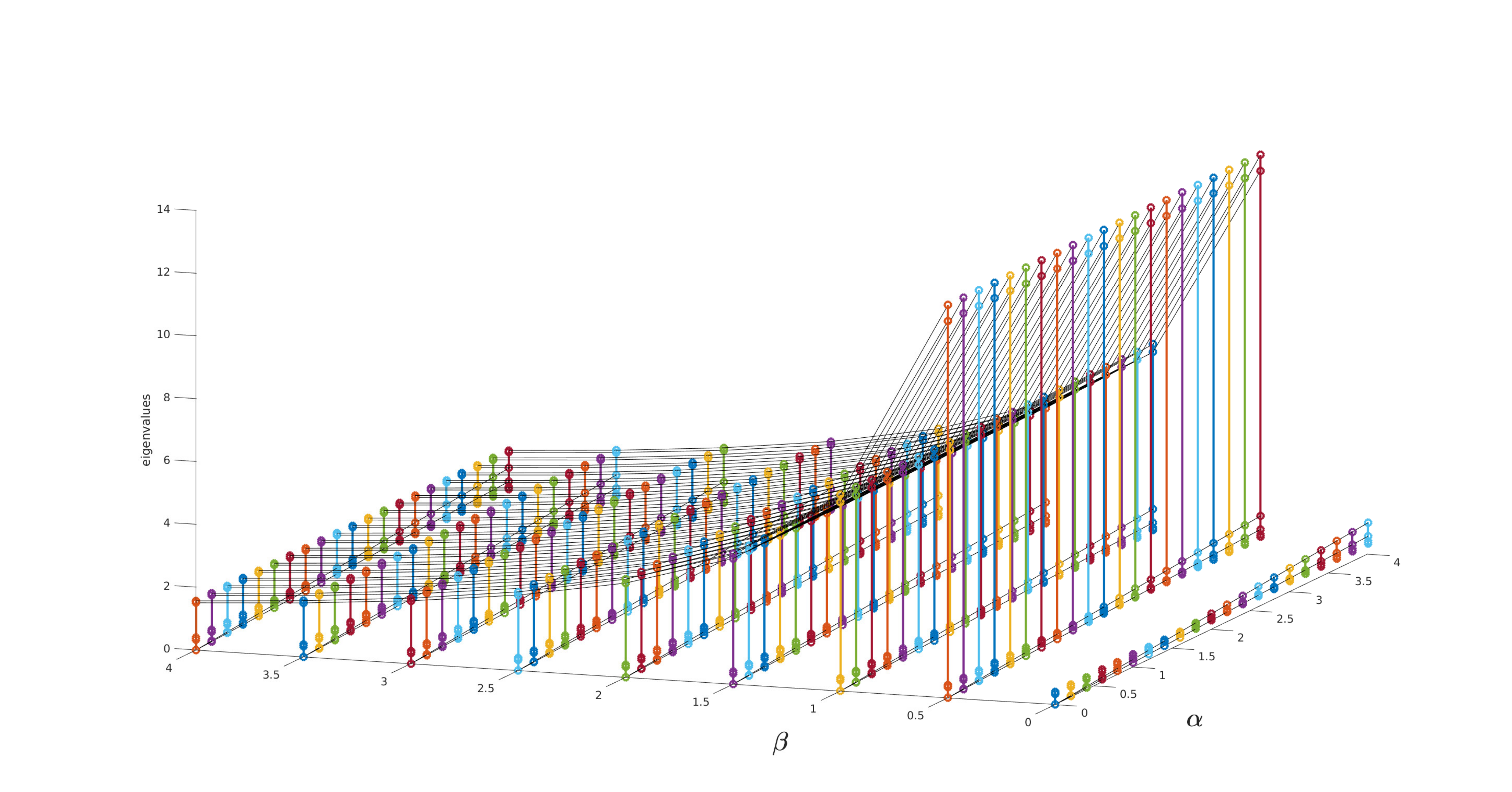}
\caption{Eigenvalues when $\KA\cap\KB=\emptyset$ as a function of $\alpha$ and
$\beta$}
\label{fig:ab4}
\end{center}
\end{figure}

In order to illustrate the case $\KA\cap\KB=\emptyset$, we report in~\cref{fig:ab4} 
the eigenvalues computed using the following diagonal matrices with entries

\[
\aligned
&\diag(\Au)=[0,0,3,4,5,6],&&\quad\diag(\Ad)=[1,2,0,0,0,0],\\
&\diag(\Bu)=[7,8,9,10,0,0],&&\quad\diag(\Bd)=[0,0,0,0,0.8,1].
\endaligned
\]
For a fixed $\alpha$, we can see in solid line the hyperbolas
$\frac{\nu_j}\beta$, $j=1,2$ while when $\beta$ is fixed we can see the
straight lines $\alpha\chi_j$, $j=1,2$. The remaining two eigenvalues are
independent of $\alpha$ and $\beta$.

\section{Virtual element method for eigenvalue problems}
\label{se:VEM}
In this section we recall how algebraic eigenvalue problems similar to the
ones discussed in the previous section can be obtained withing the framework
of the Virtual Element Method (VEM) for the discretization of elliptic
eigenvalue problems, see~\cite{GV,GMV}. 

We consider the model problem of the Laplacian operator.
Given a connected open domain with Lipschitz continuous boundary
$\Omega\subseteq\RE^d$, with $d=2,3$, we look for eigenvalues
$\lambda\in\RE$ and eigenfunctions $u\ne0$ such that
\[
\left\{
\begin{array}{ll}
-\Delta u=\lambda u &\quad\text{ in }\Omega\\
u=0&\quad\text{ on }\partial\Omega.
\end{array}
\right.
\]
In view of the application of VEM, we consider the weak form: find
$\lambda\in\RE$ and $u\in\Huo$ with $u\ne0$ such that
\begin{equation}
\label{eq:Laplace}
a(u,v)=\lambda b(u,v)\quad\forall v\in\Huo,
\end{equation}
where
\[
a(u,v)=(\nabla u,\nabla v),\quad b(u,v)=(u,v),
\]
and $(\cdot,\cdot)$ is the scalar product in $L^2(\Omega)$.

It is well-known that problem~\eqref{eq:Laplace} admits an infinite sequence
of positive eigenvalues
\[
0<\lambda_1\le \dots\le\lambda_i\le\cdots
\]
repeated according to their multiplicity, each one associated with an
eiegenfunction $u_i$ with the following properties
\begin{equation}
\label{eq:eigf}
\aligned
& a(u_i,u_j)=b(u_i,u_j)=0\quad \text{if }i\ne j\\
& b(u_i,u_i)=1,\quad a(u_i,u_i)=\lambda_i.
\endaligned
\end{equation}
Let us briefly recall the definition of the virtual element spaces and of the
discrete bilinear forms which we are going to use in this section,
see~\cite{BBCMMR,AABMR}. 
We present only the two dimensional spaces, the three dimensional ones are
obtained using the 2D virtual elements on the faces.

We decompose $\Omega$ into polygons $P$, with diameter $h_P$ and area $|P|$.
Similarly, if $e$ is an edge of an element $P$, we denote by $h_e=|e|$ its
length. Depending on the context $\partial P$ refers to either the boundary of
$P$ or the set of the edges of $P$.
The notation $\T_h$ and $\E_h$ stands for the set of the elements
and the edges, respectively. As usual, $h=\max_{P\in\T_h} h_P$.
We assume the following mesh regularity condition (see~\cite{BBCMMR}):
there exists a positive constant $\gamma$, independent of $h$, such that each
element $P\in\T_h$ is star-shaped with respect to a ball of radius greater
than $\gamma h_P$; moreover, for every element $P$ and for every edge
$e\subset\partial P$, it holds $h_e \ge \gamma h_P$.

For $k\ge1$ and $P\in\T_h$ we define 
\[
\tilde{V}_h^k(P)=\{v\in H^1(P): v|_{\partial P}\in C^0(\partial P),
v|_e\in\P_k(e)\ \forall e\subset\partial P, \Delta v\in\P_k(P)\}. 
\]
We consider the following linear forms on the space $\tilde{V}_h^k(P)$
\begin{enumerate}
\item[D1]: the values $v(V_i)$ at the vertices $V_i$ of $P$, 
\item[D2]: the scaled edge moments up to order $k-2$
\[
\dfrac{1}{|e|}\int_e vm\,\text{d}s\quad\forall m\in\mathcal{M}_{k-2}(e),\ 
\forall e\subset\partial P,
\]
\item[D3]: the scaled element moments up to order $k-2$
\[
\dfrac{1}{|P|}\int_P vm\,\text{d}x\quad\forall m\in\mathcal{M}_{k-2}(P),\ 
\]
\end{enumerate}
where $\mathcal{M}_{k-2}(\omega)$ is the set of scaled monomials on $\omega$, 
namely
\[
\mathcal{M}_{k-2}(\omega)=\Big\{\Big(\dfrac{\mathbf{x}-\mathbf{x}_\omega}
{h_\omega}\Big)^s, |s|\le k-2\Big\},
\]
with $\mathbf{x}_\omega$ the barycenter of $\omega$, and with the convention 
that $\mathcal{M}_{-1}(\omega)=\emptyset$.

From the values of the linear operators D1--D3, on each element $P$ we can
compute a projection operator $\Pinabla: \tilde{V}_h^k(P)\rightarrow
\P_k(P)$ defined as the unique solution of the following problem:
\begin{equation}
\label{eq:pinabla}
\aligned
& a^P(\Pinabla v-v,p)=0\quad\forall p\in\P_k(P)\\
& \int_{\partial P}(\Pinabla v-v)\text{d}s=0,
\endaligned
\end{equation}
where $a^P(u,v)=(\nabla u,\nabla v)_P$ and $(\cdot,\cdot)_P$ denotes the
$L^2(P)$-scalar product.

The local virtual space is defined as
\begin{equation}
\label{eq:VemP}
\VemP=\left\{v\in\tilde{V}_h^k(P):\int_P (v-\Pinabla v) p\text{d}x=0\ 
\forall p\in(\P_k\setminus\P_{k-2})(P)\right\},
\end{equation}
where $(\P_k\setminus\P_{k-2})(P)$ contains the polynomials in $\P_k(P)$
$L^2$-orthogonal to $\P_{k-2}(P)$.

We recall that by construction $\P_k(P)\subset\VemP$, so that the optimal rate
of convergence is ensured. Moreover, the linear operators D1--D3 provide a 
unisolvent set of degrees of freedom (DoFs) for $\VemP$, which allows us to
define and compute $\Pinabla$ on $\VemP$. In addition, the
$L^2$-projection operator $\Pio:\VemP\to\P_k(P)$ is also computable using the
DoFs.

The global virtual space is
\begin{equation}
\label{eq:Vem}
\V=\{v\in\Huo: v|_P\in\VemP\ \forall P\in\T_h\}.
\end{equation}
In order to discretize problem~\eqref{eq:Laplace}, we introduce the
discrete counterparts $a_h$ and $b_h$ of the bilinear forms $a$ and $b$,
respectively. Both discrete forms are obtained as sum of the following
local contributions: for all $u_h,v_h\in\V$
\begin{equation}
\label{eq:abP}
\aligned
&a_h^P(u_h,v_h)=a^P(\Pinabla u_h,\Pinabla v_h)
+S_a^P((I-\Pinabla)u_h,(I-\Pinabla)v_h)\\
&b_h^P(u_h,v_h)=b^P(\Pio u_h,\Pio v_h)+S_b^P((I-\Pio)u_h,(I-\Pio)v_h),
\endaligned
\end{equation}
where $b^P(u,v)=(u,v)_P$, and $S_a^P$ and $S_b^P$ are symmetric
positive definite bilinear forms on $\VemP\times\VemP$ such that 
\begin{equation}
\label{eq:stab}
\aligned
&c_0 a^P(v,v)\le S_a^P(v,v)\le c_1a^P(v,v) &&\quad\forall v\in\VemP
\text{ with }\Pinabla v=0\\
&c_2 b^P(v,v)\le S_b^P(v,v)\le c_3b^P(v,v) &&\quad\forall v\in\VemP
\text{ with }\Pio v=0,
\endaligned
\end{equation}
for some positive constants $c_i$ ($i=0,\dots,3$) independent of $h$.
We define
$a_h(u_h,v_h)=\sum_{P\in\T_h}a_h^P(u_h,v_h)$ and 
$b_h(u_h,v_h)=\sum_{P\in\T_h}b_h^P(u_h,v_h)$.

The virtual element counterpart of~\eqref{eq:Laplace} reads: find $\lambda_h$
and $u_h\in\V$ with $u_h\ne0$ such that
\begin{equation}
\label{eq:LaplV}
a_h(u_h,v_h)=\lambda_h b_h(u_h,v_h)\quad\forall v_h\in\V.
\end{equation}
Thanks to~\eqref{eq:stab}, the discrete problem~\eqref{eq:LaplV} admits
$N_h=\dim{\V}$ positive eigenvalues
\[
0<\lambda_{1h}\le\dots\lambda_{N_hh}
\]
and the corresponding eigenfunctions $u_{ih}$, for $i=1,\dots,N_h$, enjoy the
discrete counterpart of properties in~\eqref{eq:eigf}.

The following convergence result has been proved in~\cite{GV}.
\begin{thm}
\label{th:conv}
Let $\lambda$ be an eigenvalue of~\eqref{eq:Laplace} of multiplicity $m$ and
$\mathcal{E}_\lambda$ the corresponding eigenspace. Then there are exactly $m$
discrete eigenvalues of~\eqref{eq:LaplV} $\lambda_{j(i)h}$ ($i=1,\dots,m$)
tending to $\lambda$.
Moreover, assuming that $u\in H^{1+r}(\Omega)$, for all
$u\in\mathcal{E}_\lambda$, the following inequalities hold true:
\[
\aligned
&|\lambda-\lambda_{j(i)h}|\le Ch^{2t}\\
&\hat\delta(\mathcal{E}_\lambda,\oplus_i\mathcal{E}_{j(i)h})\le Ch^t,
\endaligned
\]
where $t=\min(k,r)$, $\hat\delta(\mathcal{E},\mathcal{F})$ represents
the gap between the spaces $\mathcal{E}$ and $\mathcal{F}$, and
$\mathcal{E}_{\ell h}$ is the eigenspace spanned by $u_{\ell h}$.
\end{thm}
\begin{remark}
\label{re:nonstab}
It is also possible to consider on the right hand side of~\eqref{eq:LaplV}
the bilinear form for  $\tb(u_h,v_h)=\sum_{P\in\T_h}b^P(\Pio u_h,\Pio v_h)$.
This leads to the following discrete eigenvalue problem: find
$(\tl,\tu)\in\RE\times\V$ with $\tu\ne0$ such that
\begin{equation}
\label{eq:nonstab}
a_h(\tu,v_h)=\tl\tb(\tu,v_h)\quad\forall v_h\in\V.
\end{equation}
The analogue of~\cref{th:conv} holds true for this partially non
stabilized discretization as well.
\end{remark}

\subsection{Computational aspects and numerical results}
\label{se:subVEM}
In order to compute the solution of problems~\eqref{eq:LaplV}
and~\eqref{eq:nonstab}, we need to describe how to obtain the matrices
associated to our bilinear forms. By construction the matrix $\Au$
(respectively, $\Bu$) associated with $\sum_P a^P(\Pinabla\cdot,\Pinabla\cdot)$
(respectively, $\sum_P b^P(\Pio\cdot,\Pio\cdot)$) has kernel corresponding to
the elements $v_h\in\V$ such that
$\Pinabla v_h$ is constant (respectively, $\Pio v_h=0$) for all $P\in\T_h$.

We observe that the local contributions of the bilinear forms displayed
in~\eqref{eq:abP} mimic the following exact relations
\begin{equation}
\label{eq:exactab}
\aligned
&a^P(u_h,v_h)=
a^P(\Pinabla u_h,\Pinabla v_h)+a^P((I-\Pinabla)u_h,(I-\Pinabla)v_h)\\
&b^P(u_h,v_h)=
b^P(\Pio u_h,\Pio v_h)+b^P((I-\Pio)u_h,(I-\Pio)v_h).
\endaligned
\end{equation}
Let us denote by $\Aul$, $\Adl$, $\Bul$ and $\Bdl$ the matrices whose entries are
given by
\begin{equation}
\label{eq:matr}
\aligned
&(\Aul)_{ij}=a^P(\Pinabla \phi_i,\Pinabla \phi_j),&\quad
&(\Adl)_{ij}=a^P((I-\Pinabla)\phi_i,(I-\Pinabla)\phi_j)\\
&(\Bul)_{ij}=b^P(\Pio \phi_i,\Pio \phi_j),&\quad
&(\Bdl)_{ij}=b^P((I-\Pio)\phi_i,(I-\Pio)\phi_j)
\endaligned
\end{equation}
with $\phi_i$ basis functions for $\VemP$.

Even if the global matrices $\m{A}$ and $\m{B}$ do not satisfy the properties
stated in~\cref{ass:C}, it turns out that~\cref{ass:C} is
fulfilled by $\m{C}=\Bul+\beta\Bdl$; moreover, $\m{C}=\Aul+\alpha\Adl$ is
characterized by the situation described in~\cref{re:intersection}.

We start with the pair $\Aul$ and $\Adl$.
The kernel $K_{\Aul}$, with abuse of notation, is characterized by
\[
K_{\Aul}=\{v\in\VemP:  a^P(\Pinabla v,\Pinabla w)=0\ \forall w\in\VemP\},
\]
that is, $K_{\Aul}$ is made of $v$ with constant $\Pinabla v$ on $P$.
Moreover, the orthogonal complement of $K_{\Aul}$, denoted by
$K_{\Aul}^\perp$ contains the elements $v\in\VemP$ such that $a^P(v,w)=0$ for
all $w\in K_{\Aul}$.

We now show that $\Adl( K_{\Aul}^\perp)=0$, that is, for all $v\in K_{\Aul}^\perp$,
$a^P((I-\Pinabla)v,(I-\Pinabla)w)=0$ for all $w\in\VemP$. 
We recall that, if $v\in K_{\Aul}^\perp$, then $a^P(v,w)=0$ for 
all $w\in K_{\Aul}$. This implies that for $v\in K_{\Aul}^\perp$ and
$w\in K_{\Aul}$, it holds true that
$a^P(v,w)=a^P((I-\Pinabla)v,(I-\Pinabla)w)=0$.
Now we can write for all $w\in\VemP$
\[
\aligned
&a^P((I-\Pinabla)v,(I-\Pinabla)w)\\
&\quad=a^P((I-\Pinabla)v,(I-\Pinabla)(I-\Pinabla)w)+
a^P((I-\Pinabla)v,(I-\Pinabla)\Pinabla w)=0.
\endaligned
\]
Indeed, $\Pinabla(I-\Pinabla)w=0$ implies that $(I-\Pinabla)w\in K_{\Aul}$, and
thus the first term vanishes,
while for the second term it is enough to observe that
$\Pinabla(\Pinabla w)=\Pinabla w$.
Thus property c) of~\cref{ass:C} is verified for $\m{C}=\m{A}$.

Concerning property b) of~\cref{ass:C}, we have 
by construction, that $a^P((I-\Pinabla)v,(I-\Pinabla)v)\ge0$
for all $v\in\VemP$, see~\eqref{eq:exactab}. On the other hand, if $v$ is
constant on $P$, then $\Pinabla v=v$ is constant, therefore
$v\in K_{\Aul}$ and $(I-\Pinabla)v=0$ so that $v$ belongs also to the
kernel of $\Adl$. Hence the pair $\Aul$ and $\Adl$ does not satisfy
property b), but it is in the situation described in~\cref{re:intersection}. 

Let us now consider the pair $\Bul$ and $\Bdl$. We observe that the kernel of
$\Bul$ is characterized by $\Pio v=0$. The analysis performed for the pair
$\Aul$ and $\Adl$ can be repeated and gives that in this case~\cref{ass:C} 
is verified for $\m{C}=\m{B}$.

As a consequence of the assembling of the local matrices, the global
matrices $\Au$ and $\Ad$ ($\Bu$ and $\Bd$, respectively) do not satisfy anymore
the properties listed in~\cref{ass:C}. In particular, for $k=1$ we
shall see that the matrices $\Au$ and $\Bu$ are not singular. Nevertheless, we
are going to show that the numerical results look pretty much similar to the
ones reported in~\cref{se:param}.

Moreover, in practice the matrices $\Adl$ and $\Bdl$ are not available and
they are replaced by using the local bilinear forms $S_a^P$ and $S_b^P$ given
in~\eqref{eq:abP} as follows.

Let us denote by $\mathbf{u}_h,\mathbf{v}_h\in\RE^{N_P}$ the vectors containing
the values of the $N_P$ local DoFs associated to $u_h,v_h\in\VemP$. Then, we
define the local stabilized forms as
\[
S_a^P(u_h,v_h)=\sigma_P \mathbf{u}_h^\top\mathbf{v}_h,\quad
S_b^P(u_h,v_h)=\tau_P h_P^2\mathbf{u}_h^\top\mathbf{v}_h
\]
where the stability parameters $\sigma_P$ and $\tau_P$ are positive constants
which might depend on $P$ but are independent of $h$. We point out that this
choice implies the stability requirements in~\eqref{eq:stab}. 
In the applications, the parameter $\sigma_P$ is usually
chosen depending on the mean value of the eigenvalues of the matrix
stemming from the term $a_P(\Pinabla\cdot,\Pinabla\cdot)$, and
$\tau_P$ as the mean value of the eigenvalues of the matrix resulting from
$\frac1{h^2_P}(\Pio\cdot,\Pio\cdot)_P$.
The choice of the stabilized form $S_a^P$ 
is discussed in some papers concerning the source problem, see, e.g.,
\cite{BLR} and the references therein. One can find an analysis of the
stabilization parameters $\sigma_P$ in~\cite{CMRS}.
% hence, from now on, we
%assume that the matrix $\m{A}$ associated with the bilinear
%form $a_h$ is a fixed symmetric and positive definite matrix. The bilinear
%form $b$ is more specific for eigenvalue problems and we focus on the choice of
%its stabilization. 

If $\sigma_P$ and $\tau_P$ vary in a small range, it is reasonable to take 
$\sigma_P=\alpha$ and $\tau_P=\beta$
for all $P$ and this is the situation which we discuss further. Therefore, the
structure of the matrices is $\m{A}=\Au+\alpha\Ad$ and $\m{B}=\Bu+\beta\Bd$
where $\Ad$ and $\Bd$ are the matrices with local contribution given by
$\mathbf{u}_h^\top\mathbf{v}_h$ and $h_P^2\mathbf{u}_h^\top\mathbf{v}_h$,
respectively. 
We study the behavior of the eigenvalues as $\alpha$ and $\beta$ vary in given
ranges.

\begin{center}
\begin{figure}[ht]
\includegraphics[width=0.6\textwidth]{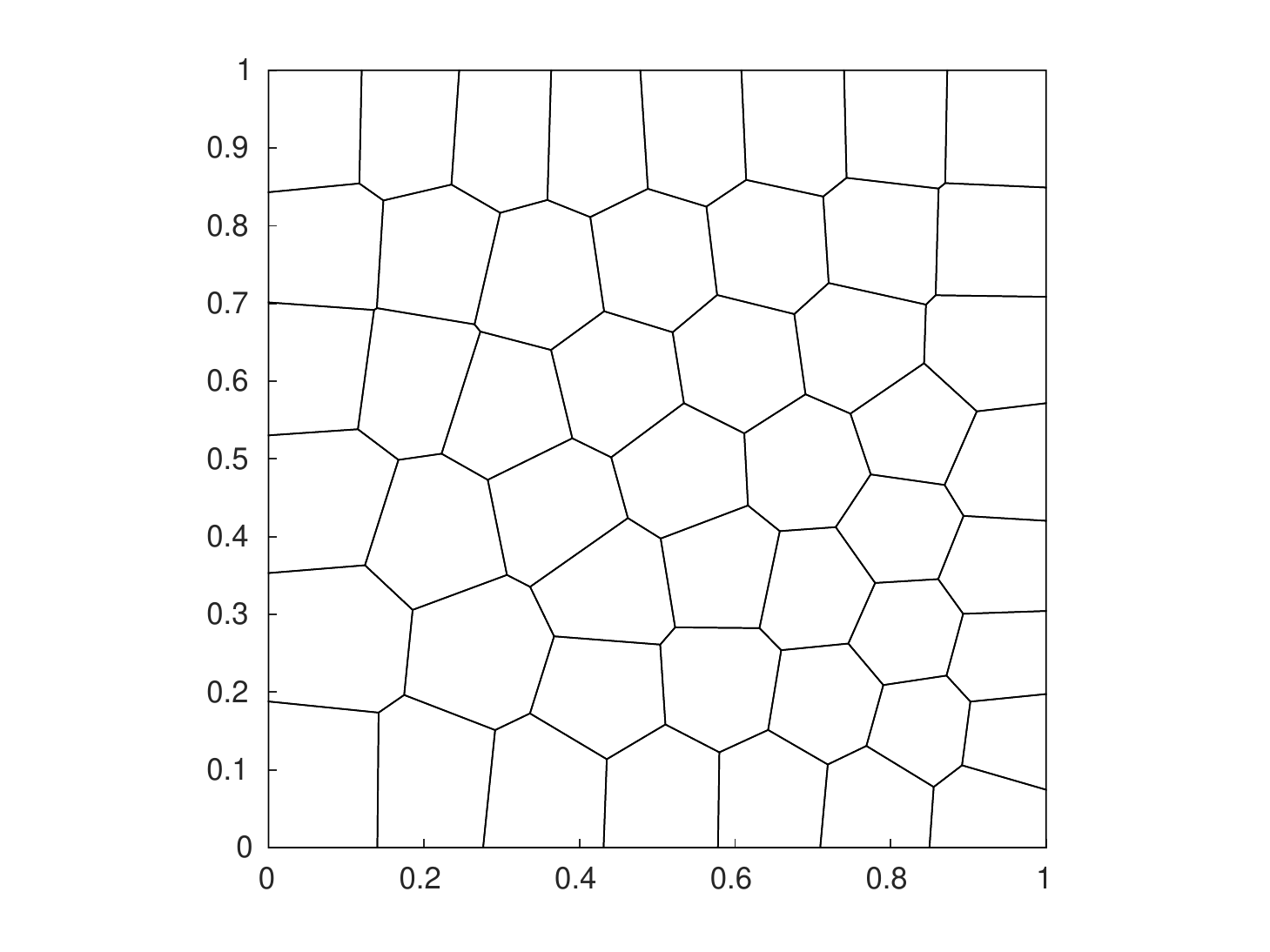}
\caption{Voronoi mesh with 50 polygons.}
\label{fig:Voronoi}
\end{figure}
\end{center}

In the following tests $\Omega$ is the unit square partitioned using a sequence
of Voronoi meshes with a given number of elements. In~\cref{fig:Voronoi} we
report the coarsest mesh with 50 elements ($h=0.2350$, 151 edges, 102
vertices). We recall that the exact
eigenvalues are given by $(i^2+j^2)\pi^2$ for $i,j\in\mathbb{N}\setminus\{0\}$
with eigenfunctions $\sin(i\pi x)\sin(j\pi y)$.
The following numerical results have been obtained using \textsc{Matlab} and, in
particular, the routine \texttt{eig} for the computation of the eigenvalues.
In the following figures, we shall always report the computed eigenvalues
divided by $\pi^2$.

\Cref{tb:kerAu} and~\cref{tb:kerBu} display the dimension of the kernel of the
matrices $\Au$
and $\Bu$ for $k=1,2,3$, and for different numbers $N$ of the elements in
the mesh.
\begin{table}
\caption{Dimension of $\KA$ with respect to $k$ and the number of elements} 
\centering
\begin{tabular}{cc cc cc cc cc cc }\toprule
$k$ && $N=50$ && $N=100$ &&$N=200$ &&$N=400$ &&$N=800$\\
\midrule
1 && 0 && 0 && 0 && 0 && 0\\
2 && 3 && 30 && 99 && 258 && 565\\
3 && 27 && 94 && 246 && 588 && 1312\\
\bottomrule
\end{tabular} 
\label{tb:kerAu}
\end{table}
\begin{table}
\caption{Dimension of $\KB$ with respect to $k$ and the number of elements} 
\centering
\begin{tabular}{cc cc cc cc cc cc }\toprule
$k$ && $N=50$ && $N=100$ &&$N=200$ &&$N=400$ &&$N=800$\\
\midrule
1 && 0 && 0 && 0 && 0 && 0\\
2 && 0 && 0 && 0 && 0 && 0\\
3 && 0 && 1 && 43 && 182 && 504\\
\bottomrule
\end{tabular} 
\label{tb:kerBu}
\end{table}
In particular we see that for $k=1$ the matrix $\Au$ is nonsingular.

We have computed the lowest eigenvalue of $\Au x=\lambda \Bu x$, 
which gives an estimate of the \emph{inf-sup} constant of the discrete problem~\eqref{eq:LaplV}. 
The results, presented in~\cref{tb:infsup}, show that the first eigenvalue is
decreasing, and this behavior corresponds to the fact that the bilinear form
$\sum_P a^P(\Pinabla\cdot,\Pinabla\cdot)$ is not stable.
\begin{table}
\caption{First eigenvalues of $\Au x=\lambda \Bu x$ for different meshes}
\centering
\begin{tabular}{c c c c c}\toprule
$N=50$ & $N=100$ &$N=200$ &$N=400$ &$N=800$\\
\midrule
1.92654e+00 & 1.74193e+00 & 1.06691e+00 & 6.81927e-01 & 5.54346e-01\\
\bottomrule
\end{tabular}
\label{tb:infsup}
\end{table}

We now discuss some tests, where we present the behavior of the eigenvalues as
the parameters $\alpha$ and $\beta$ vary, for the mesh with $N=200$ and
different degree $k$ of the polynomials in the space $\V$.
  
\begin{figure}[h]
\begin{center}

\includegraphics[width=\textwidth]{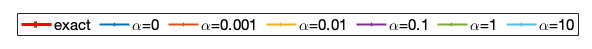}

\subfigure[\tiny{$k=1$, $\beta=0$}]
{\includegraphics[width=.32\textwidth]{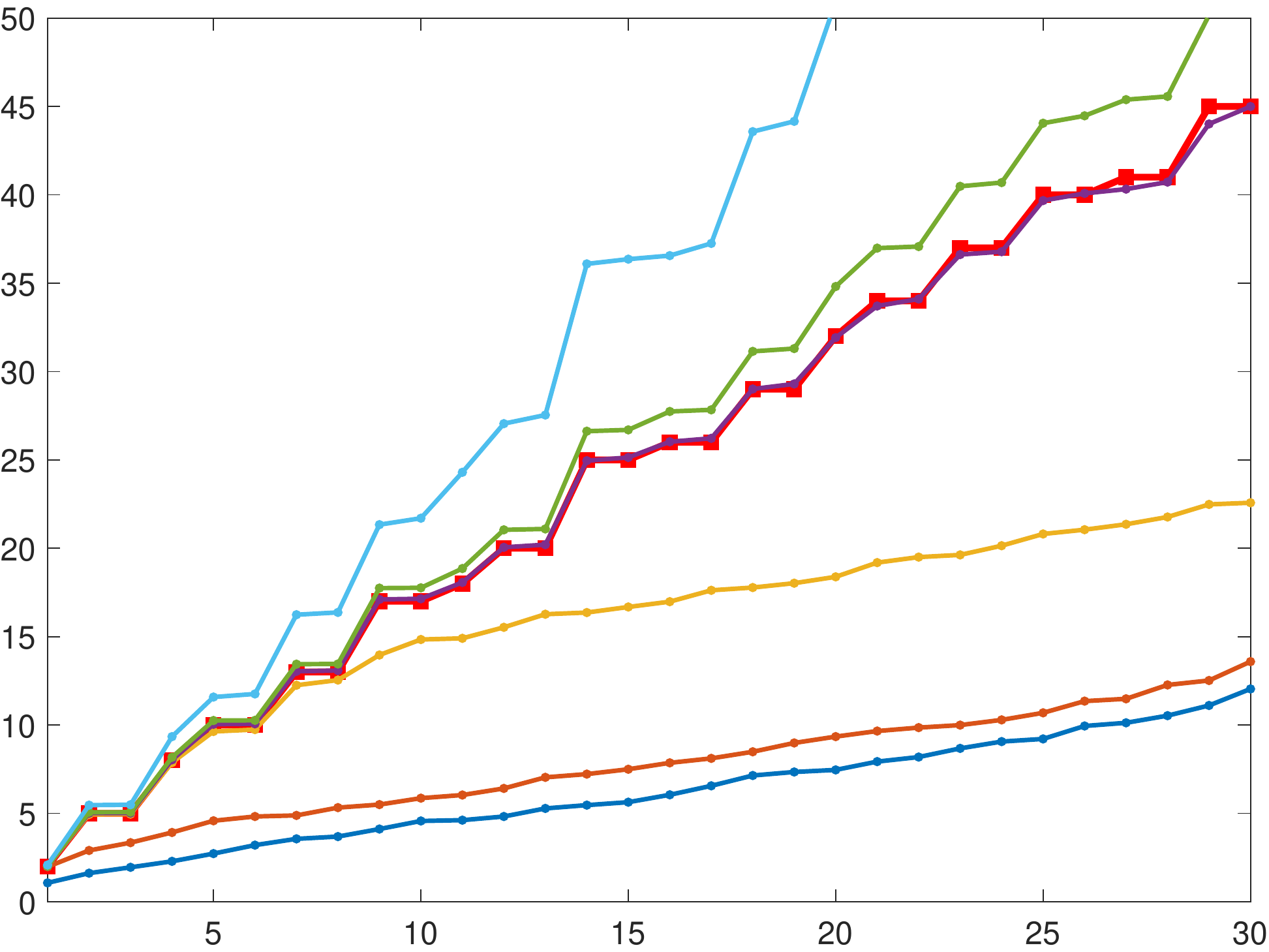}}
\subfigure[\tiny{$k=1$, $\beta=1$}]
{\includegraphics[width=.32\textwidth]{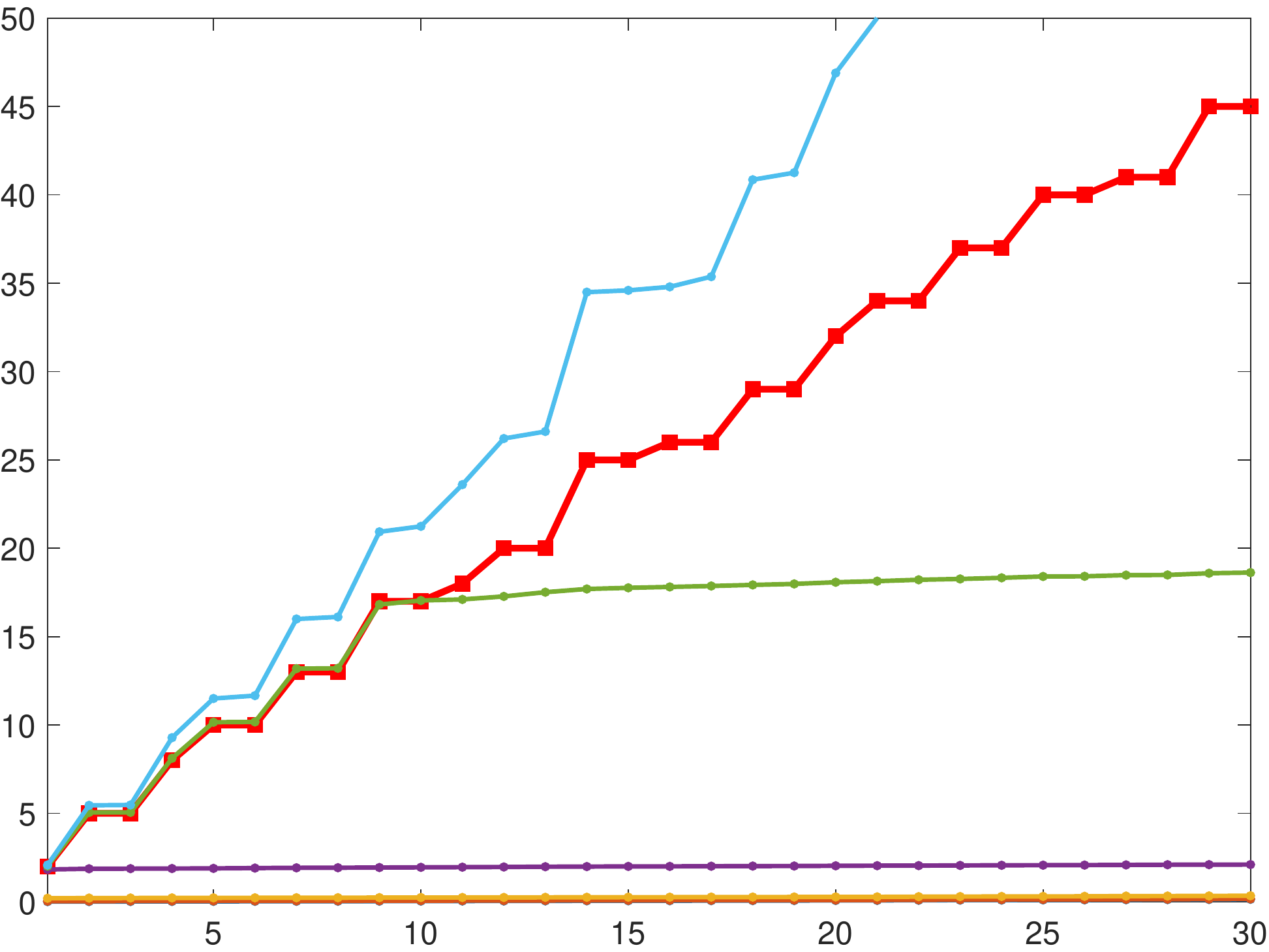}}
\subfigure[\tiny{$k=1$, $\beta=5$}]
{\includegraphics[width=.32\textwidth]{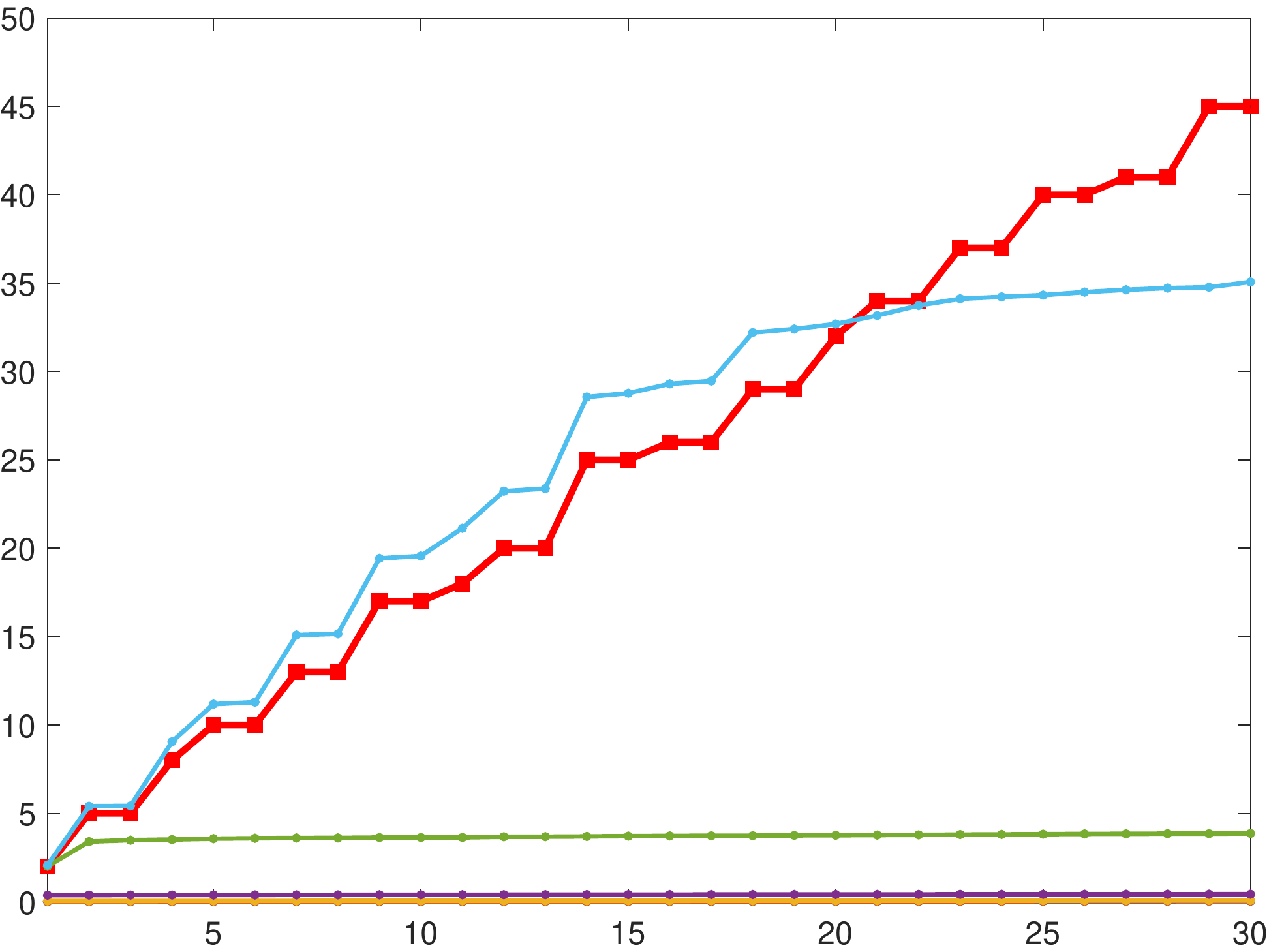}}

\subfigure[\tiny{$k=2$, $\beta=0$}]
{\includegraphics[width=.32\textwidth]{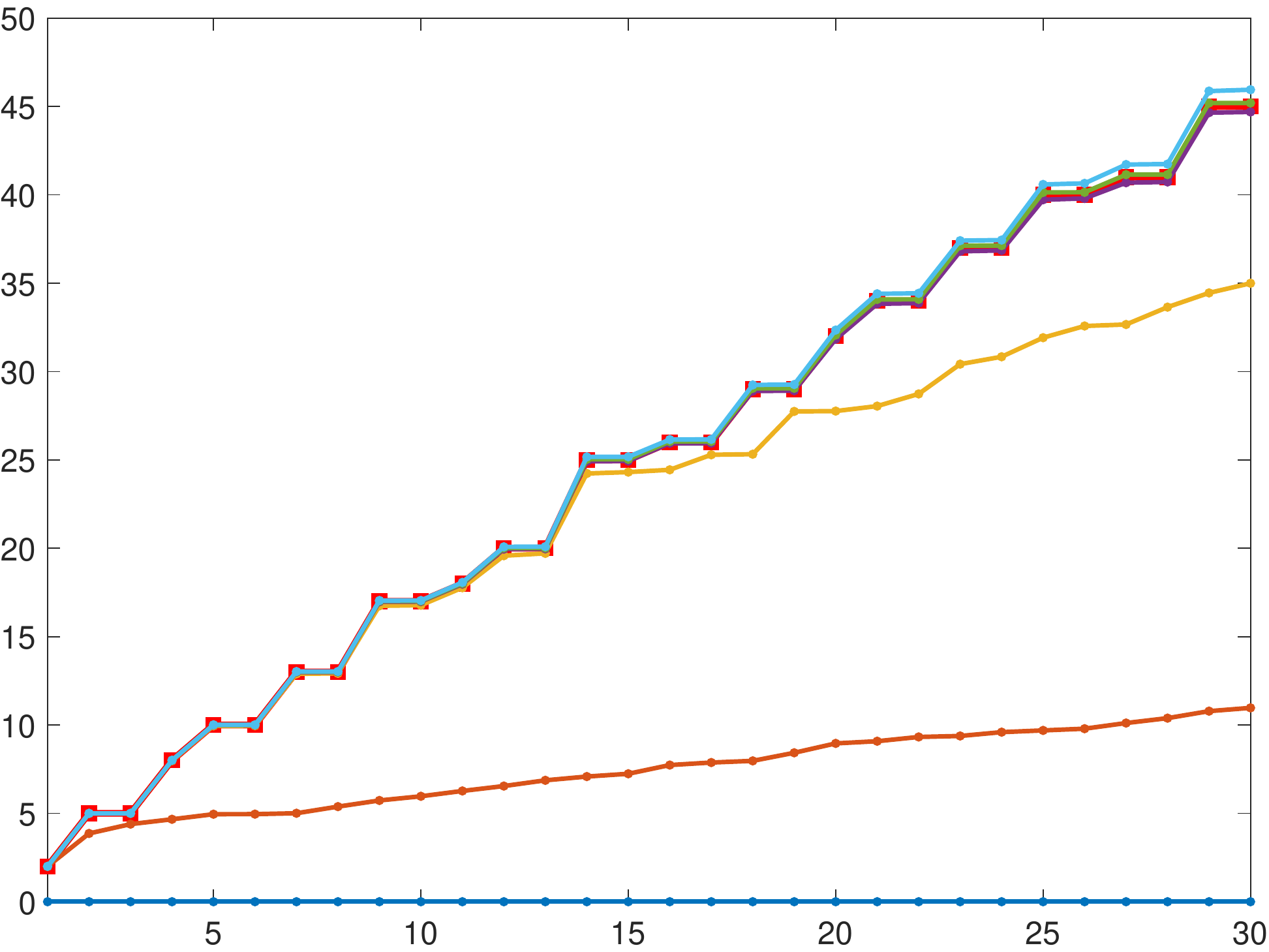}}
\subfigure[\tiny{$k=2$, $\beta=1$}]
{\includegraphics[width=.32\textwidth]{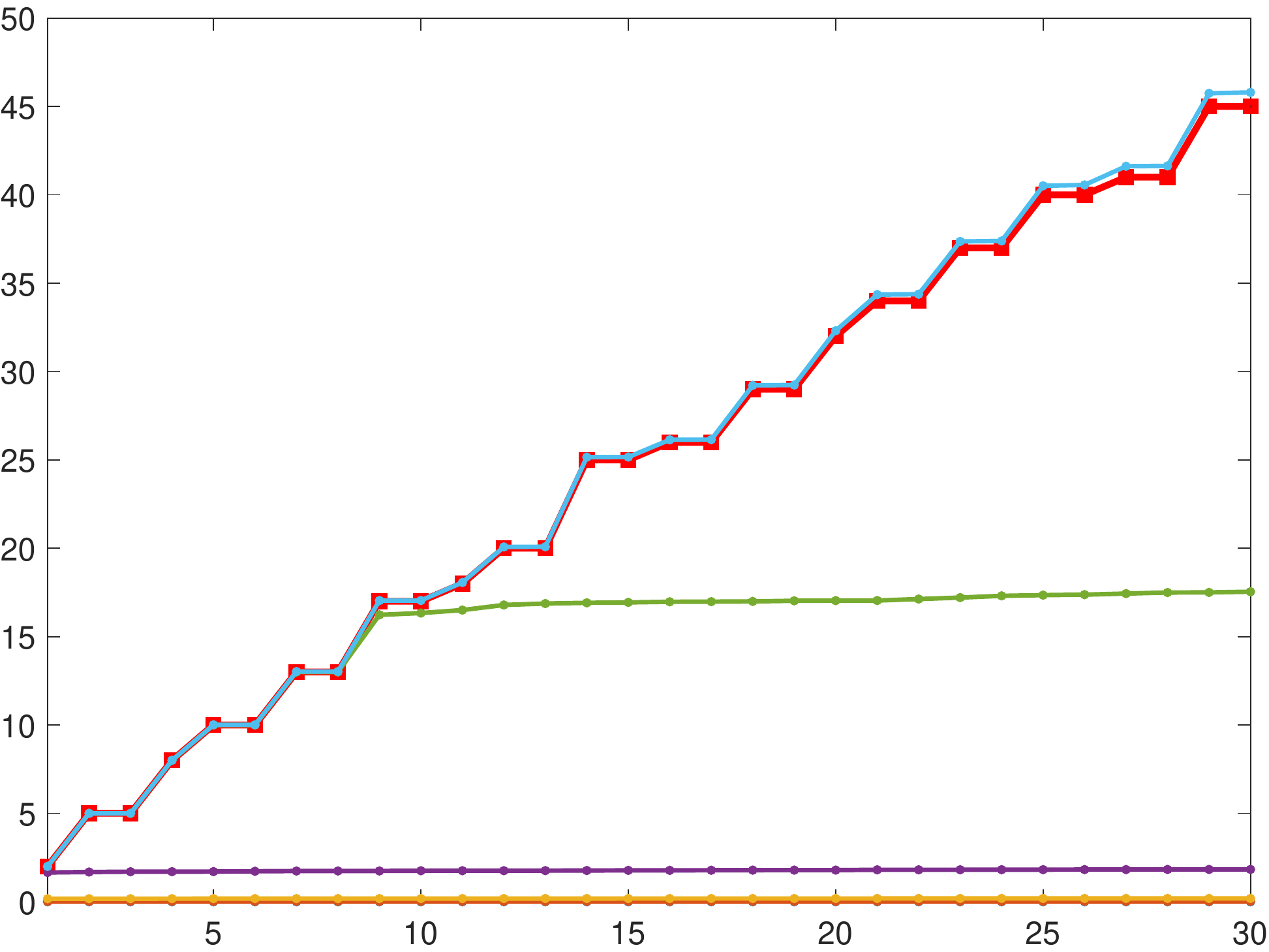}}
\subfigure[\tiny{$k=2$, $\beta=5$}]
{\includegraphics[width=.32\textwidth]{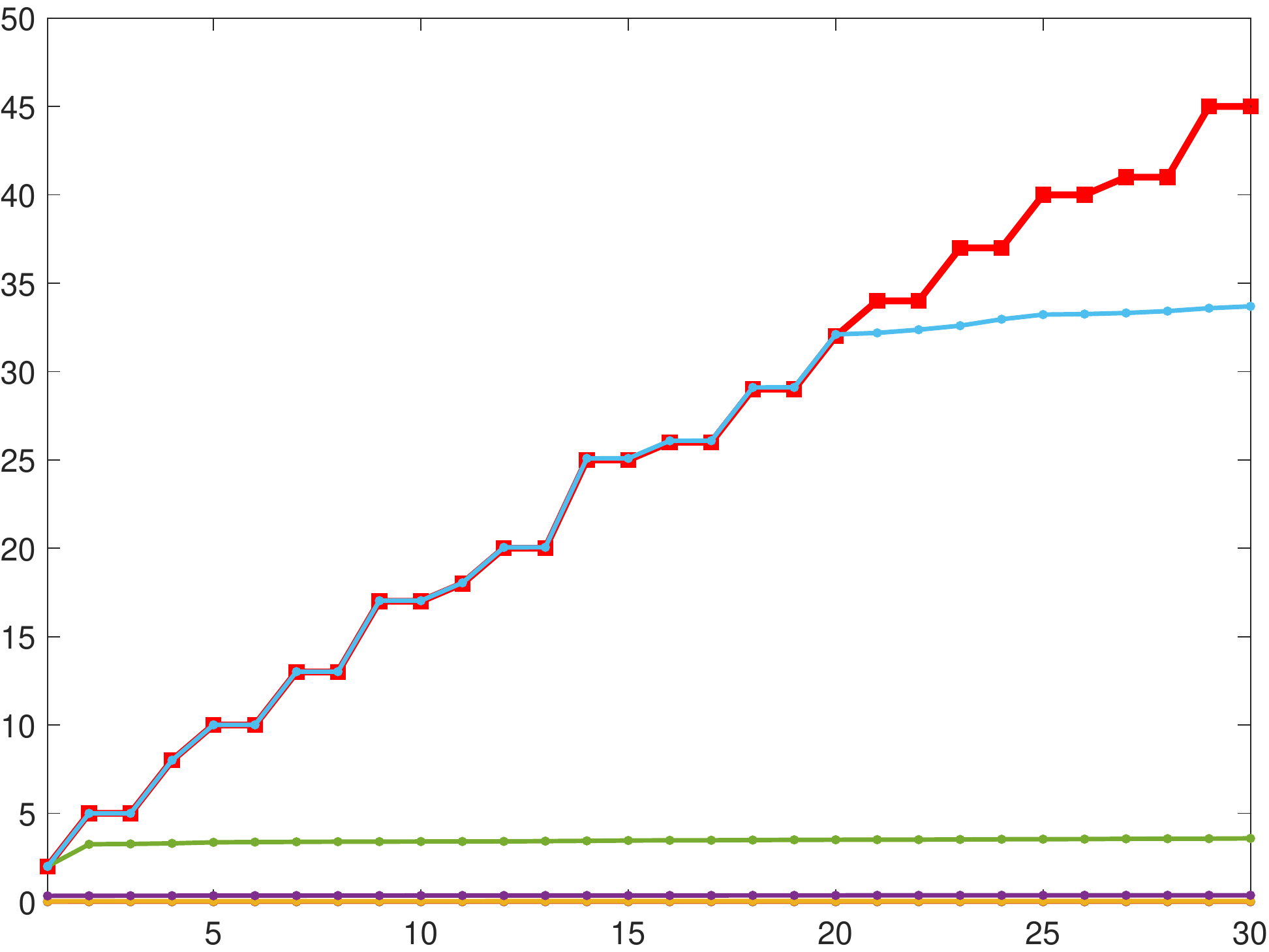}}

\subfigure[\tiny{$k=3$, $\beta=0$}]
{\includegraphics[width=.32\textwidth]{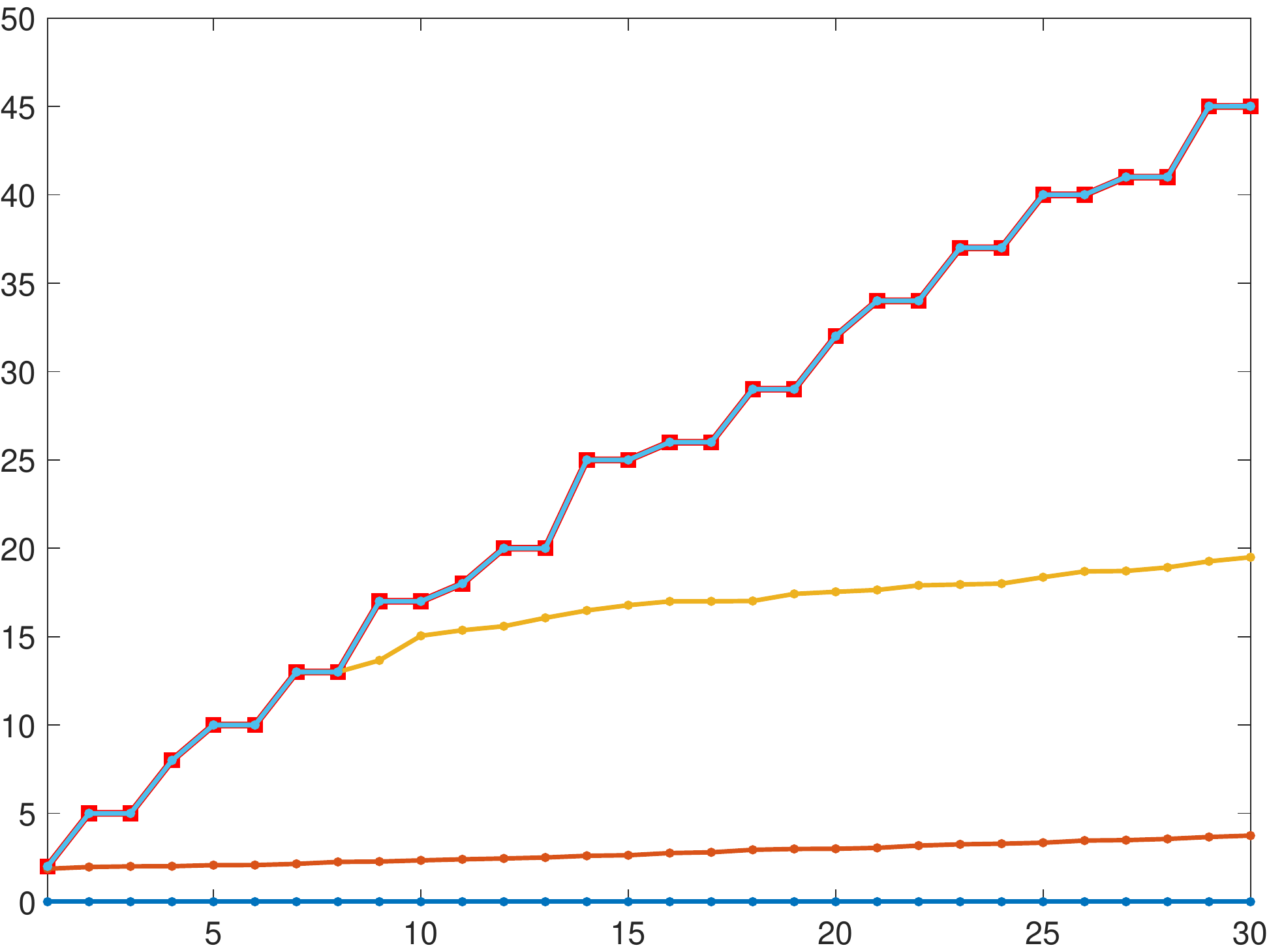}}
\subfigure[\tiny{$k=3$, $\beta=1$}]
{\includegraphics[width=.32\textwidth]{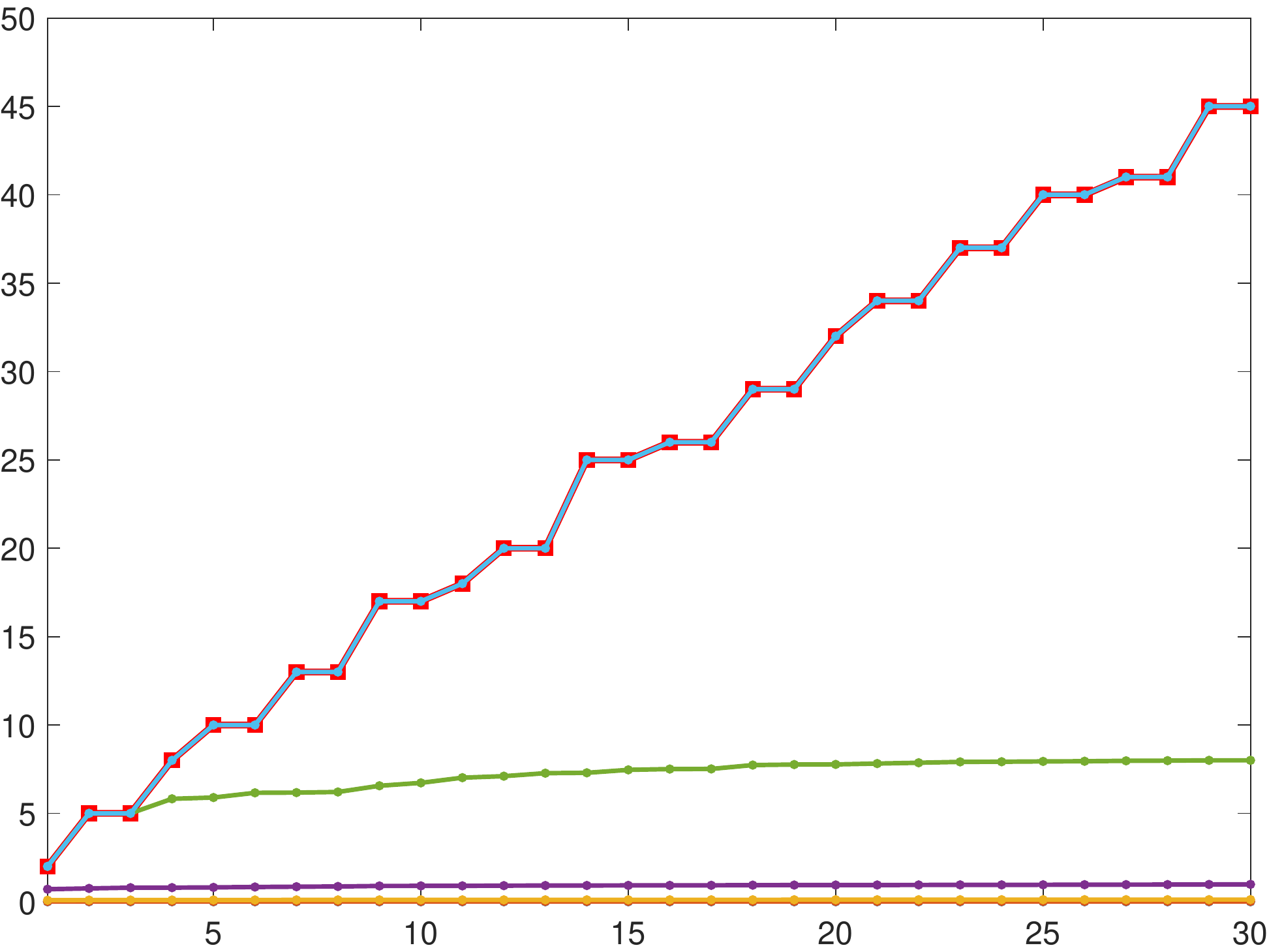}}
\subfigure[\tiny{$k=3$, $\beta=5$}]
{\includegraphics[width=.32\textwidth]{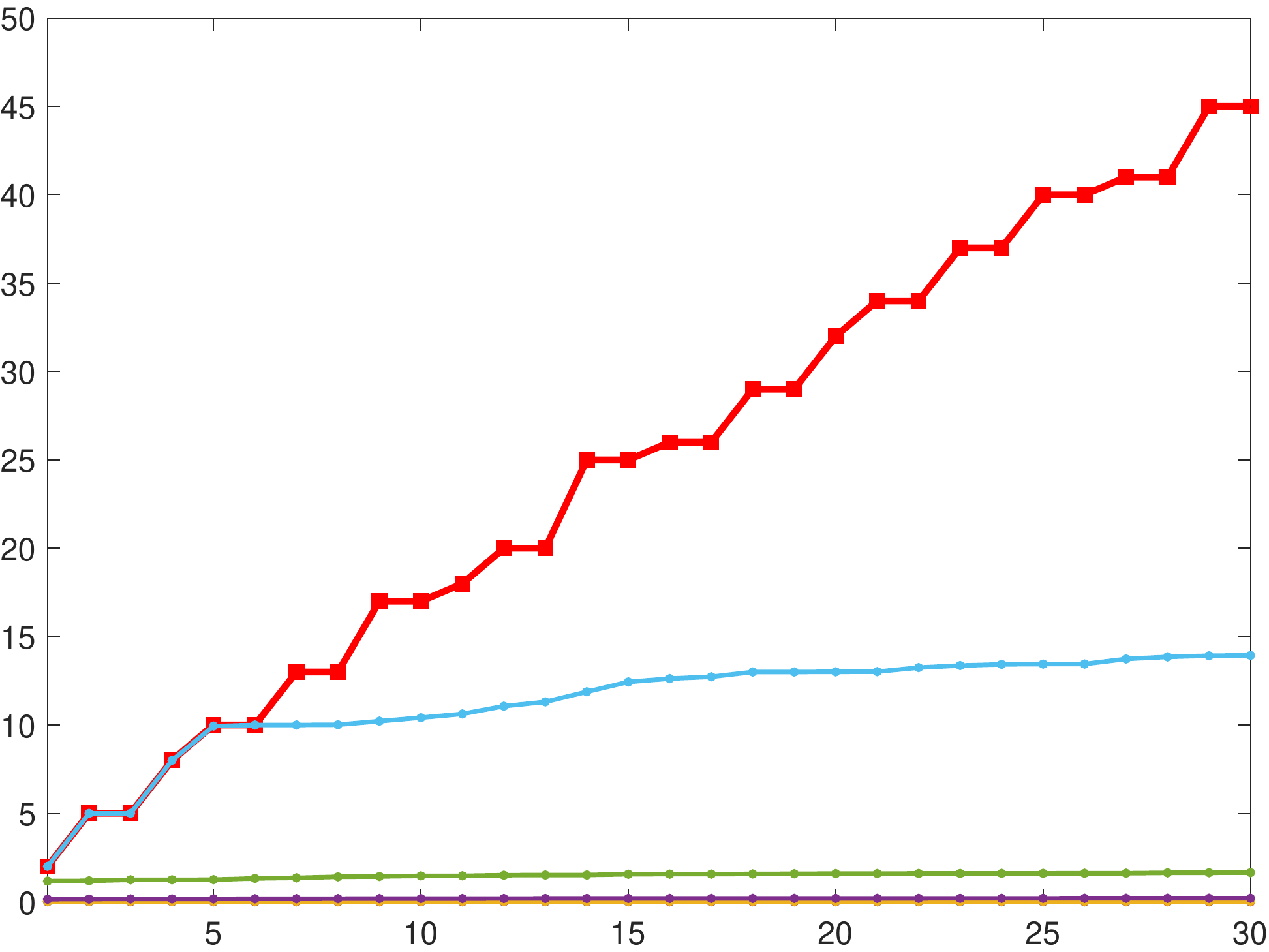}}

\end{center}
\caption{First 30 eigenvalues for different values of $k$, $\alpha$ and
$\beta$}
\label{fig:girato}
\end{figure}

The rows of~\cref{fig:girato} contain the results for fixed $k$ and the
values $\beta=0,1,5$, while, in the columns, $\beta$ is fixed and $k$ varies. 
In each picture, we plot in red the exact eigenvalues and
with different colors those corresponding to $\alpha=10^r$ with $r=-3,\dots,1$.

These plots clearly confirm that the choice of the parameters for optimal
performance is not so immediate. Consider, in particular, that we are solving
the Laplace eigenvalue problem (isotropic diffusion) on a domain as simple as a
square. For an arbitrary elliptic problem and more general domains the
situation could be much more complicated.
For $\beta=0$, the first 30 eigenvalues are well approximated with
higher degree of polynomials whenever $\alpha\ge0.1$. The value $\alpha=0.1$
seems to be the best choice in the case $k=1$. Increasing $\beta$ does not
produce much improvement. All the pictures seem to indicate that higher values
of $\alpha$ might give better results. In particular, for $k=2,3$ the first 30
eigenvalues are approximated with a reasonable accuracy for $\alpha=10$ and
$\beta=1$. Increasing $\beta$ and keeping $\alpha=10$, we see that a smaller
number of eigenvalues are captured.

\begin{figure}[h]
\begin{center}
\subfigure[\tiny{$k=1$, $\beta=0$, $\alpha\in[0,10]$}]
{\includegraphics[width=.32\textwidth]{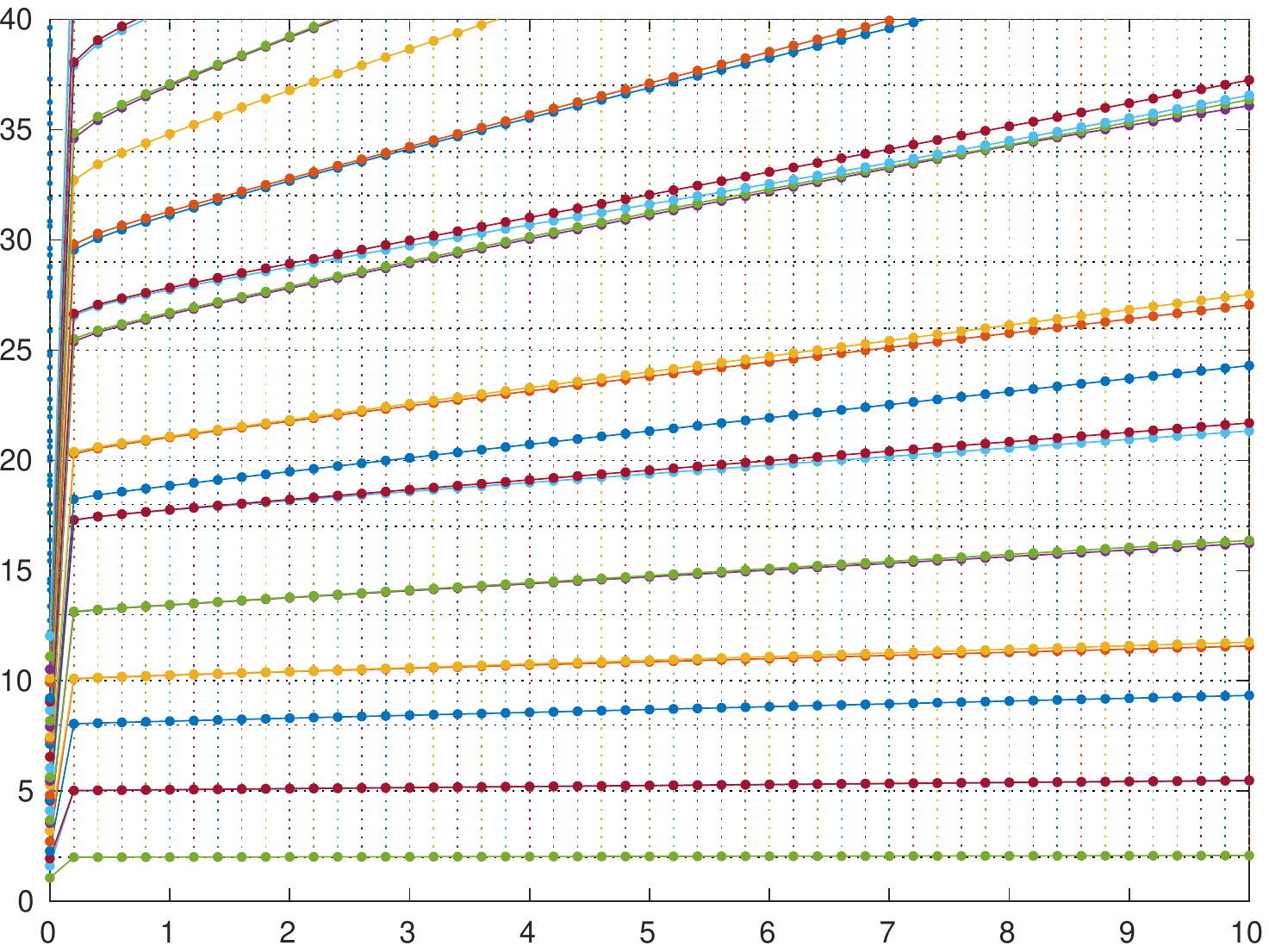}}
\subfigure[\tiny{$k=1$, $\beta=1$, $\alpha\in[0,10]$}]
{\includegraphics[width=.32\textwidth]{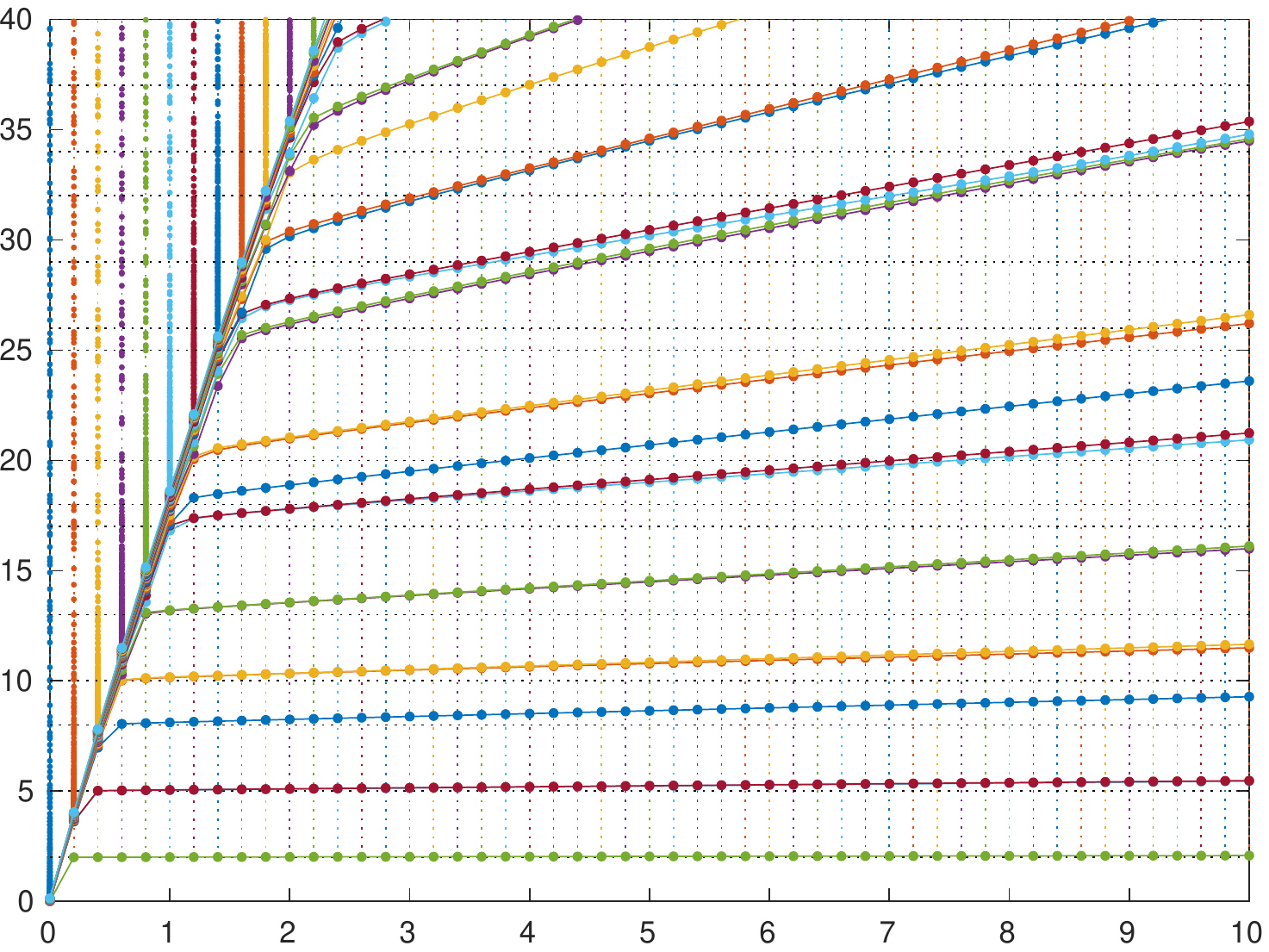}}
\subfigure[\tiny{$k=1$, $\beta=5$, $\alpha\in[0,10]$}]
{\includegraphics[width=.32\textwidth]{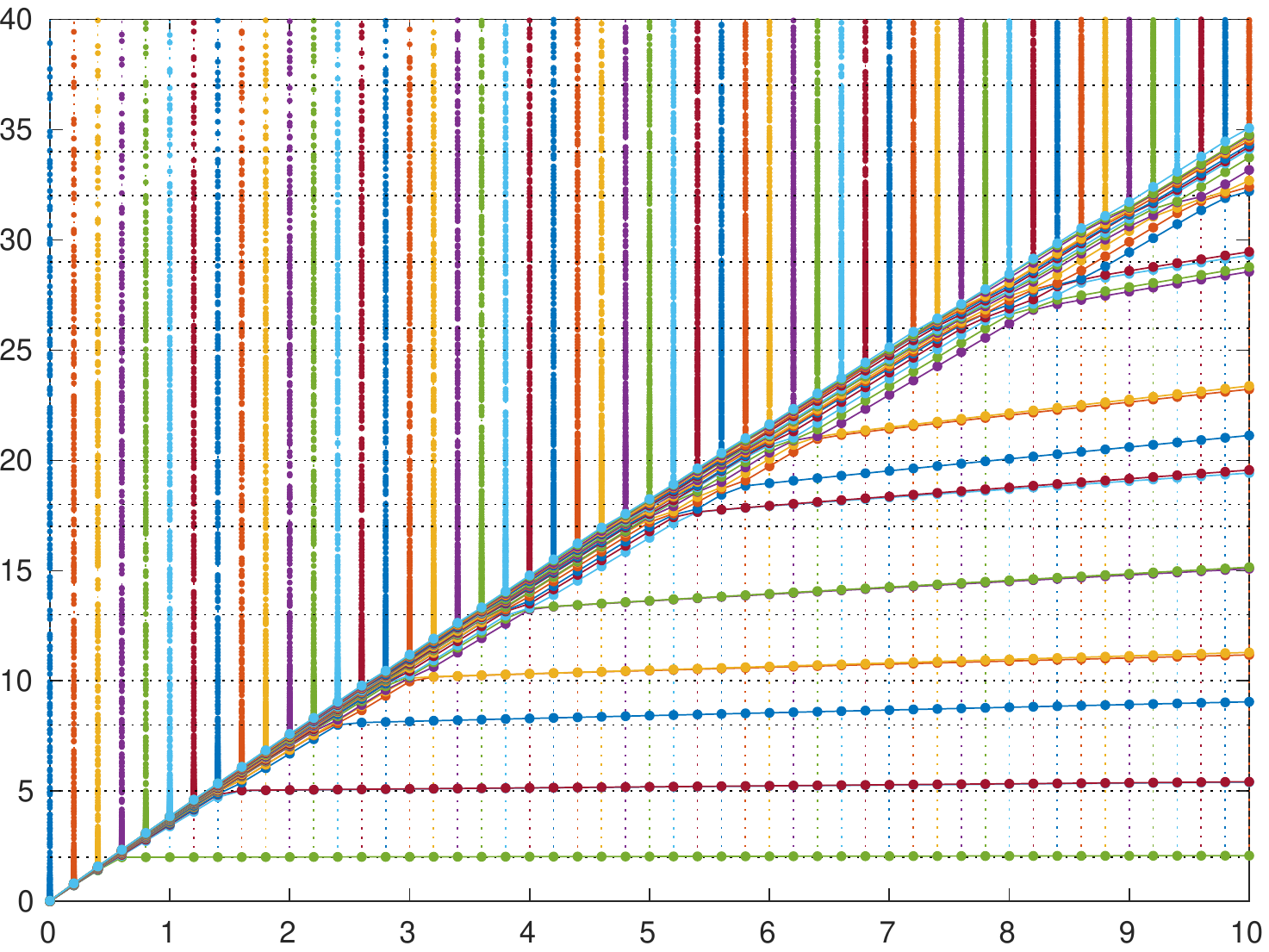}}

\subfigure[\tiny{$k=2$, $\beta=0$, $\alpha=[0,10]$}]
{\includegraphics[width=.32\textwidth]{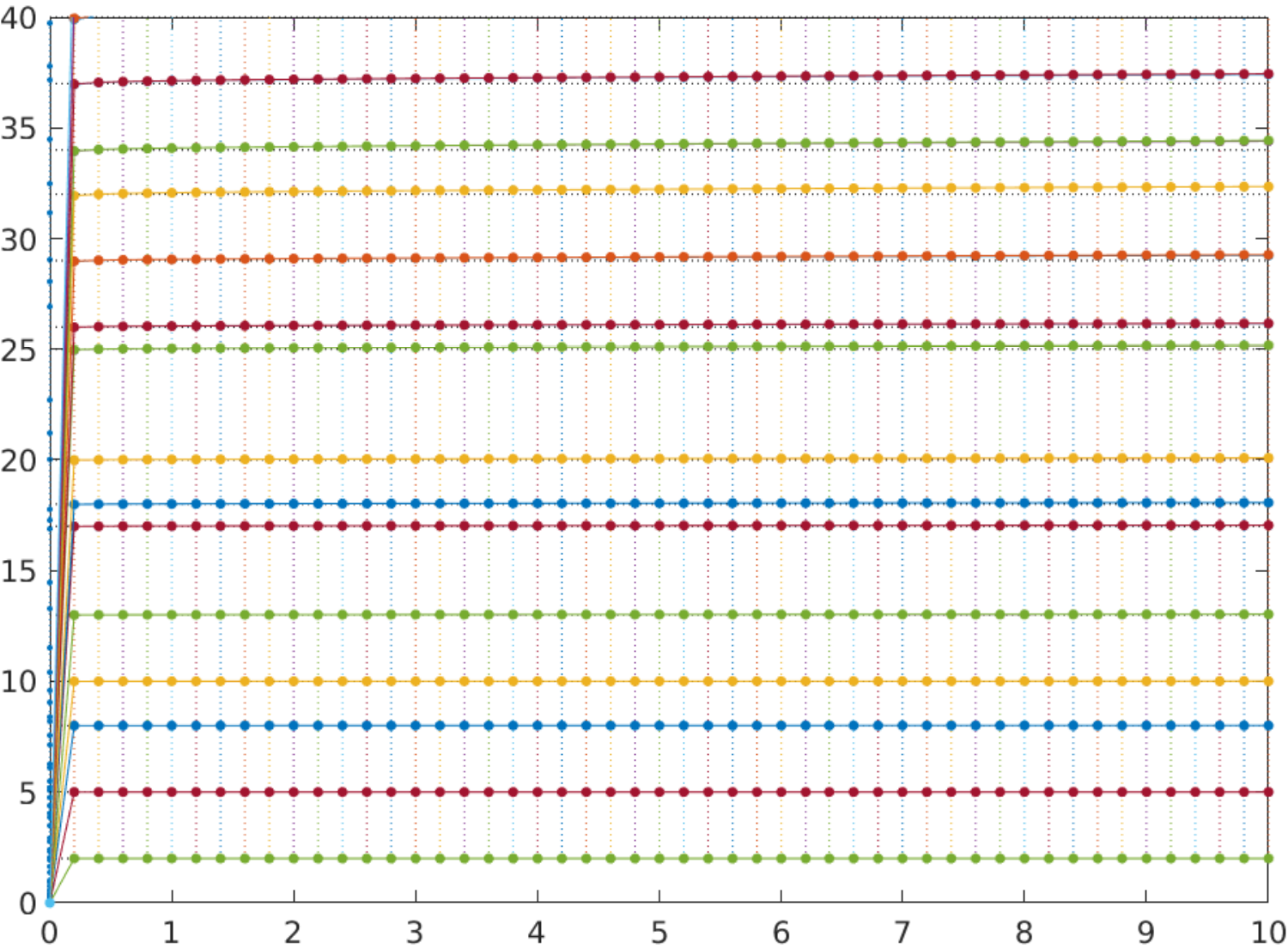}}
\subfigure[\tiny{$k=2$, $\beta=1$, $\alpha=[0,10]$}]
{\includegraphics[width=.32\textwidth]{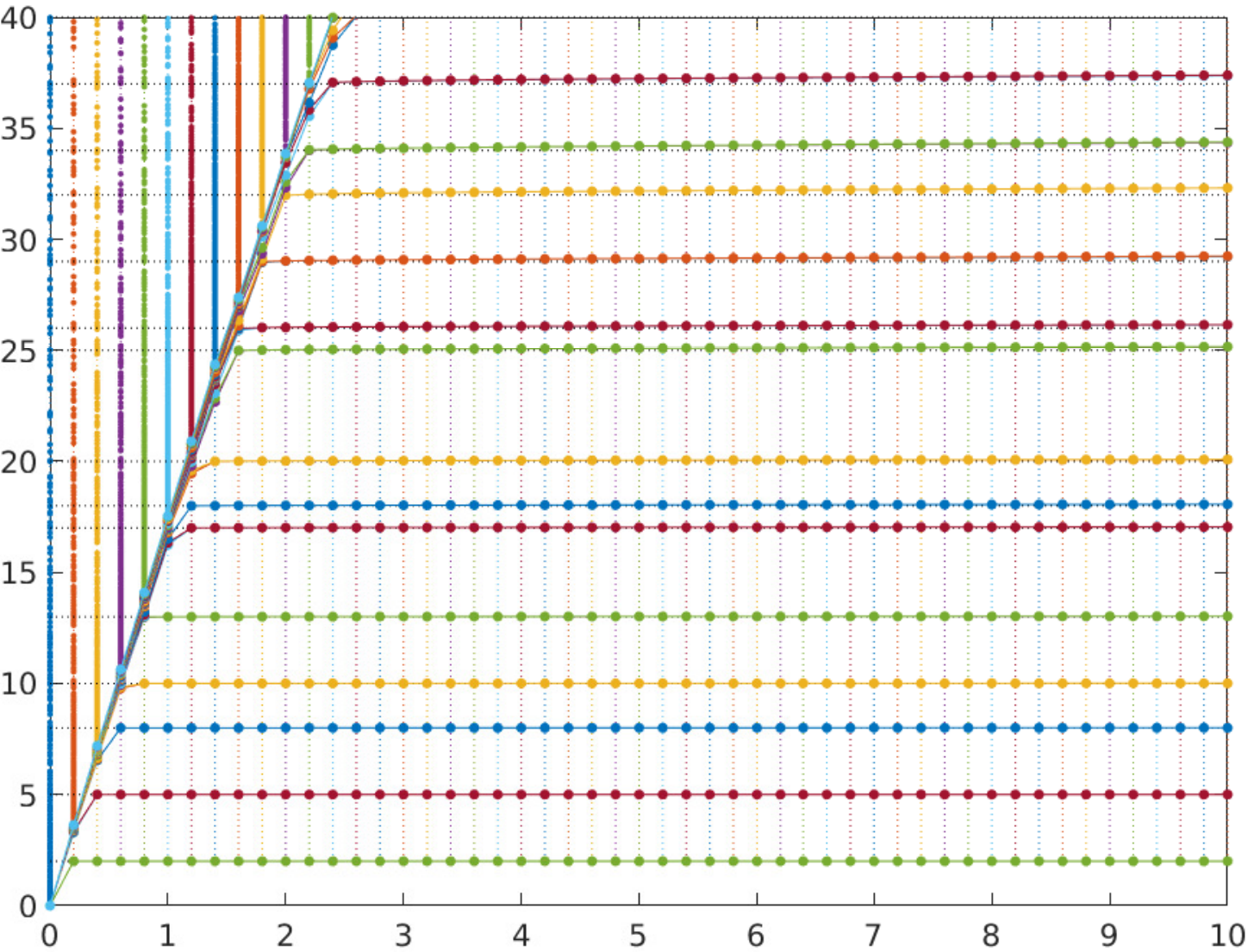}}
\subfigure[\tiny{$k=2$, $\beta=5$, $\alpha=[0,10]$}]
{\includegraphics[width=.32\textwidth]{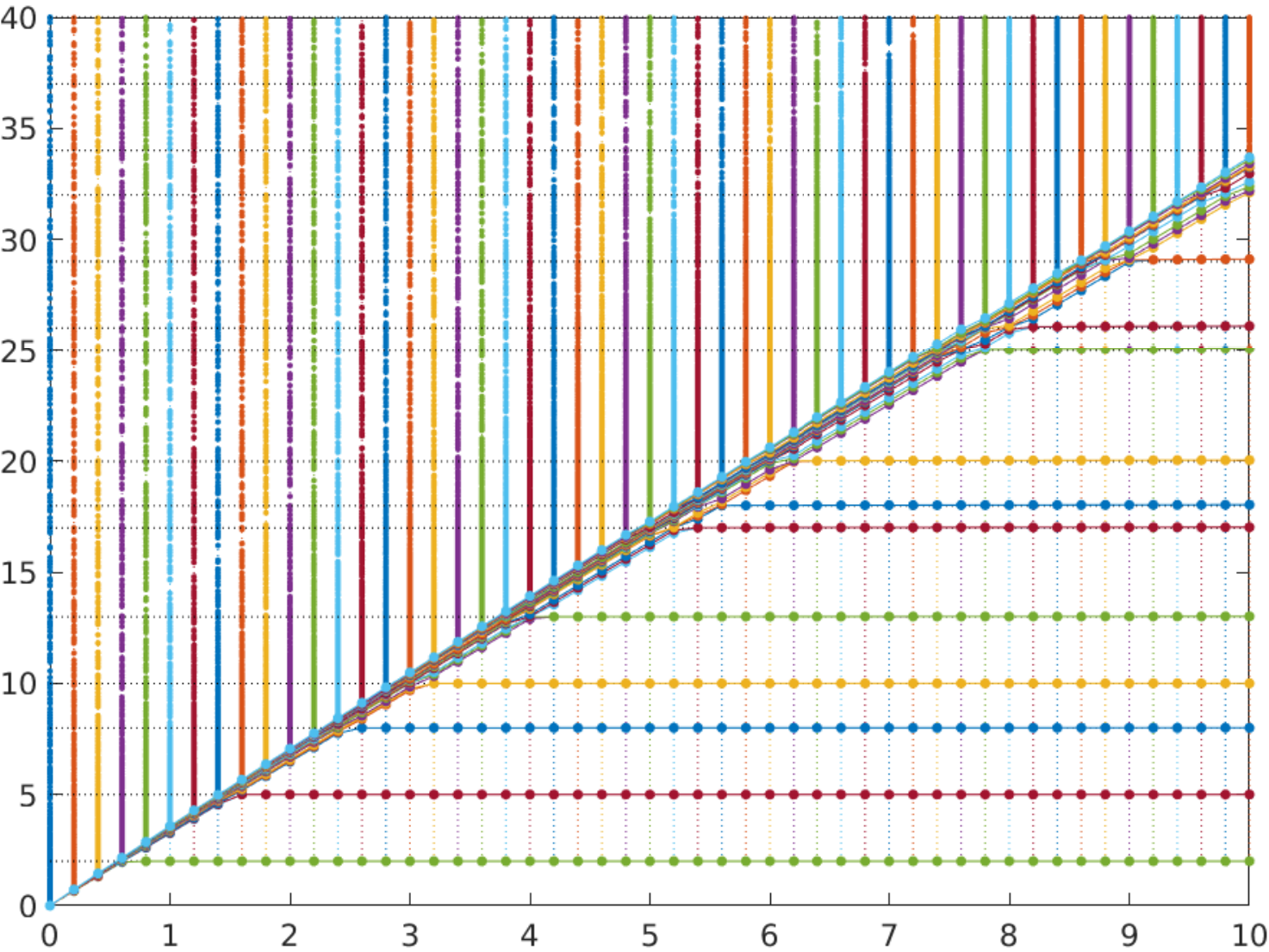}}

\subfigure[\tiny{$k=3$, $\beta=0$, $\alpha=[0,10]$}]
{\includegraphics[width=.32\textwidth]{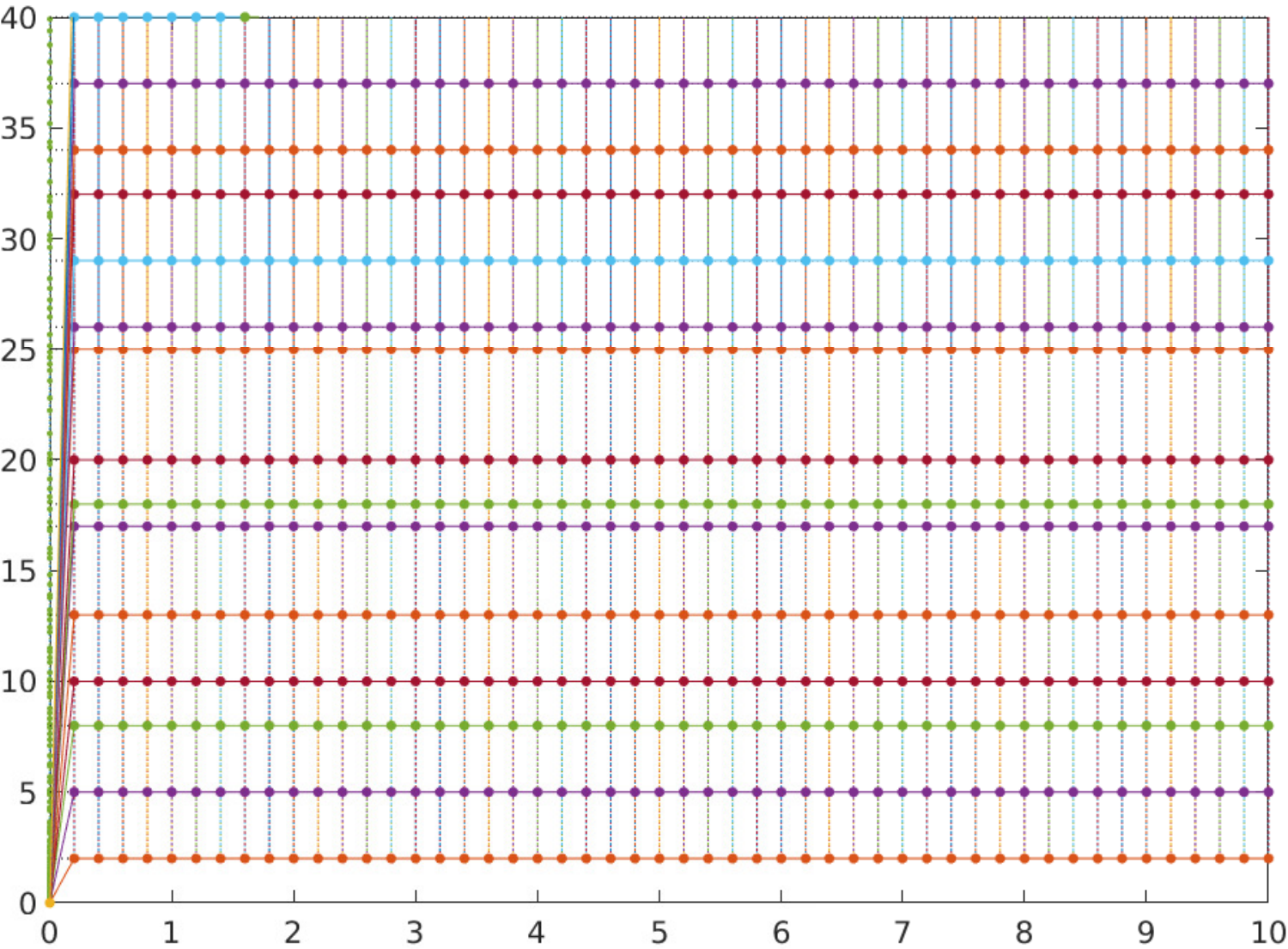}}
\subfigure[\tiny{$k=3$, $\beta=1$, $\alpha=[0,10]$}]
{\includegraphics[width=.32\textwidth]{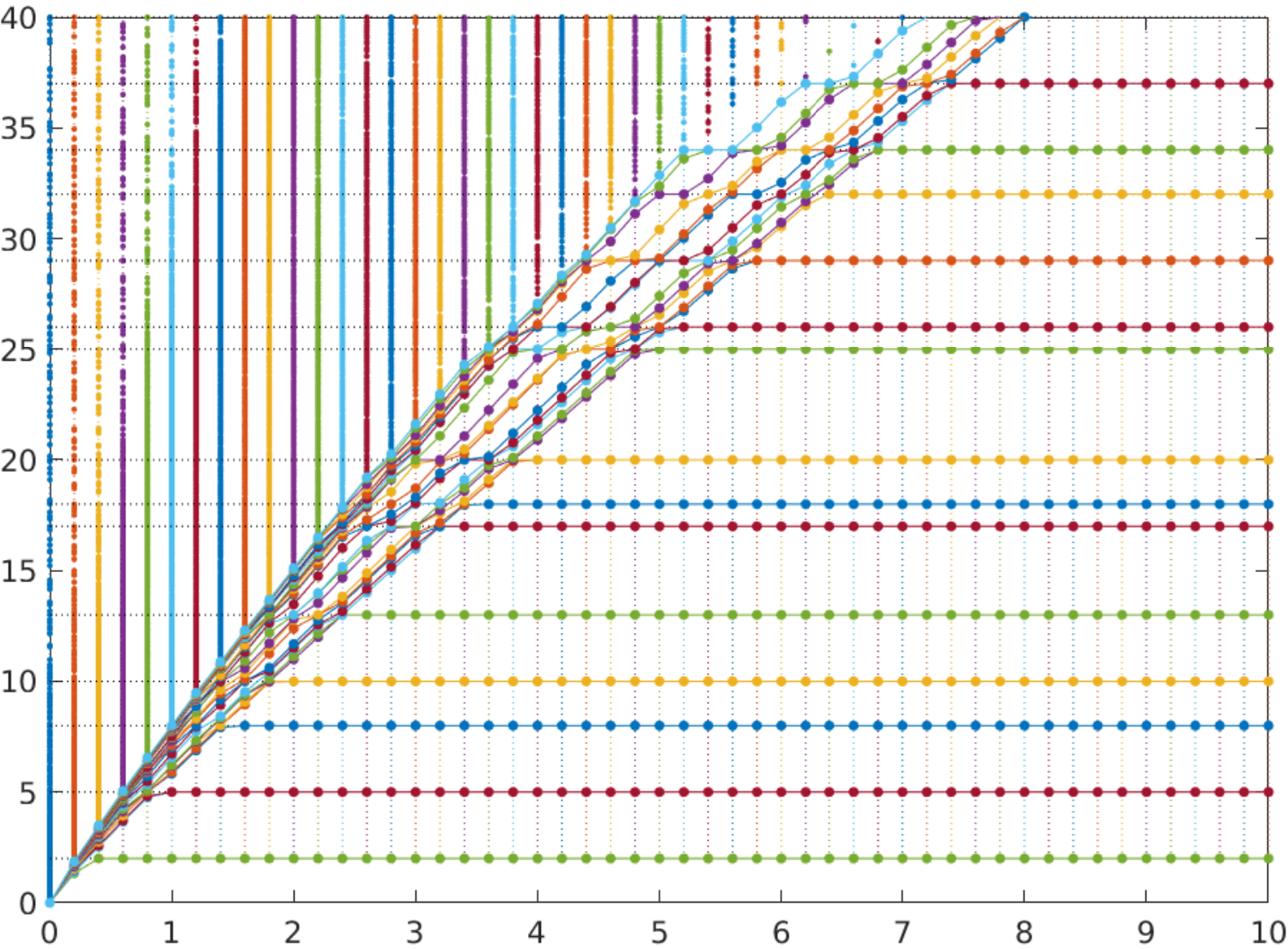}}
\subfigure[\tiny{$k=3$, $\beta=5$, $\alpha=[0,10]$}]
{\includegraphics[width=.32\textwidth]{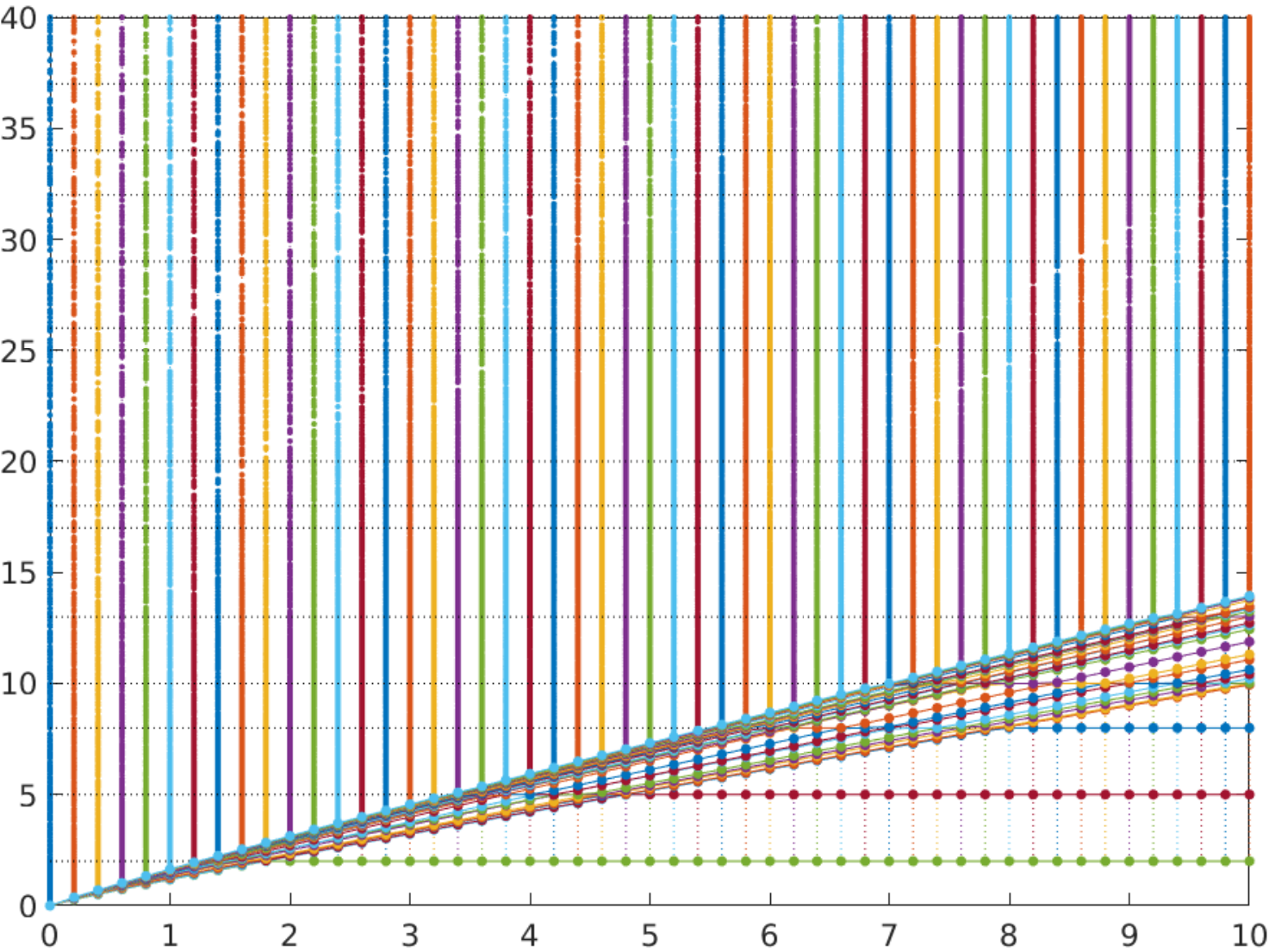}}
\end{center}
\caption{Eigenvalues versus $\alpha$ for different values of $k$ and $\beta$}
\label{fig:alfa}
\end{figure}

\Cref{fig:alfa} shows the behavior of the eigenvalues as $\alpha$ varies
from $0$ to $10$. At a first glance the pictures remind of~\cref{fig:case1} (left) even if, 
as it has been explained before, the situation is not exactly matching what we 
discussed in~\cref{se:param}.

Each subplot reports all computed eigenvalues between $0$ and $40$; the dotted
horizontal lines represent the exact solutions. The first $30$ computed
eigenvalues are connected together with lines of different colors in an
automated way. An ``ideal'' good approximation would correspond to a series of
colored lines matching the dotted lines of the exact eigenvalues.
It is interesting to look at the differences between various degrees ($k$ from
$1$ to $3$ moving from the top to the bottom) and values of $\beta$ (equal to
$0$, $1$, and $5$ from left to right).

\begin{figure}
\includegraphics[width=.70\textwidth]{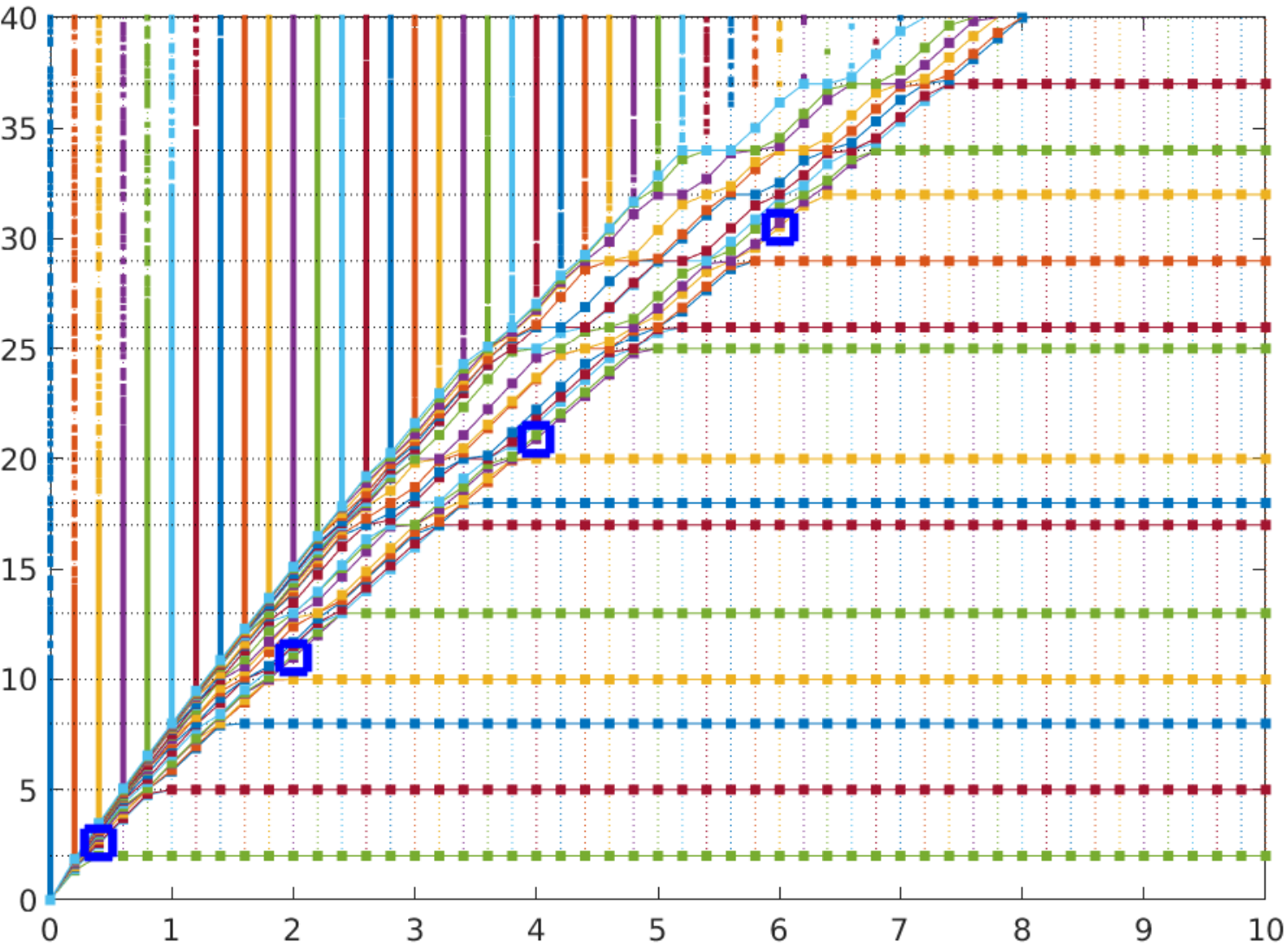}
\caption{Same plot as in~\cref{fig:alfa}(h) with four marked (spurious)
eigenvalues}
\label{fig:marked}
\end{figure}

\begin{figure}
\begin{center}
\includegraphics[width=.45\textwidth]{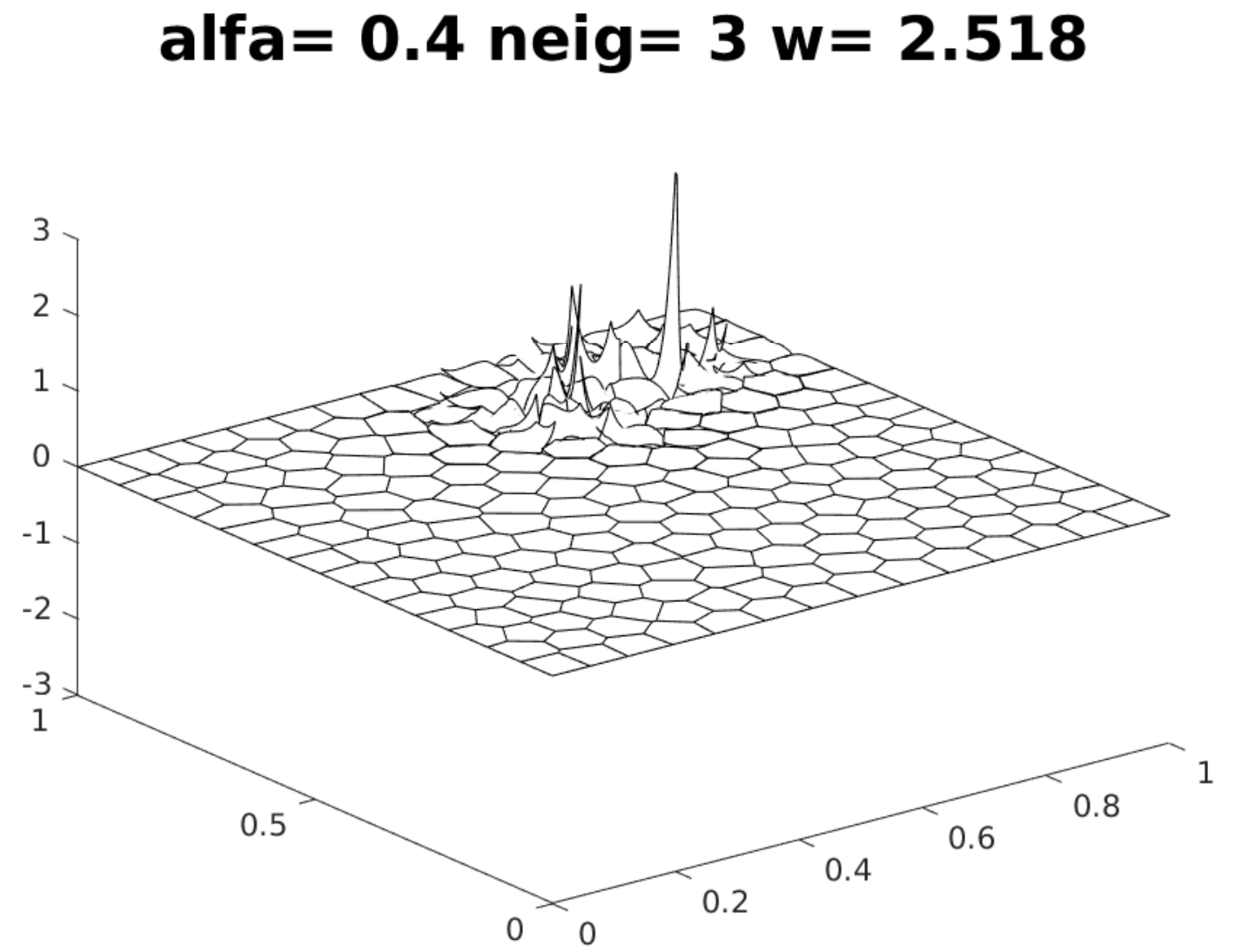}
\includegraphics[width=.45\textwidth]{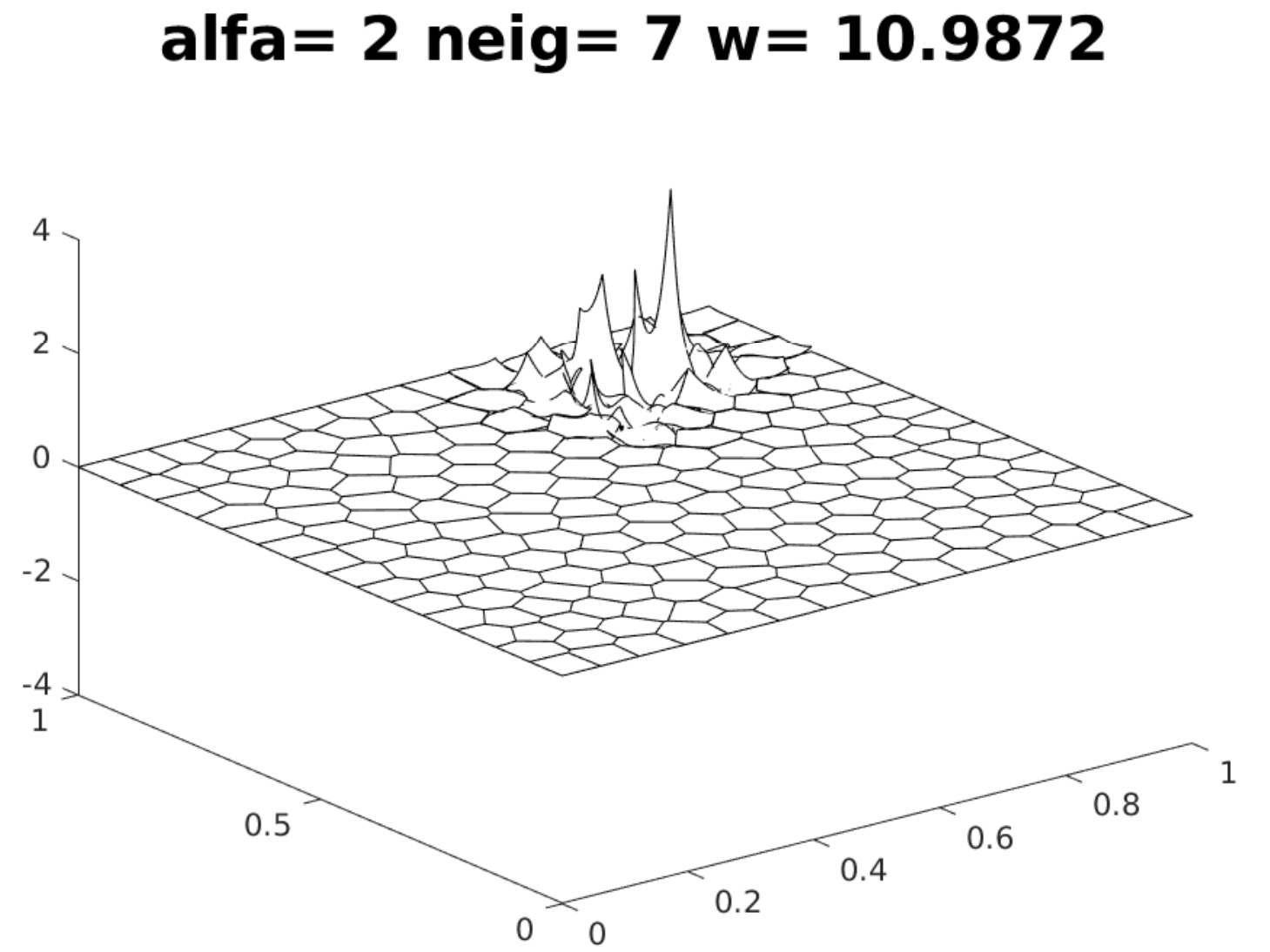}

\includegraphics[width=.45\textwidth]{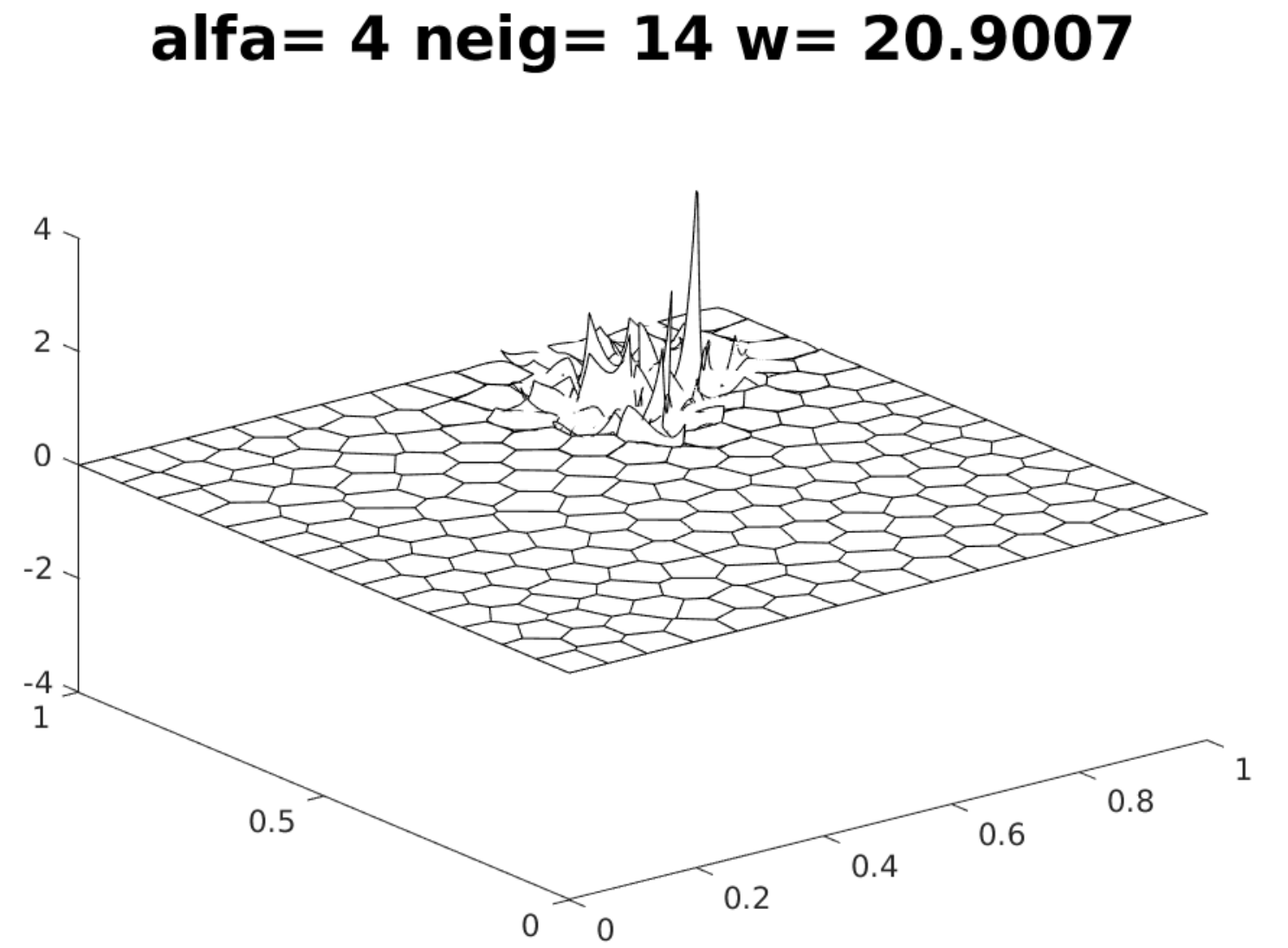}
\includegraphics[width=.45\textwidth]{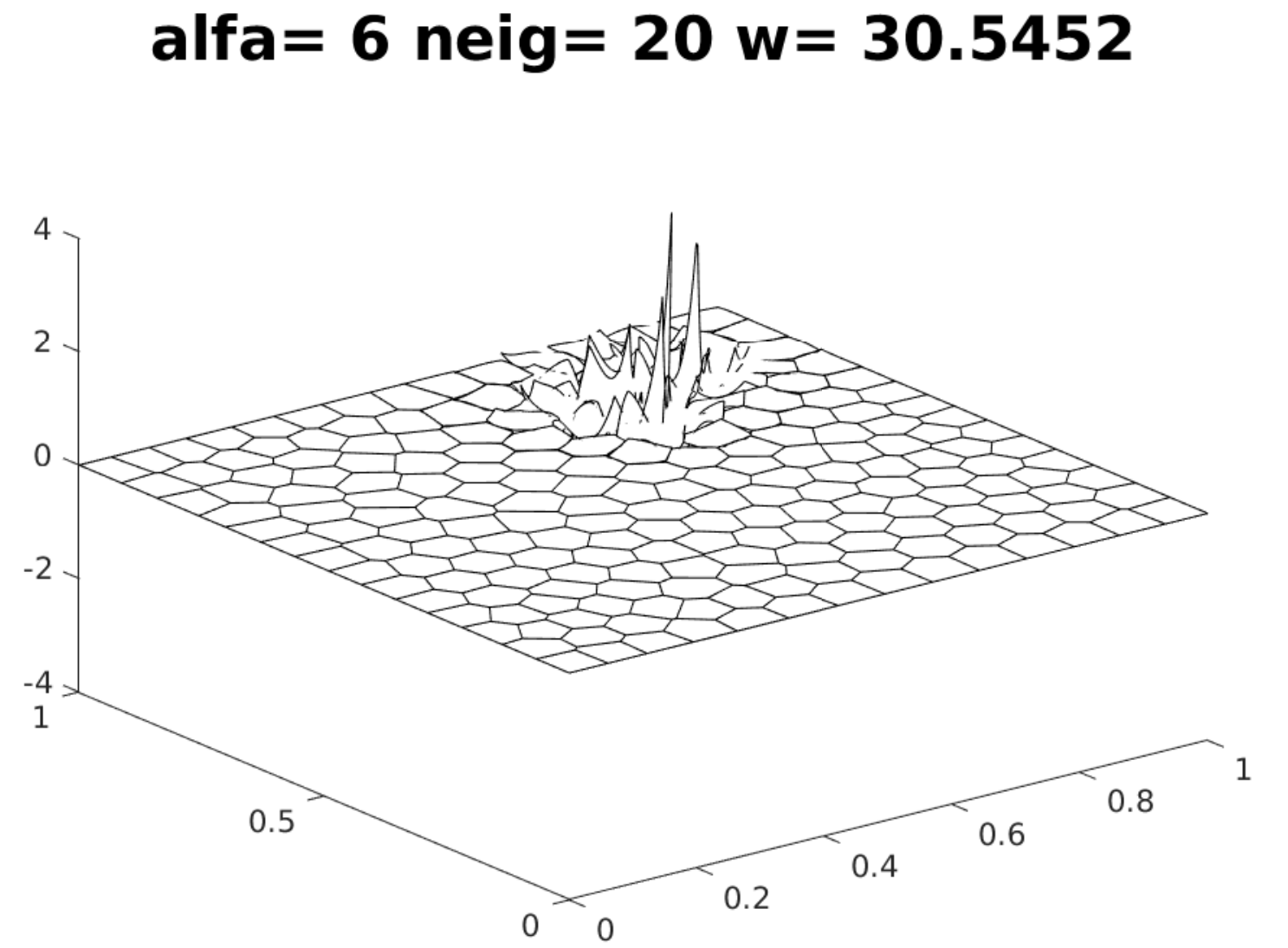}
\end{center}
\caption{Eigenfunctions corresponding to the eigenvalues marked 
in~\cref{fig:marked}}
\label{fig:autof}
\end{figure}

More reliable results seem to be obtained for large $k$ and small $\beta$.
Actually, the limit case of $\beta=0$ appears to be the safest choice. This is
in agreement with the claim of~\cite{BMRR} where the authors remark that
``even the value $\sigma_E=0$ yields very accurate results, in spite of the
fact that for such a value of the parameter the stability estimate and hence
most of the proofs of the theoretical results do not hold'' (note that
$\sigma_E=0$ in~\cite{BMRR} has the same meaning as $\beta$ in our paper). It
is interesting to observe that the analysis of~\cite{GV}, summarized in~\cref{th:conv}, 
covers the case $\beta=0$ as well.
On the other hand $\beta=0$ may produce a singular matrix $\m{B}$ and this
could be not convenient from the computational point of view.

In order to better understand the behavior of the eigenvalues reported in~\cref{fig:alfa}(h), 
we highlight in~\cref{fig:marked} four eigenvalues that are apparently aligned along an oblique line. 
The corresponding eigenfunctions are reported in~\cref{fig:autof}. The four
eigenfunctions look similar, so that the analogy with~\cref{fig:case1}
(left) is even more evident.

\begin{figure}
\begin{center}
\includegraphics[width=.45\textwidth]{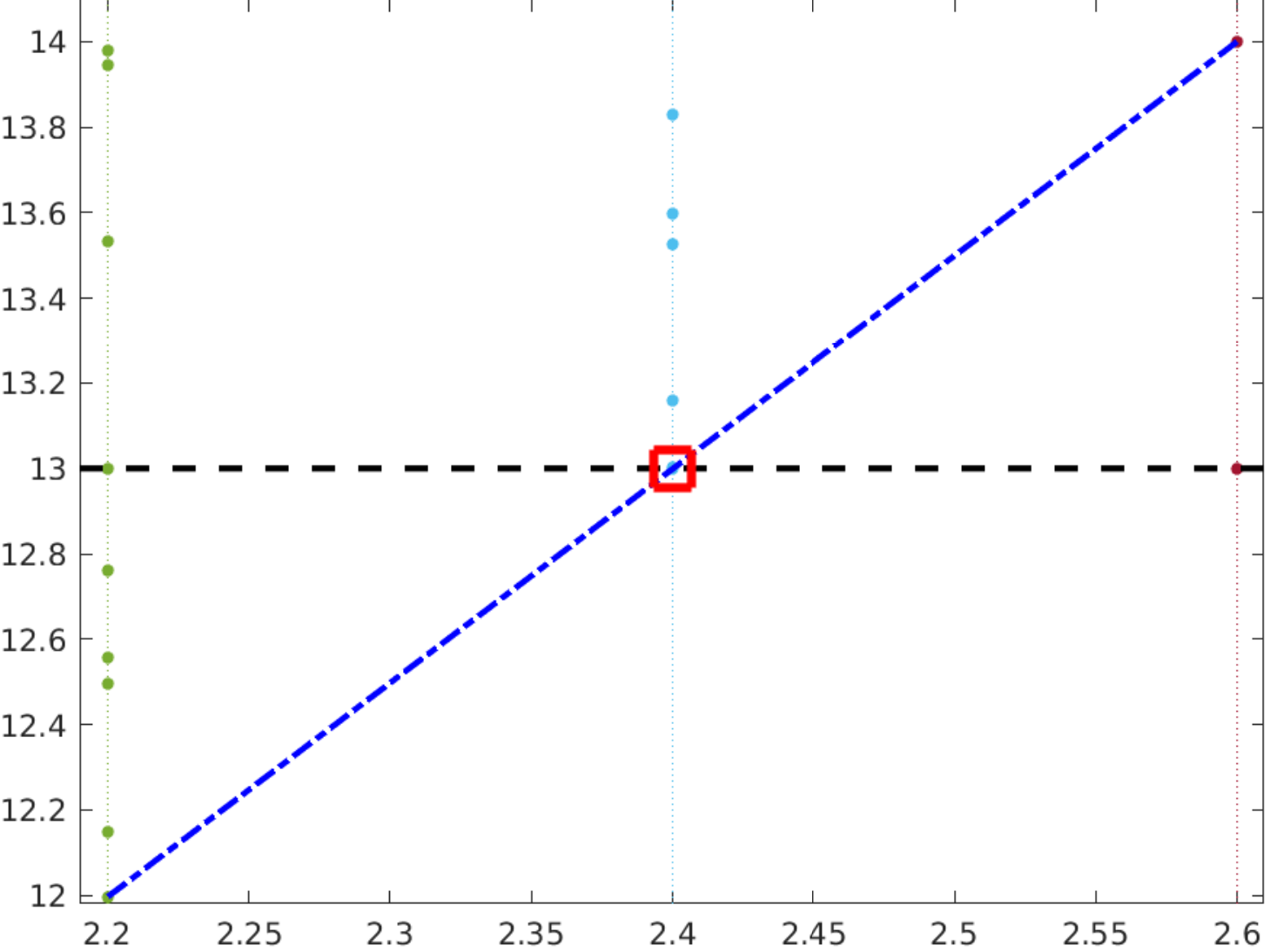}
\includegraphics[width=.45\textwidth]{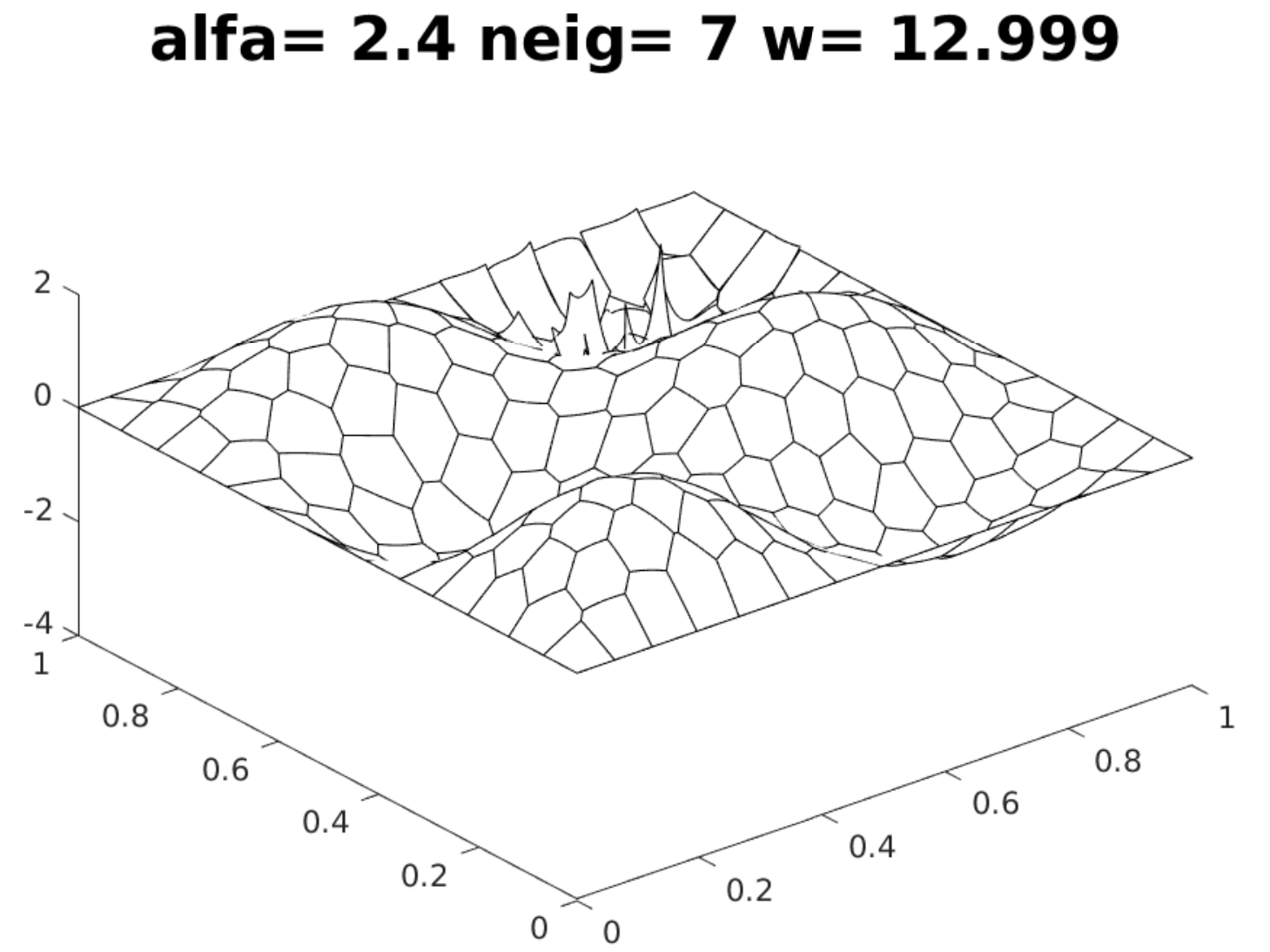}

\medskip

\includegraphics[width=.45\textwidth]{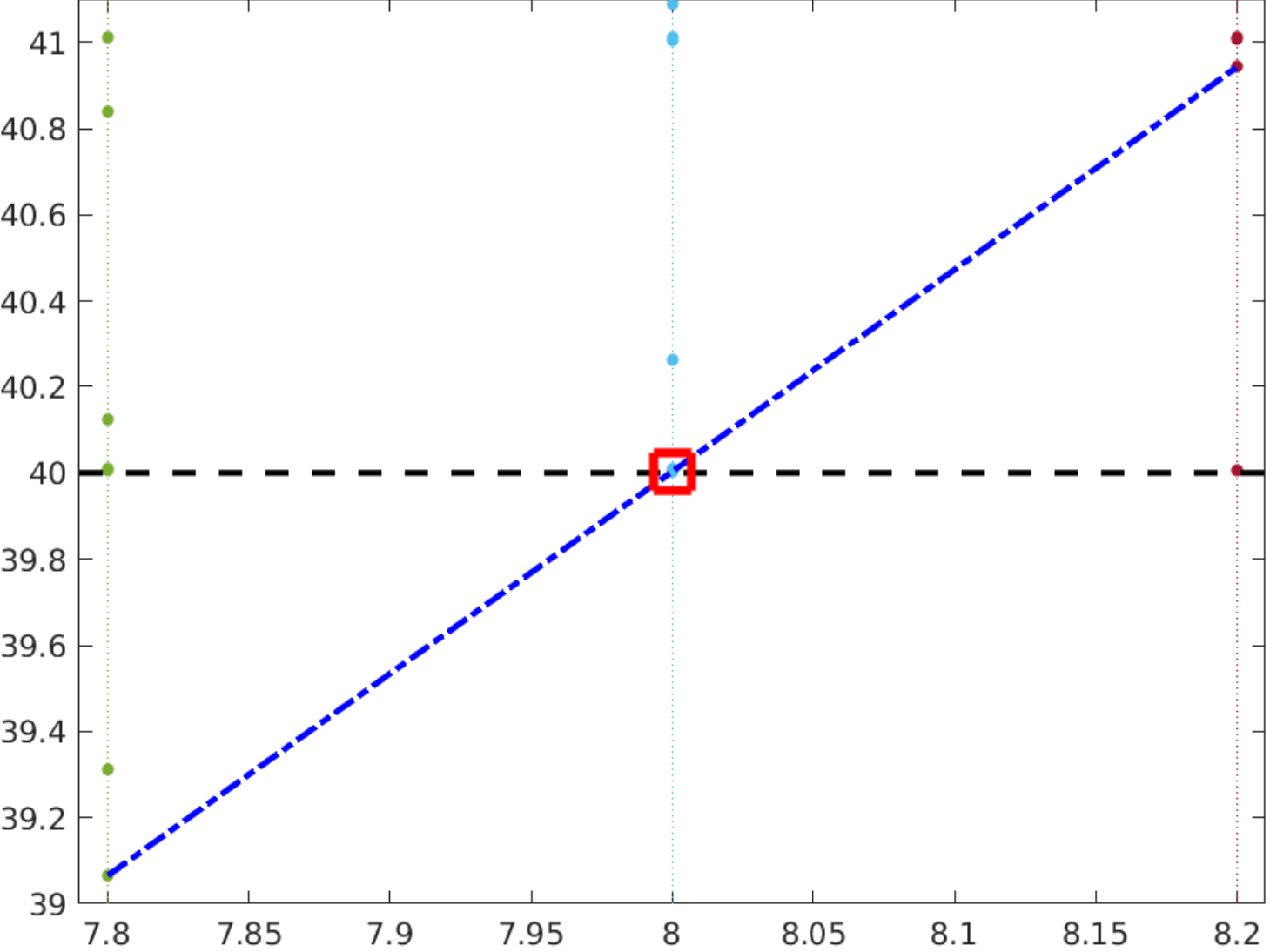}
\includegraphics[width=.45\textwidth]{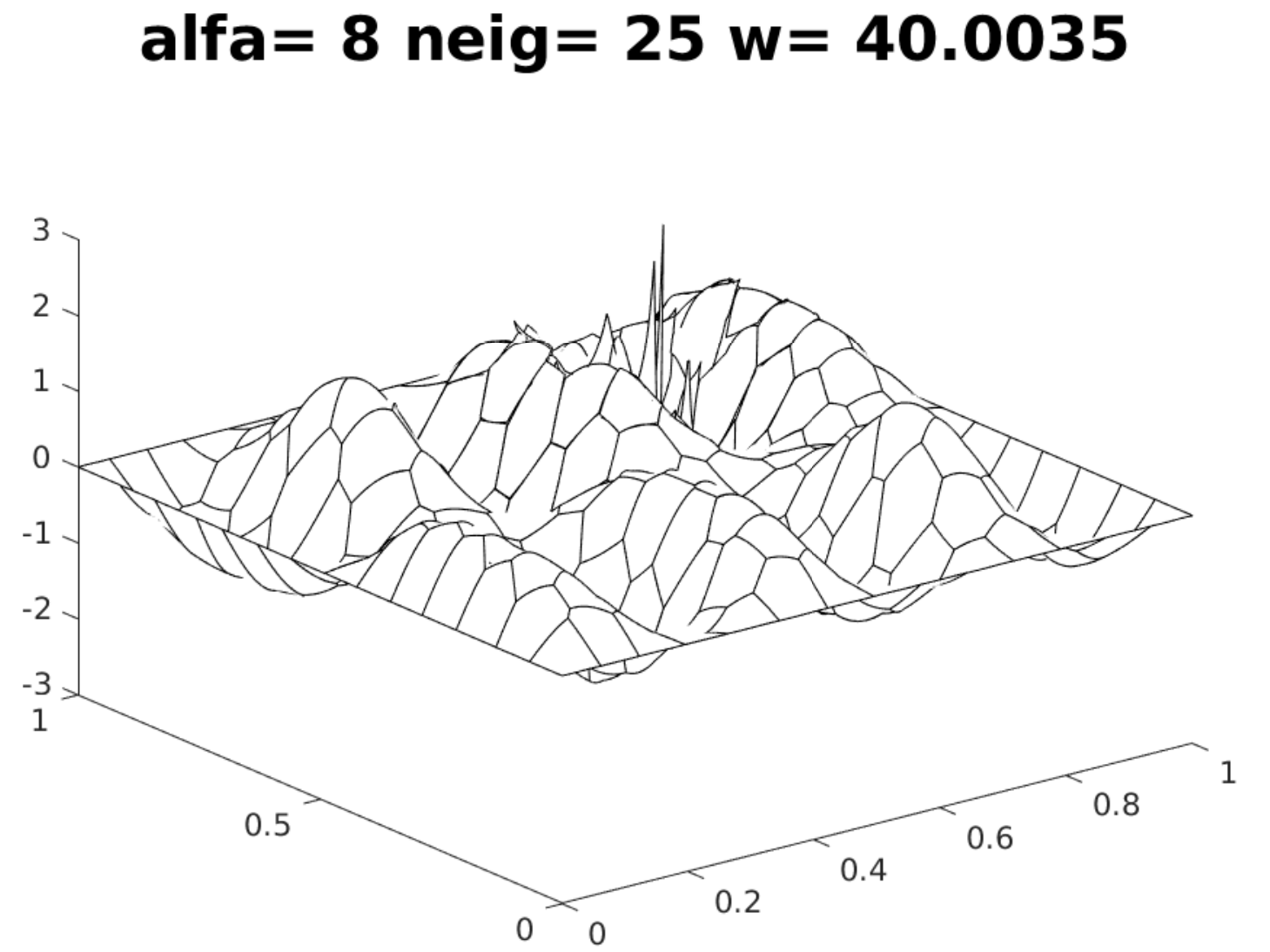}
\end{center}
\caption{Intersections of good and spurious eigenvalues}
\label{fig:autof-marked}
\end{figure}

We conclude this discussion with an example where, for a given value of
$\alpha$, a good eigenvalue (i.e., an eigenvalue corresponding to a correct
approximation) is crossing a spurious one (i.e., an eigenvalue belonging to an
oblique line). In this case it may happen that the two eigenfunctions mix
together, thus yielding to an even more complicated situation. This behavior
is reported in~\cref{fig:autof-marked}, where a region of the plot shown
in~\cref{fig:alfa}(h) is blown-up close to an intersection point:
actually three eigenvalues (a spurious one and two corresponding to good ones)
are clustered at the marked intersection points.

\begin{figure}[h]
\begin{center}
\subfigure[\tiny{$k=1$, $\alpha=0.1$, $\beta\in[0,5]$}]
{\includegraphics[width=.32\textwidth]{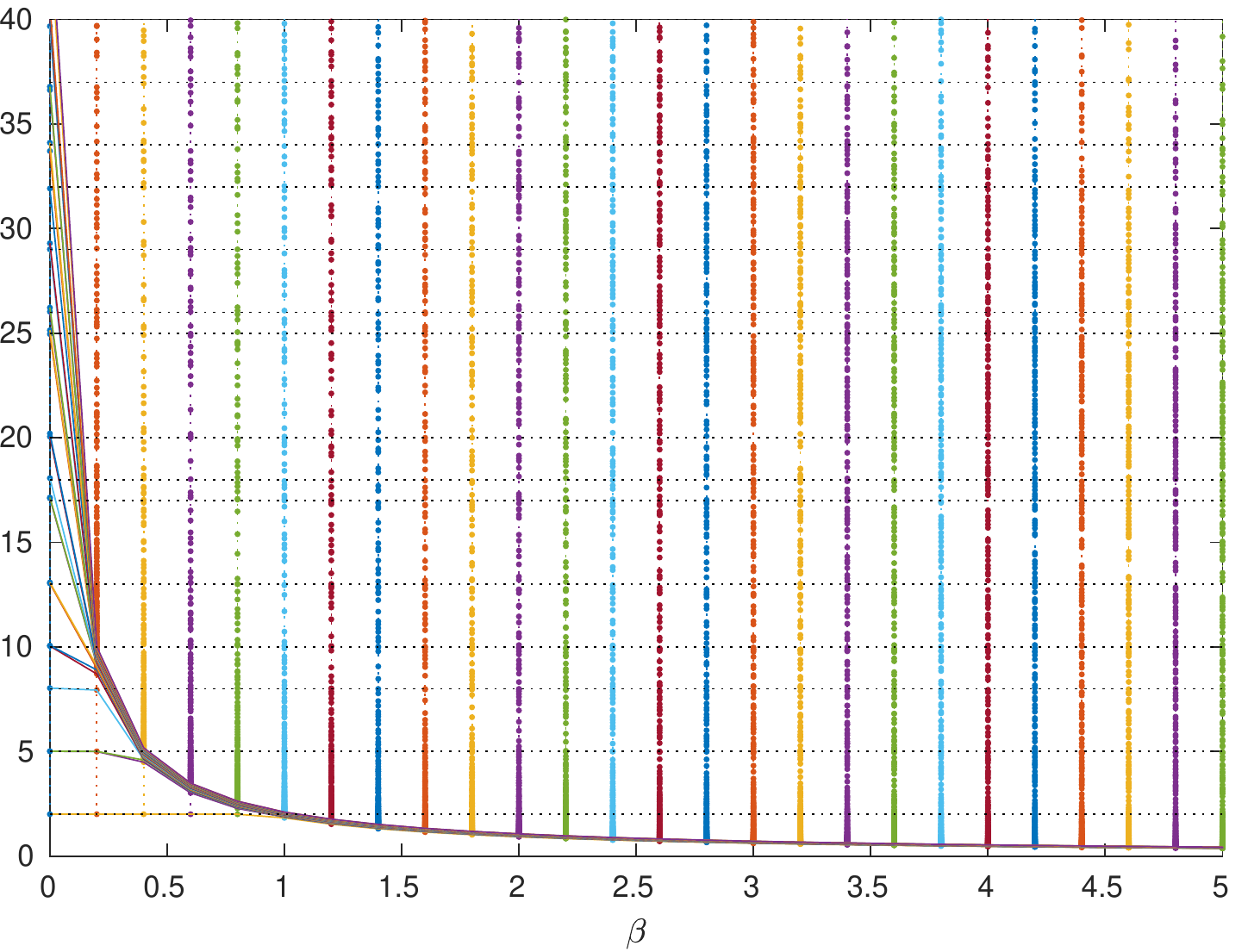}}
\subfigure[\tiny{$k=1$, $\alpha=1$, $\beta\in[0,5]$}]
{\includegraphics[width=.32\textwidth]{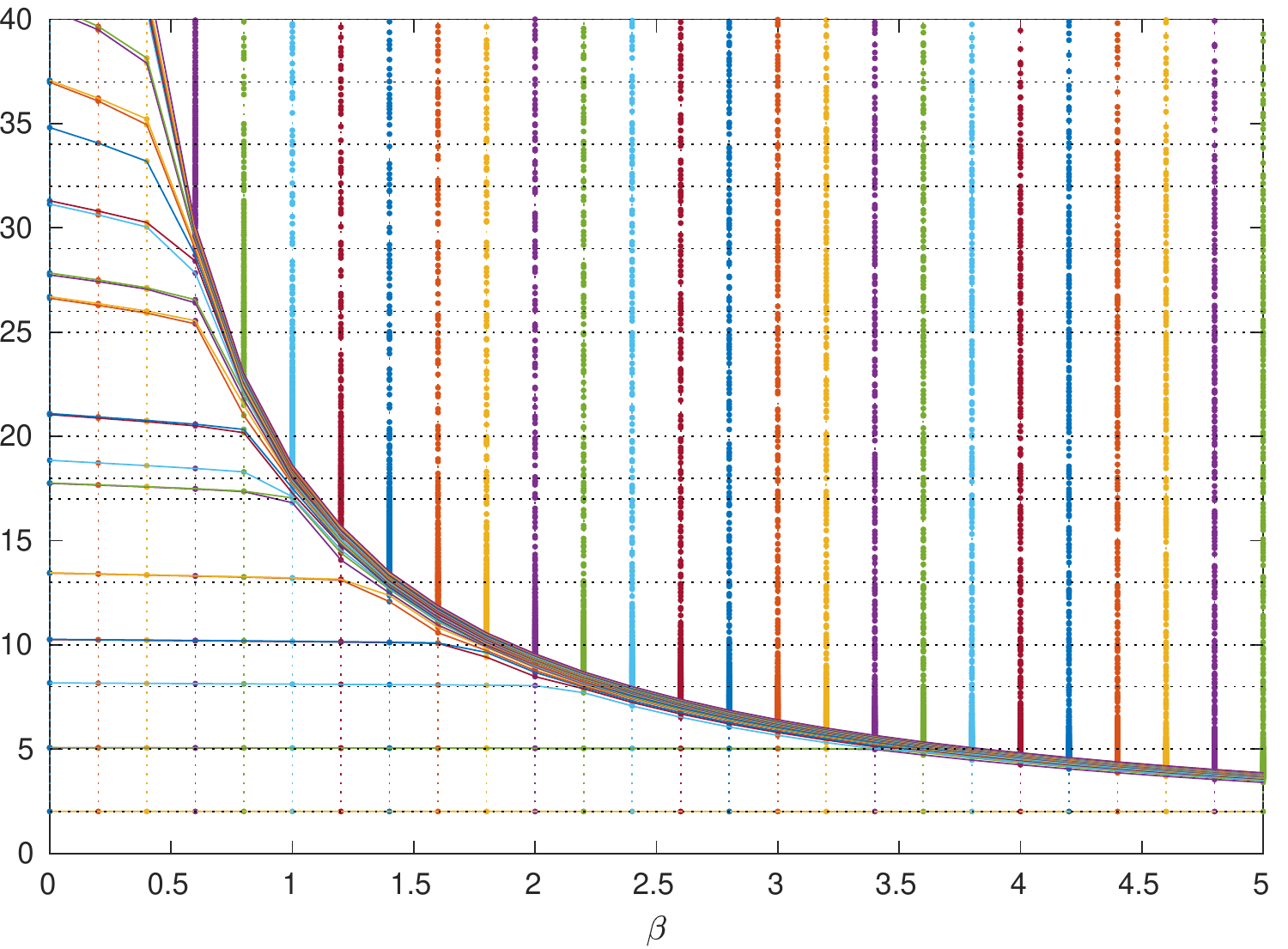}}
\subfigure[\tiny{$k=1$, $\alpha=10$, $\beta\in[0,5]$}]
{\includegraphics[width=.32\textwidth]{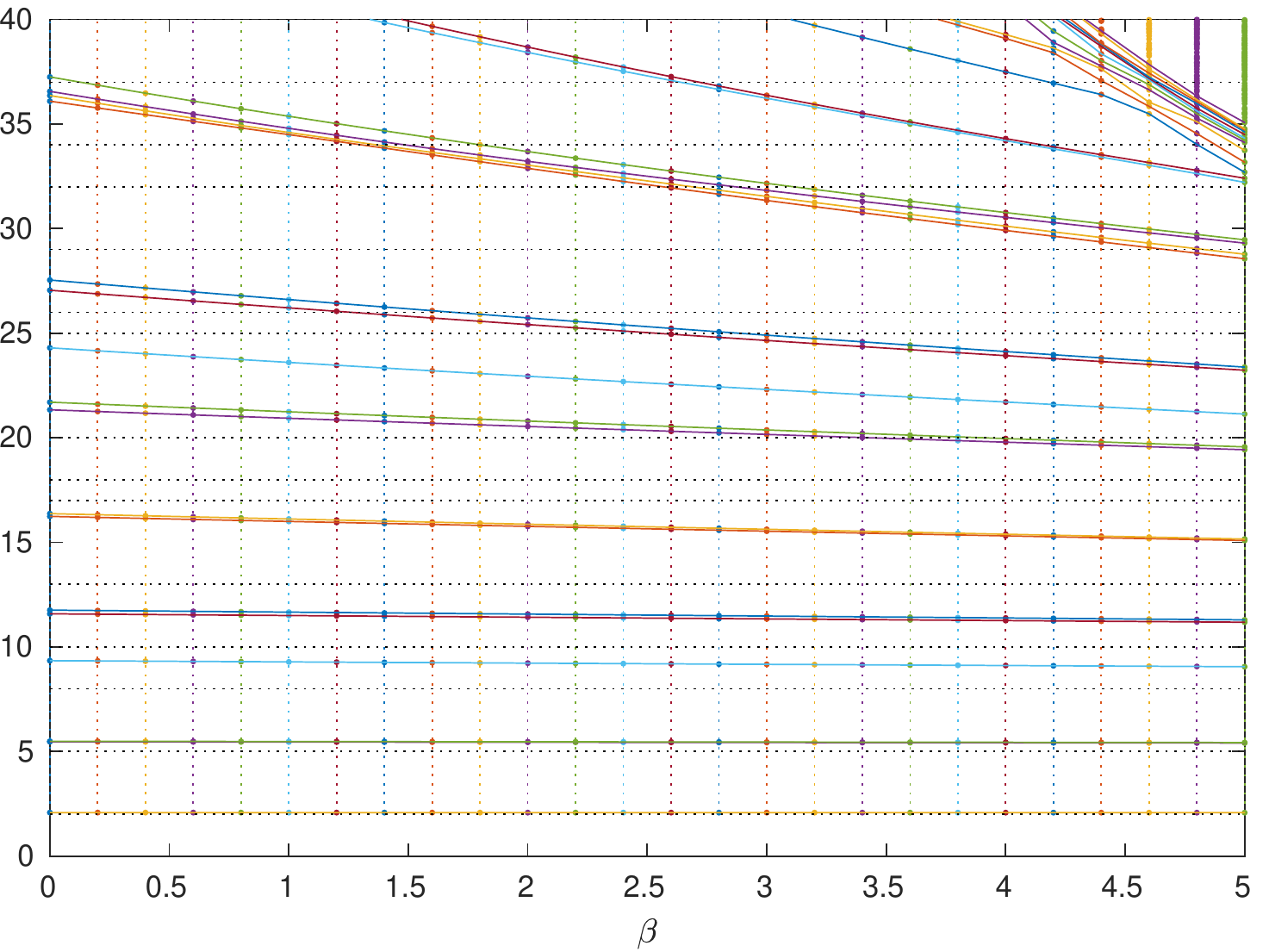}}

\subfigure[\tiny{$k=2$, $\alpha=0.1$, $\beta=[0,5]$}]
{\includegraphics[width=.32\textwidth]{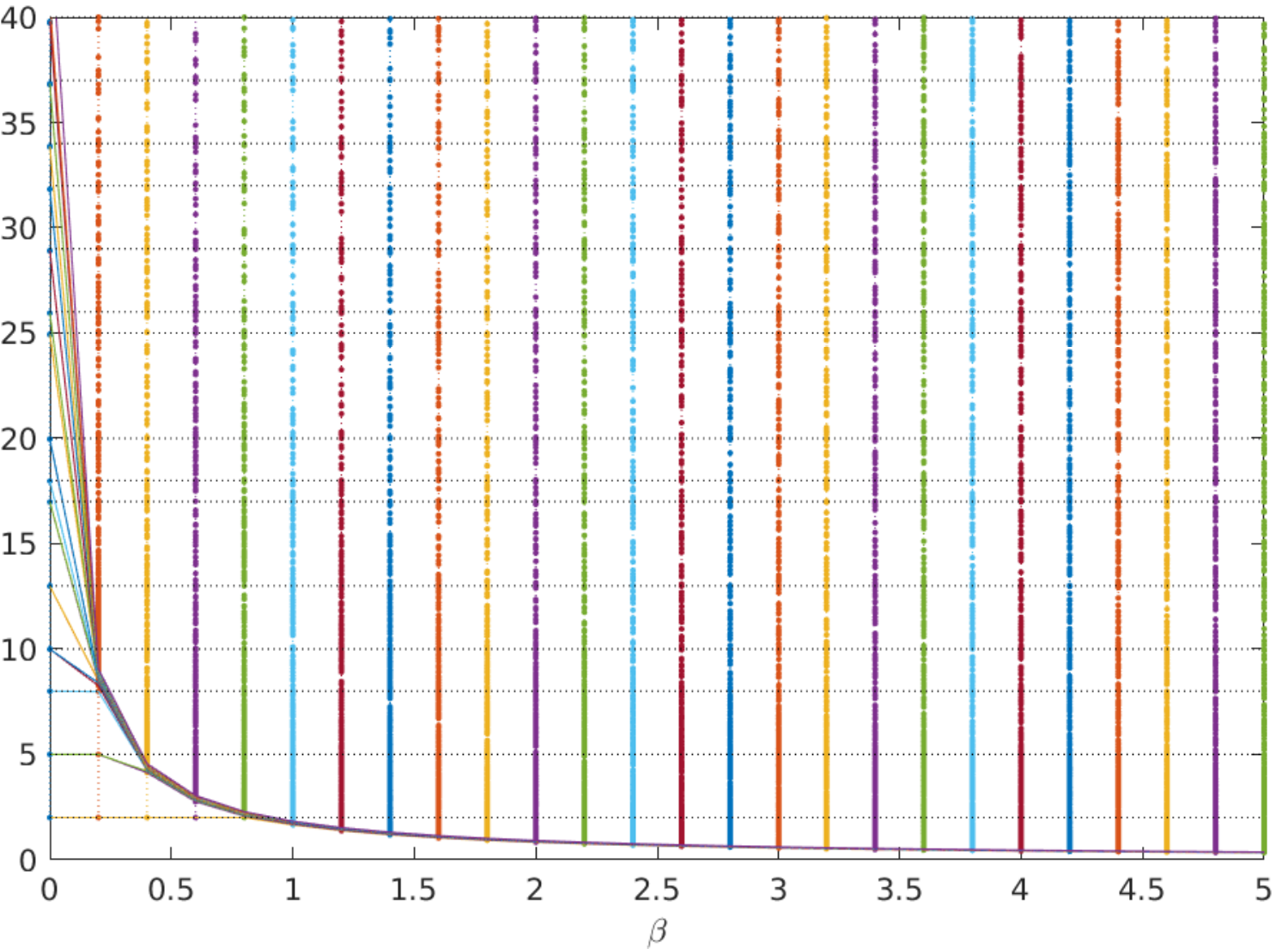}}
\subfigure[\tiny{$k=2$, $\alpha=1$, $\beta=[0,5]$}]
{\includegraphics[width=.32\textwidth]{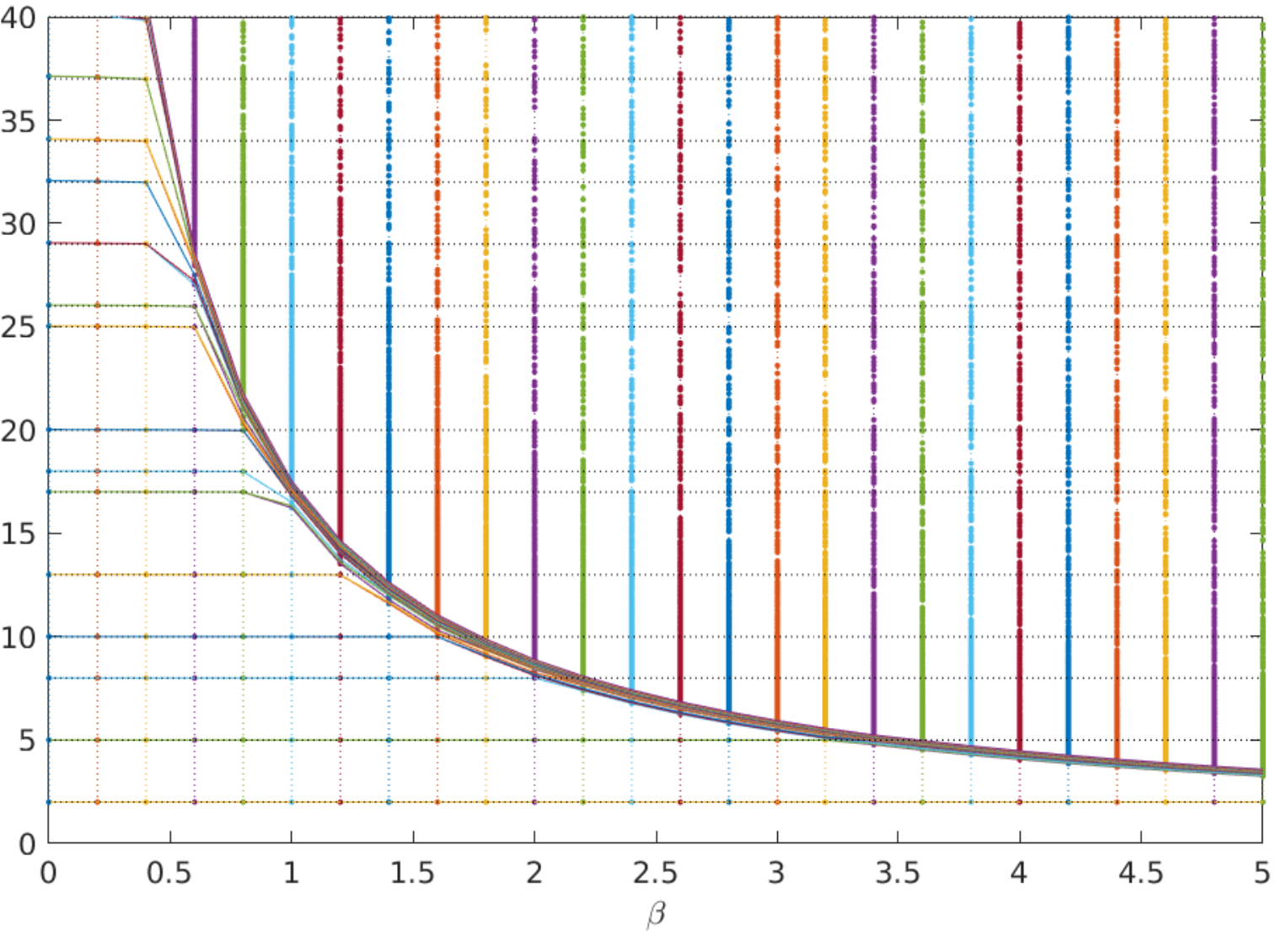}}
\subfigure[\tiny{$k=2$, $\alpha=10$, $\beta=[0,5]$}]
{\includegraphics[width=.32\textwidth]{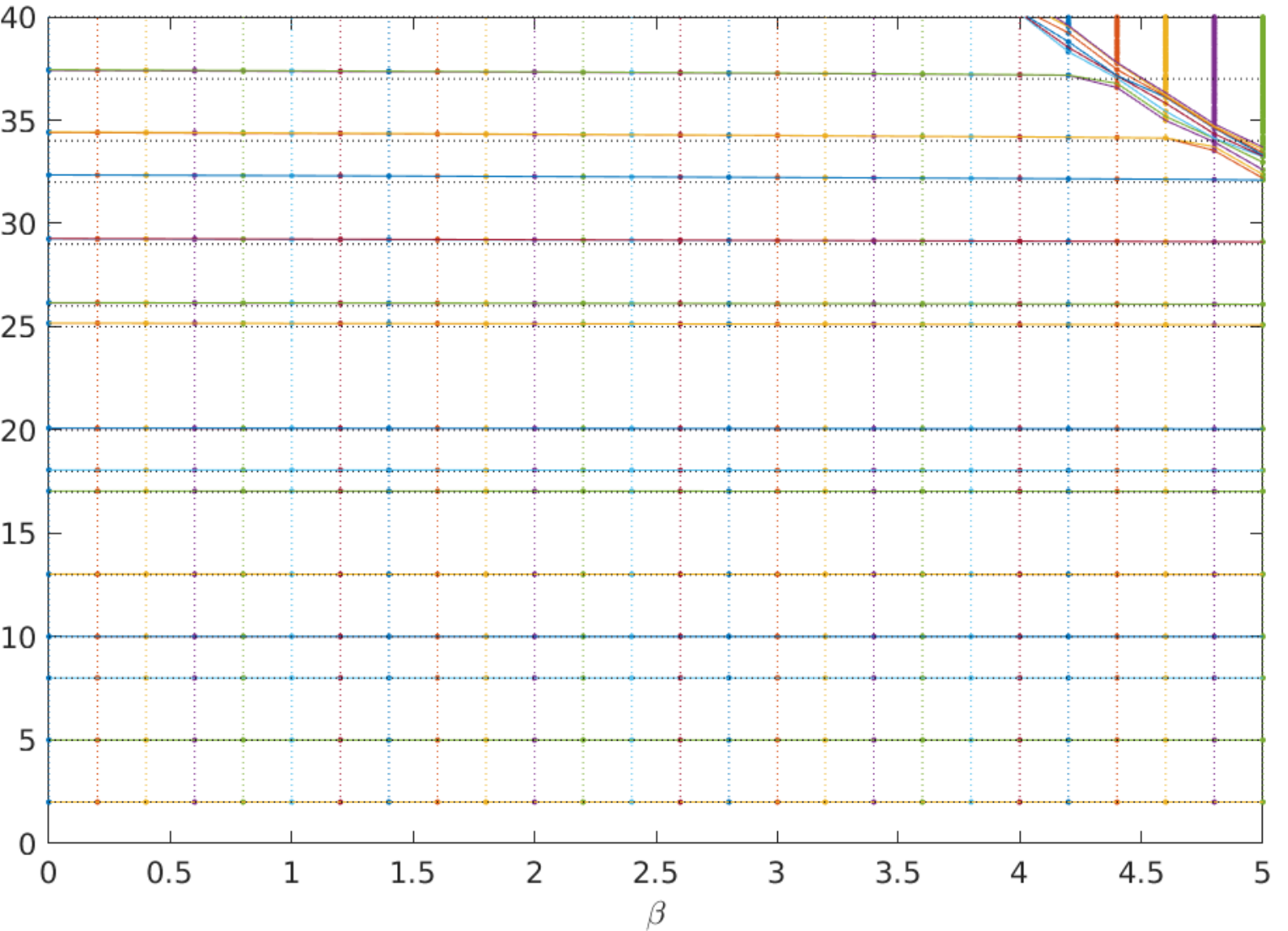}}

\subfigure[\tiny{$k=3$, $\alpha=0.1$, $\beta=[0,5]$}]
{\includegraphics[width=.32\textwidth]{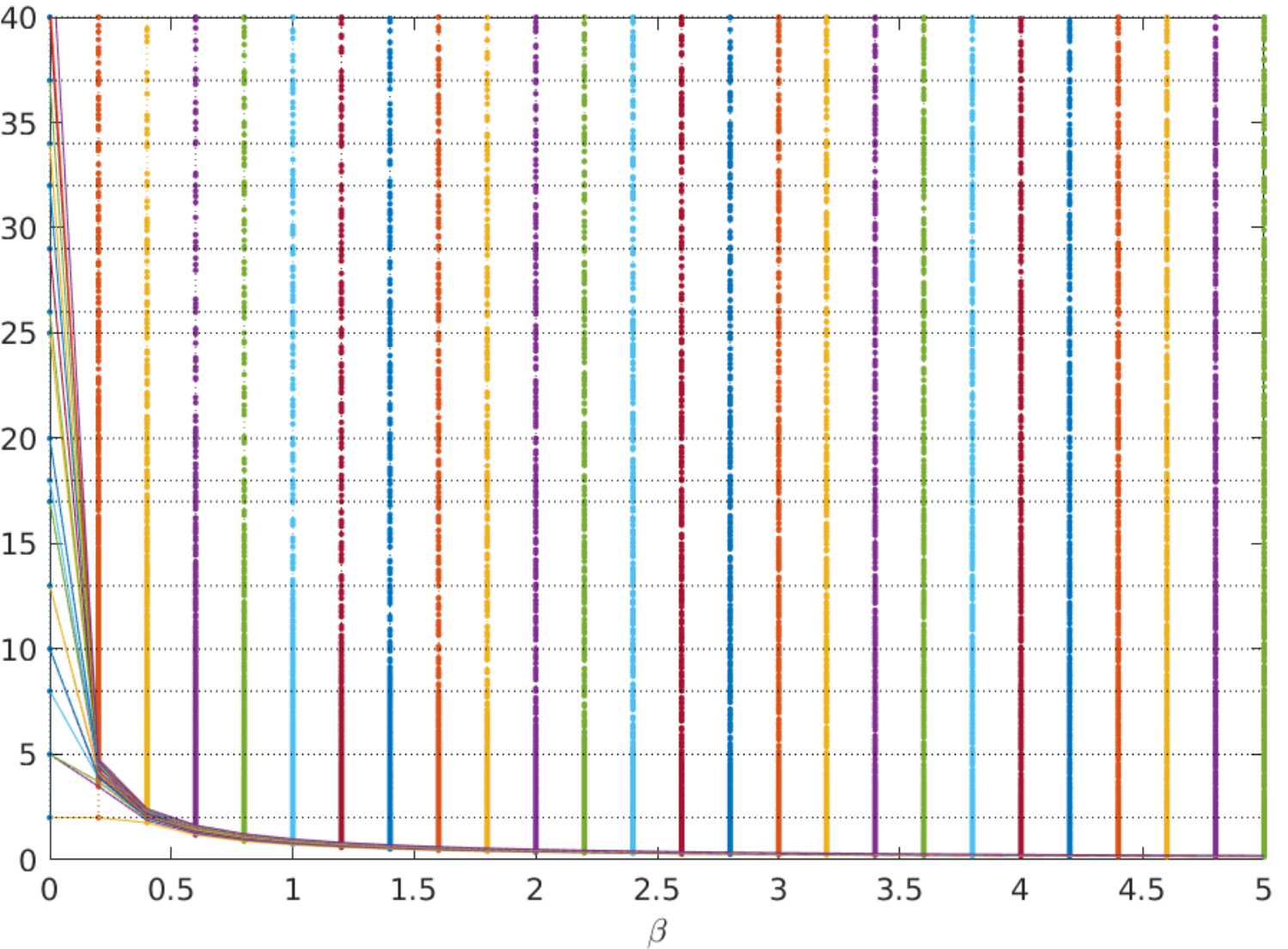}}
\subfigure[\tiny{$k=3$, $\alpha=1$, $\beta=[0,5]$}]
{\includegraphics[width=.32\textwidth]{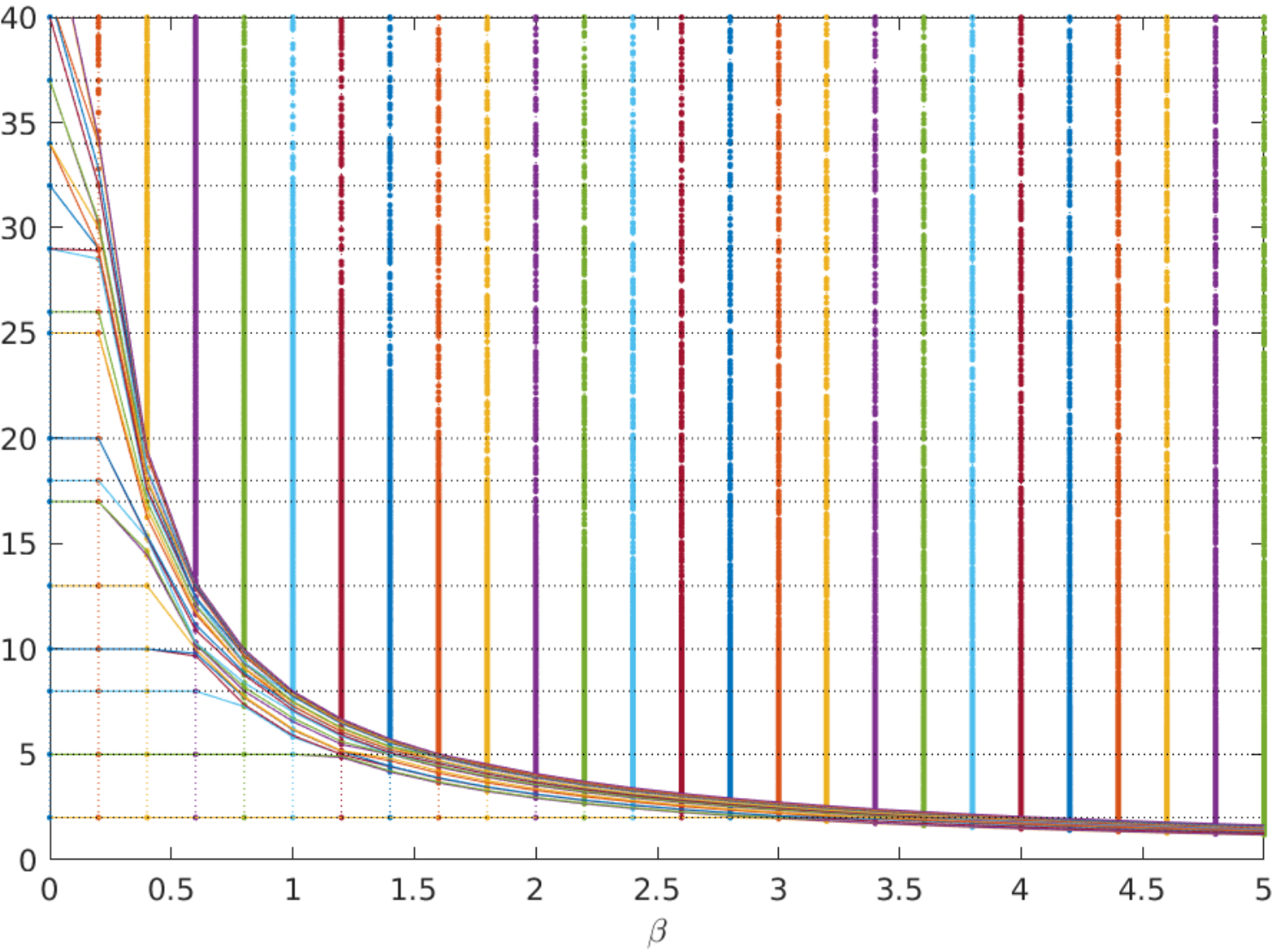}}
\subfigure[\tiny{$k=3$, $\alpha=10$, $\beta=[0,5]$}]
{\includegraphics[width=.32\textwidth]{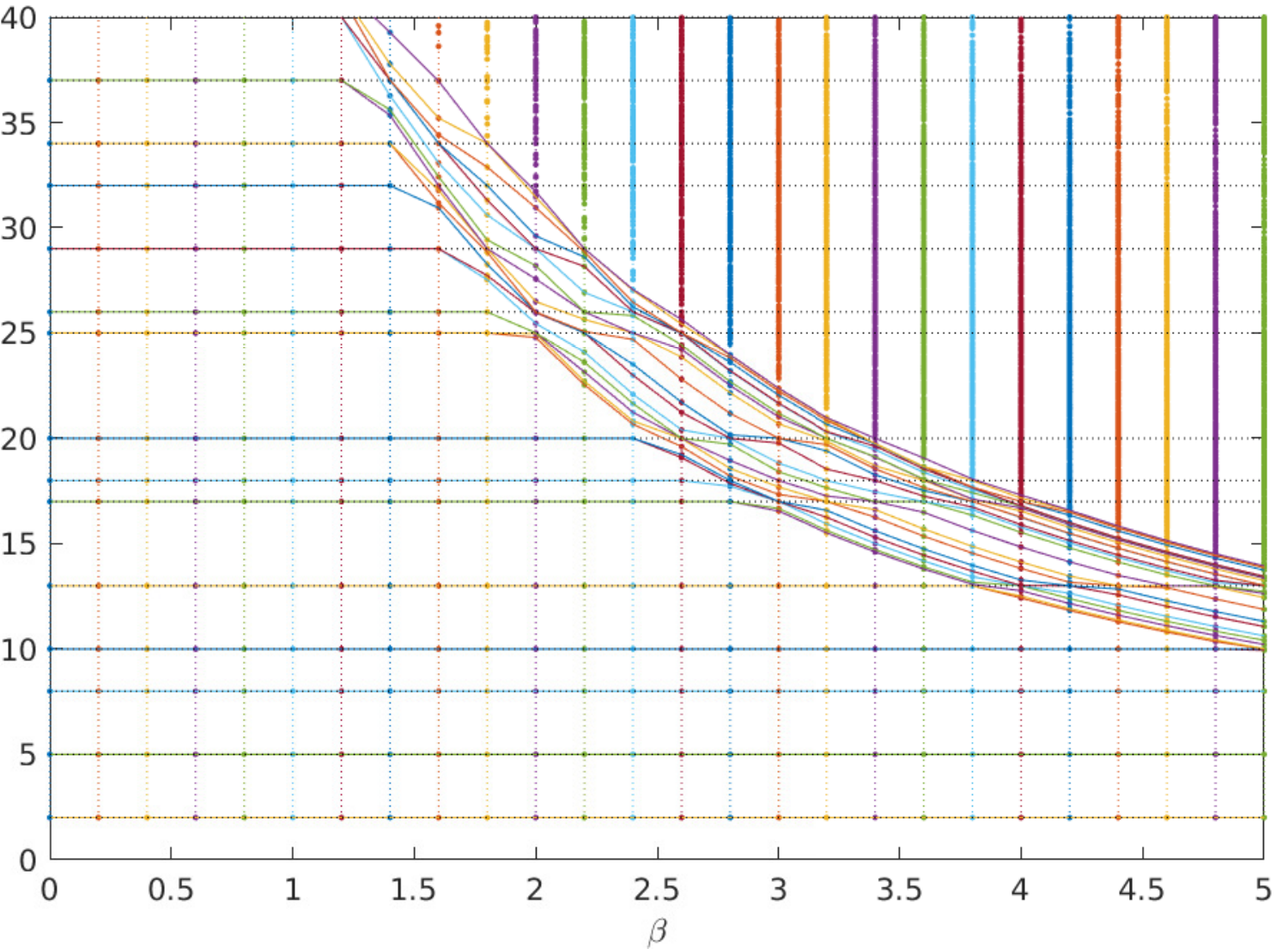}}
\end{center}
\caption{Eigenvalues versus $\beta$ for different values of $k$ and $\alpha$}
\label{fig:beta}
\end{figure}

\Cref{fig:beta} shows the computed eigenvalues smaller that $40$ when
$\beta$ varies from $0$ to $5$ and for a fixed value of $\alpha$.
As in~\cref{fig:girato} and in analogy with~\cref{fig:alfa}, the
rows correspond to the degree $k$ of
polynomials, while the columns refer to different values of $\alpha$. The
dotted horizontal lines represent the exact eigenvalues.
The lines with different colors in each picture follow the $n$-th
eigenvalue for $n=1,\dots,30$. 
It turns out that all lines are originating from curves that look like
hyperbolas when $\beta$ is large. Following each of these hyperbolas from
$\beta=+\infty$ backwards, it happens that when the hyperbola meets a correct
approximation of an eigenvalue of the continuous problem, it deviates from its
trajectory and becomes a (almost horizontal) straight line.
In the case $k=1$, we see that the higher eigenvalues are computed
with decreasing accuracy as $\beta$ approaches $0$.

We recognize in these pictures the situation presented in~\cref{se:caso2}, 
corresponding to the behavior of the eigenvalues when
the parameter $\beta$ in matrix $\m{B}$ varies.
In this test, the kernel of matrix $\Bu$ is not empty only for $k=3$.
Nevertheless, we can see that when $\beta$ approaches $0$, there are several
eigenvalues going to $\infty$.  On the other side, for greater values of
$\beta$ we obtain several spurious eigenvalues.
The range of $\beta$, which gives eigenvalues close to
the exact ones, clearly depends on $k$ and $\alpha$.    

\begin{figure}[h]
\begin{center}
\includegraphics[width=\textwidth]{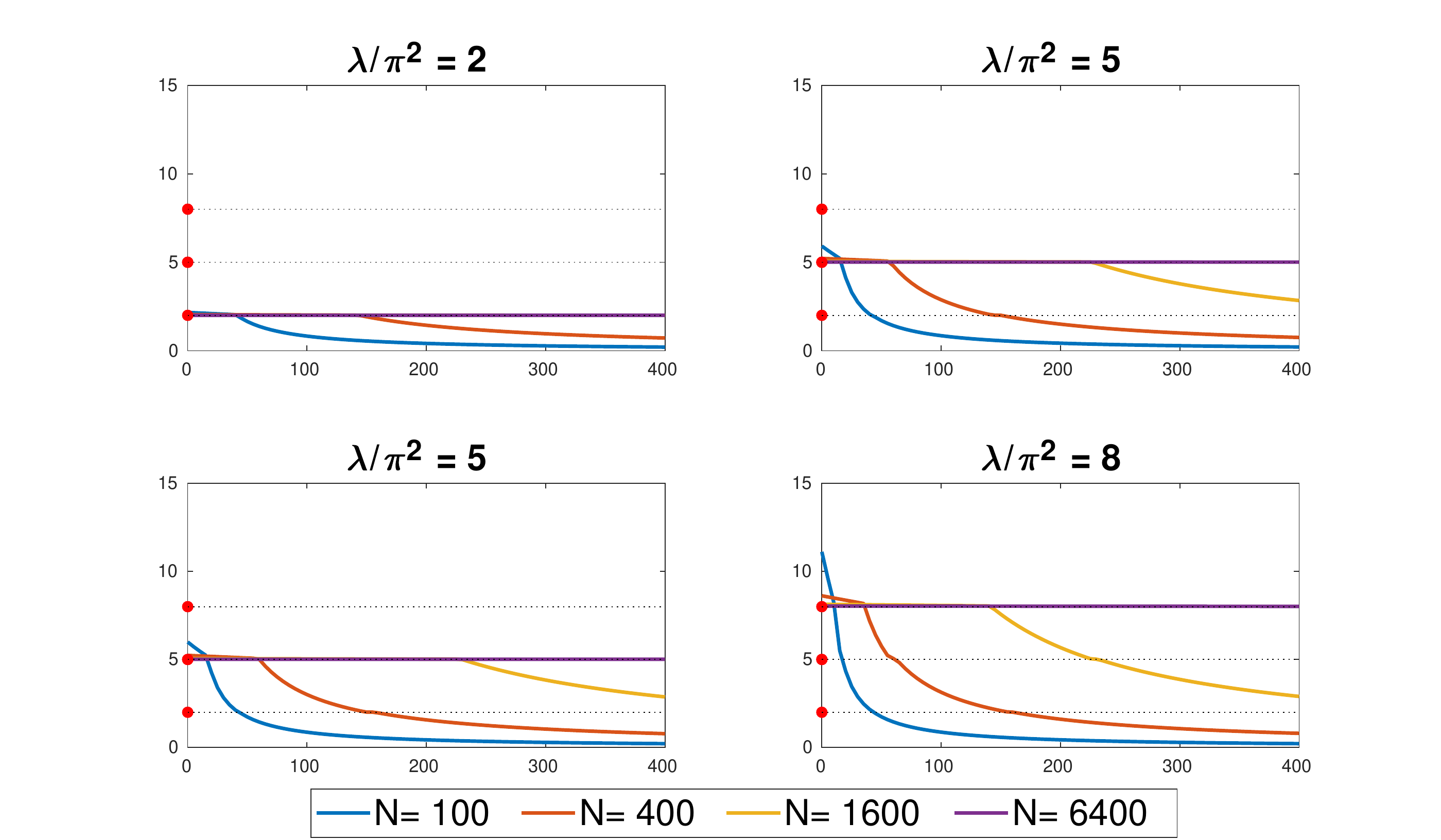}
\end{center}
\caption{First four eigenvalues}
\label{fig:6400}
\end{figure}

\Cref{fig:6400} displays, in separate pictures, the first four eigenvalues, 
with $k=1$, $\alpha=10$, different values of $h$, and $0\le\beta\le400$. Taking
into account that the routine \texttt{eig} sorts the eigenvalues in ascending
order, the four pictures display, in lexicographical order, the first, second,
third and fourth computed eigenvalues. In each subplot, each line refers to a
particular mesh.  We can see that the eigenvalues computed with the finest mesh
seem to be insensitive with respect to the value of $\beta$. 
On the opposite side the coarsest mesh gives approximations 
of the correct values only when $\beta$ is very small and, furthermore,
the accuracy is rather low.
For each eigenvalue and each fixed mesh we recognize a
critical value of the parameter such that greater values of $\beta$ produce
spurious eigenvalues.
The behavior of these eigenvalues clearly reproduces that of the eigenvalues 
in~\cref{fig:case1} (right) referring to Case 2. The results are plotted with
a different perspective depending on the fact that the results now depend also
on the computational mesh. 
The right bottom plot of~\cref{fig:6400} highlights a phenomenon which
already appears in~\cref{fig:beta}(i). Indeed, we see that the red line
corresponding to the fourth computed eigenvalue for $N=400$ lies along an
hyperbola until $\beta=65$ where it reaches the value $5$ associated with
second and third exact eigenvalues. Between $\beta=65$ and $\beta=55$ the red
line remains close to $5$, then decreasing $\beta$ it follows a different
hyperbola until it reaches the expected value for $\beta=35$. 

\section*{Conclusions}
In this paper we have discussed how numerically computed eigenvalues can
depend on discretization parameters. \Cref{se:param}
shows the dependence on $\alpha$ and $\beta$ of the eigenvalues
of~\eqref{eq:eig} when $\m{A}$ and $\m{B}$ have the forms~\eqref{eq:A}
and~\eqref{eq:B}, respectively. In~\cref{se:VEM} we have studied the
behavior of the eigenvalues of the Laplace operator computed with the Virtual
Element Method. The presence of two parameters resembles the abstract setting
of~\cref{se:param}; even if assumptions satisfied by the VEM matrices are
more complicated than the ones previously discussed, the numerical results are
pretty much in agreement. The present work opens the question of a viable
choice of the parameters for eigenvalue computations when the discretization
scheme depends on a suitable tuning of them (such as in the case of VEM).

 \section*{Acknowledgments}
The authors are members of INdAM Research group GNCS and their
research is supported by PRIN/MIUR. The research of the first and third
authors is partially supported by IMATI/CNR.

\bibliography{ref}

\begin{thebibliography}{10}

\bibitem{AABMR}
{\sc B.~Ahmad, A.~Alsaedi, F.~Brezzi, L.~D. Marini, and A.~Russo}, {\em
  Equivalent projectors for virtual element methods}, Comput. Math. Appl., 66
  (2013), pp.~376--391, \url{https://doi.org/10.1016/j.camwa.2013.05.015},
  \url{https://doi.org/10.1016/j.camwa.2013.05.015}.

\bibitem{BaCo}
{\sc S.~Badia and R.~Codina}, {\em A nodal-based finite element approximation
  of the {M}axwell problem suitable for singular solutions}, SIAM J. Numer.
  Anal., 50 (2012), pp.~398--417, \url{https://doi.org/10.1137/110835360}.

\bibitem{BBCMMR}
{\sc L.~Beir\~{a}o~da Veiga, F.~Brezzi, A.~Cangiani, G.~Manzini, L.~D. Marini,
  and A.~Russo}, {\em Basic principles of virtual element methods}, Math.
  Models Methods Appl. Sci., 23 (2013), pp.~199--214,
  \url{https://doi.org/10.1142/S0218202512500492}.

\bibitem{BLR}
{\sc L.~Beir\~{a}o~da Veiga, C.~Lovadina, and A.~Russo}, {\em Stability
  analysis for the virtual element method}, Math. Models Methods Appl. Sci., 27
  (2017), pp.~2557--2594, \url{https://doi.org/10.1142/S021820251750052X}.

\bibitem{BMRR}
{\sc L.~Beir\~{a}o~da Veiga, D.~Mora, G.~Rivera, and R.~Rodr\'{\i}guez}, {\em A
  virtual element method for the acoustic vibration problem}, Numer. Math., 136
  (2017), pp.~725--763, \url{https://doi.org/10.1007/s00211-016-0855-5}.

\bibitem{bfg}
{\sc D.~Boffi, M.~Farina, and L.~Gastaldi}, {\em On the approximation of
  {M}axwell's eigenproblem in general 2{D} domains}, Computers \& Structures,
  79 (2001), pp.~1089 -- 1096.

\bibitem{BoGue}
{\sc A.~Bonito and J.-L. Guermond}, {\em Approximation of the eigenvalue
  problem for the time harmonic {M}axwell system by continuous {L}agrange
  finite elements}, Math. Comp., 80 (2011), pp.~1887--1910,
  \url{https://doi.org/10.1090/S0025-5718-2011-02464-6}.

\bibitem{2006}
{\sc A.~Buffa and I.~Perugia}, {\em Discontinuous {G}alerkin approximation of
  the {M}axwell eigenproblem}, SIAM J. Numer. Anal., 44 (2006), pp.~2198--2226,
  \url{https://doi.org/10.1137/050636887}.

\bibitem{CMRS}
{\sc A.~Cangiani, G.~Manzini, A.~Russo, and N.~Sukumar}, {\em Hourglass
  stabilization and the virtual element method}, Internat. J. Numer. Methods
  Engrg., 102 (2015), pp.~404--436, \url{https://doi.org/10.1002/nme.4854}.

\bibitem{CoDaMax}
{\sc M.~Costabel and M.~Dauge}, {\em {M}axwell and {L}am\'{e} eigenvalues on
  polyhedra}, Math. Methods Appl. Sci., 22 (1999), pp.~243--258.

\bibitem{CoDareg}
{\sc M.~Costabel and M.~Dauge}, {\em Weighted regularization of {M}axwell
  equations in polyhedral domains. {A} rehabilitation of nodal finite
  elements}, Numer. Math., 93 (2002), pp.~239--277,
  \url{https://doi.org/10.1007/s002110100388}.

\bibitem{CoDaDurham}
{\sc M.~Costabel and M.~Dauge}, {\em Computation of resonance frequencies for
  {M}axwell equations in non-smooth domains}, in Topics in computational wave
  propagation, vol.~31 of Lect. Notes Comput. Sci. Eng., Springer, Berlin,
  2003, pp.~125--161, \url{https://doi.org/10.1007/978-3-642-55483-4_4}.

\bibitem{ElsnerSun}
{\sc L.~Elsner and J.~G. Sun}, {\em Perturbation theorems for the generalized
  eigenvalue problem}, Linear Algebra Appl., 48 (1982), pp.~341--357,
  \url{https://doi.org/10.1016/0024-3795(82)90120-3}.

\bibitem{GMV}
{\sc F.~Gardini, G.~Manzini, and G.~Vacca}, {\em The nonconforming virtual
  element method for eigenvalue problems}, ESAIM Math. Model. Numer. Anal., 53
  (2019), pp.~749--774, \url{https://doi.org/10.1051/m2an/2018074}.

\bibitem{GV}
{\sc F.~Gardini and G.~Vacca}, {\em Virtual element method for second-order
  elliptic eigenvalue problems}, IMA J. Numer. Anal., 38 (2018),
  pp.~2026--2054, \url{https://doi.org/10.1093/imanum/drx063}.

\bibitem{GVL}
{\sc G.~H. Golub and C.~F. Van~Loan}, {\em Matrix computations}, Johns Hopkins
  Studies in the Mathematical Sciences, Johns Hopkins University Press,
  Baltimore, MD, fourth~ed., 2013.

\bibitem{greenbaum2019firstorder}
{\sc A.~Greenbaum, R.~cang Li, and M.~L. Overton}, {\em First-order
  perturbation theory for eigenvalues and eigenvectors}, 2019,
  \url{https://arxiv.org/abs/1903.00785}.

\bibitem{LiStewart}
{\sc R.-C. Li and G.~W. Stewart}, {\em A new relative perturbation theorem for
  singular subspaces}, Linear Algebra Appl., 313 (2000), pp.~41--51,
  \url{https://doi.org/10.1016/S0024-3795(00)00074-4}.

\bibitem{MRR}
{\sc D.~Mora, G.~Rivera, and R.~Rodr\'{\i}guez}, {\em A virtual element method
  for the {S}teklov eigenvalue problem}, Math. Models Methods Appl. Sci., 25
  (2015), pp.~1421--1445, \url{https://doi.org/10.1142/S0218202515500372}.

\bibitem{MRV}
{\sc D.~Mora, G.~Rivera, and I.~Vel\'{a}squez}, {\em A virtual element method
  for the vibration problem of {K}irchhoff plates}, ESAIM Math. Model. Numer.
  Anal., 52 (2018), pp.~1437--1456, \url{https://doi.org/10.1051/m2an/2017041}.

\bibitem{MV}
{\sc D.~Mora and I.~Vel\'{a}squez}, {\em A virtual element method for the
  transmission eigenvalue problem}, Math. Models Methods Appl. Sci., 28 (2018),
  pp.~2803--2831, \url{https://doi.org/10.1142/S0218202518500616}.

\bibitem{CGMMV}
{\sc O.{\v C}ert{\'i}k, F.~Gardini, G.~Manzini, L.~Mascotto, and G.~Vacca},
  {\em The p- and hp-versions of the virtual element method for elliptic
  eigenvalue problems}, Computers \& Mathematics with Applications,  (2019),
  \url{https://doi.org/https://doi.org/10.1016/j.camwa.2019.10.018}.

\bibitem{2010}
{\sc D.~S\'{a}rm\'{a}ny, F.~Izs\'{a}k, and J.~J.~W. van~der Vegt}, {\em Optimal
  penalty parameters for symmetric discontinuous {G}alerkin discretisations of
  the time-harmonic {M}axwell equations}, J. Sci. Comput., 44 (2010),
  pp.~219--254, \url{https://doi.org/10.1007/s10915-010-9366-1}.

\bibitem{MR1066108}
{\sc B.~Simon}, {\em Fifty years of eigenvalue perturbation theory}, Bull.
  Amer. Math. Soc. (N.S.), 24 (1991), pp.~303--319,
  \url{https://doi.org/10.1090/S0273-0979-1991-16020-9}.

\bibitem{StewartSun}
{\sc G.~W. Stewart and J.~G. Sun}, {\em Matrix perturbation theory}, Computer
  Science and Scientific Computing, Academic Press, Inc., Boston, MA, 1990.

\bibitem{WarburtonEmbree}
{\sc T.~Warburton and M.~Embree}, {\em The role of the penalty in the local
  discontinuous {G}alerkin method for {M}axwell's eigenvalue problem}, Comput.
  Methods Appl. Mech. Engrg., 195 (2006), pp.~3205--3223,
  \url{https://doi.org/10.1016/j.cma.2005.06.011}.

\end{thebibliography}
\bibliographystyle{siamplain}
\end{document}